\documentclass[10pt]{article}
\usepackage[T1]{fontenc}
\usepackage[utf8]{inputenc}
\usepackage{lmodern}
\usepackage[a4paper,margin=1in]{geometry}
\usepackage{amsmath,amssymb,amsthm,mathtools,mathrsfs}
\usepackage{microtype}
\usepackage[skip=0.45\baselineskip plus 2pt,indent=0pt]{parskip}
\usepackage{enumitem,booktabs,array,longtable}
\usepackage{xcolor,xurl}
\usepackage{tikz}
\usetikzlibrary{arrows.meta,positioning}
\usepackage{hyperref}
\usepackage{aliascnt}
\usepackage[capitalise,nameinlink,noabbrev]{cleveref}
\setlist[enumerate,1]{label=\textnormal{(\roman*)},leftmargin=2.35em}
\setlist[itemize]{leftmargin=1.8em}
\allowdisplaybreaks
\numberwithin{equation}{section}
\newtheorem{theorem}{Theorem}[section]
\newaliascnt{lemma}{theorem}
\newtheorem{lemma}[lemma]{Lemma}\aliascntresetthe{lemma}
\newaliascnt{proposition}{theorem}
\newtheorem{proposition}[proposition]{Proposition}\aliascntresetthe{proposition}
\newaliascnt{corollary}{theorem}
\newtheorem{corollary}[corollary]{Corollary}\aliascntresetthe{corollary}
\newaliascnt{claim}{theorem}
\newtheorem{claim}[claim]{Claim}\aliascntresetthe{claim}
\crefname{claim}{Claim}{Claims}
\Crefname{claim}{Claim}{Claims}
\theoremstyle{definition}
\newaliascnt{definition}{theorem}
\aliascntresetthe{definition}
\newaliascnt{remark}{theorem}
\aliascntresetthe{remark}

\newcommand{\Q}{\mathbb Q}

\newcommand{\cA}{\mathcal A}\newcommand{\cB}{\mathcal B}
\newcommand{\cD}{\mathcal D}\newcommand{\cE}{\mathcal E}\newcommand{\cF}{\mathcal F}
\newcommand{\cG}{\mathcal G}\newcommand{\cH}{\mathcal H}

\newcommand{\cQ}{\mathcal Q}\newcommand{\cR}{\mathcal R}\newcommand{\cS}{\mathcal S}

\newcommand{\Forb}{\operatorname{Forb}}\newcommand{\Ext}{\operatorname{Ext}}
\newcommand{\ex}{\operatorname{ex}}
\newcommand{\Def}{\operatorname{Def}}
\newcommand{\Sep}{\operatorname{Sep}}
\newcommand{\BB}{\operatorname{BB}}
\newcommand{\Improve}{\operatorname{Improve}}
\newcommand{\eps}{\varepsilon}\newcommand{\dotcup}{\mathbin{\dot\cup}}
\newcommand{\x}{\boldsymbol{x}}\newcommand{\y}{\boldsymbol{y}}
\newcommand{\zero}{\mathbf 0}
\newcommand{\Cand}{\operatorname{Cand}}

\hypersetup{
colorlinks=true, linkcolor=blue!55!black, citecolor=green!40!black, urlcolor=blue!60!black, pdftitle={Turing Universality, Computability, and Incompleteness in Hypergraph Tur\'an Theory}, pdfauthor={Heng Li and Xizhi Liu}, pdfsubject={Computability, structural universality, and Goedel--Rosser incompleteness in ordinary hypergraph Turan problems}, pdfkeywords={hypergraph Turan problem, Turing machine, computability, Goedel incompleteness, extremal limit, stability, phase transition} }
\title{\fontsize{16}{20}\selectfont\bfseries \makebox[\textwidth][c]{Turing universality, computability, and incompleteness}\\ in hypergraph Tur\'an theory}
\author{Heng~Li\and Xizhi~Liu}
\date{\today}
\begin{document}
\maketitle
\begin{abstract}
The hypergraph Tur\'an problem asks for the largest asymptotic edge density of hypergraphs avoiding a prescribed family and for the structure of nearly extremal hypergraphs. We show that both questions can encode arbitrary computation. Fix a universal Turing machine $\mathsf U$. For every fixed $r\ge3$, there is a rational threshold $\tau_r\in(0,1)$ such that each binary word $\beta$ can be effectively mapped to a finite family $\cF_{r,\beta}$ of $r$-graphs whose Tur\'an density equals $\tau_r$ when $\mathsf U$ does not halt on $\beta$ and exceeds $\tau_r$ when it halts. For $r=3$, the family consists of a single forbidden $3$-graph $P_\beta$.

For every sufficiently large fixed $r$, we also construct finite forbidden families whose extremal structures reflect the same halting distinction. In the nonhalting case, there is a unique extremal limit, and every nearly extremal hypergraph is close in edit distance to a fixed finite-profile construction. In the halting case, the extremal limit space splits into 2 nonempty compact phases, separated by the opposite signs of a continuous statistic depending only on $r$. Each phase is homeomorphic to the entire extremal limit space of $P_\beta$. Thus, uniqueness, connectedness, and stability with respect to a single finite profile are undecidable.

The reductions can be verified in ZFC. Consequently, for every consistent computably axiomatized theory extending ZFC and every fixed $r\ge3$, we construct a finite family $\cF$ for which the true equality $\pi(\cF)=\tau_r$ is neither provable nor refutable in that theory. A single forbidden $3$-graph suffices when $r=3$. Analogous independence statements hold for the structural properties above in sufficiently large uniformity.

Although finite-family Tur\'an densities can be approximated effectively, exact comparison is undecidable even for single forbidden $3$-graphs. We determine its complexity for finite and effectively enumerated forbidden families, and characterize the densities attainable by the latter. Quantitatively, the smallest finite witnesses of a density improvement have Busy-Beaver growth, and positive density gaps admit no computable lower bound in terms of the encoded input.
\end{abstract}

\tableofcontents

\section{Introduction}\label{sec:intro}
An $r$-graph is a finite family of $r$-subsets of a finite vertex set. For a family $\cF$ of forbidden $r$-graphs, let $\ex(n,\cF)$ be the largest size of an $n$-vertex $\cF$-free $r$-graph. Its Tur\'an density is
\[
 \pi(\cF)\coloneqq\lim_{n\to\infty}\frac{\ex(n,\cF)}{\binom nr}.
\]
The limit exists by averaging~\cite{KatonaNemetzSimonovits1964}. Forbidden hypergraphs are assumed to have edges.

For graphs, the Erd\H{o}s--Stone--Simonovits theorem~\cite{ErdosStone1946,ErdosSimonovits1966} determines the density from the smallest chromatic number among the forbidden graphs. Hypergraph densities have a richer arithmetic structure. Baber and Talbot~\cite{BaberTalbot2012} and Pikhurko~\cite{Pikhurko2014} independently constructed finite families with irrational density. Liu and Pikhurko~\cite{LiuPikhurko2023} proved that the algebraic degrees of finite-family densities are unbounded, and Li, Liu, and Liu~\cite{LiLiuLiu2026} constructed finite families with transcendental density.

The structural problem asks which configurations can attain or nearly attain the extremal density. For graph Tur\'an problems of positive density, classical stability~\cite{Erdos1967Stability,Simonovits1968Stability} gives complete multipartite approximations to near-extremal graphs. Hypergraph Tur\'an problems can have several extremal structures separated in edit distance. Liu and Mubayi~\cite{LiuMubayi2022}, Liu, Mubayi, and Reiher~\cite{LiuMubayiReiher2026}, and Hou et al.~\cite{HouEtAl2023} constructed examples with two, arbitrarily many, and infinitely many extremal constructions, respectively. Liu and Pikhurko~\cite{LiuPikhurko2025} obtained further examples by mixing patterns.

For each fixed $n$, exhaustive search computes $\ex(n,\cF)$, but the density and structural questions concern hypergraphs of arbitrarily large order. From an algorithmic perspective, it is useful to distinguish three questions.
\begin{enumerate}
\item Can $\pi(\cF)$ be approximated to any prescribed accuracy?
\item Can one decide exact comparisons such as $\pi(\cF)=q$ for a given rational $q$, or $\pi(\cF)=\pi(\cH)$ for two finite forbidden families?
\item Can one decide whether every nearly extremal hypergraph follows one fixed construction, or whether several separated classes of extremal structures occur?
\end{enumerate}
The first question has a positive answer for finite forbidden families, as recalled in \Cref{prop:approx}. Successively finer approximations eventually detect a strict inequality, but an interval containing both candidate values does not certify equality. We study exact comparison and structural classification by encoding the halting problem in these questions.

Fix a deterministic universal machine $\mathsf U$. Its input $\beta$ is a finite binary word, on which it either halts after finitely many steps or continues forever. Universality means that it can simulate any algorithm from a finite description of that algorithm and its input. No algorithm correctly decides, for every $\beta$, whether $\mathsf U$ halts on $\beta$~\cite{Turing1936,Rogers1967}. We choose the machine described in Section~\ref{ttm:sec:direct-tableau}, for which $\beta_\infty=\mathtt 0$ is provably nonhalting.

Write $\pi(P)\coloneqq\pi(\{P\})$. We first state the numerical result for a single forbidden $3$-graph.

\begin{theorem}\label{thm:singleton}
There are an effectively specified rational $\tau_3\in(0,1)$ and a total computable map sending each binary word $\beta$ to a single ordinary simple $3$-graph $P_\beta$ such that $\pi(P_\beta)=\tau_3$ if and only if $\mathsf U(\beta)$ does not halt, and $\pi(P_\beta)>\tau_3$ if and only if $\mathsf U(\beta)$ halts. The map is primitive recursive. The threshold is $\tau_3=6B_*$, where the rational value $B_*$ is specified in \Cref{lem:wi-rational-baseline}.
\end{theorem}
Li, Liu, and Liu~\cite[Theorem 1.2]{LiLiuLiu2026} already established rational-threshold undecidability. \Cref{thm:singleton} gives this conclusion for a single forbidden $3$-graph and one fixed rational threshold. Only the forbidden $3$-graph varies with $\beta$, so deciding equality with $\tau_3$ for these hypergraphs would decide the halting problem.

We next extend the density dichotomy to every $r\ge3$ and give a fixed construction attaining the threshold. A \emph{blowup} replaces each vertex by an independent class, possibly empty, and each edge by all transversals of its classes. A \emph{probability weighting} $\x$ on a finite $r$-graph $G$ assigns nonnegative weights to its vertices with total weight $1$. For any nonnegative vertex weighting $\x$, write
\[
 \lambda(G;\x)\coloneqq\sum_{e\in G}\prod_{v\in e}x_v,
 \qquad \lambda(G)\coloneqq\max_{\x}\lambda(G;\x),
 \qquad \lambda_r(G)\coloneqq r!\lambda(G).
\]
The maximum is over probability weightings. The normalized value $\lambda_r(G)$ is the largest limiting edge density among blowups of $G$. An $r$-graph is $\cF$-\emph{hom-free} if no member of $\cF$ has an edge-preserving map into it, equivalently if every blowup is $\cF$-free.

\begin{theorem}\label{thm:numerical}
For every fixed $r\ge3$, there are an effectively specified rational $\tau_r\in(0,1)$ and a finite $r$-graph $K_{*,r}$ with the following properties. From a word $\beta$ one can effectively construct a finite family $\cF_{r,\beta}$ satisfying both statements below.
\begin{enumerate}
\item We have $\pi(\cF_{r,\beta})=\tau_r$ if and only if $\mathsf U(\beta)$ does not halt, and $\pi(\cF_{r,\beta})>\tau_r$ if and only if $\mathsf U(\beta)$ halts.
\item Every blowup of $K_{*,r}$ is $\cF_{r,\beta}$-free for every input $\beta$. The uniform weighting of $K_{*,r}$ is optimal, and $\lambda_r(K_{*,r})=\tau_r$.
\end{enumerate}
For $r=3$, the output is the one-member family $\cF_{3,\beta}=\{P_\beta\}$. More generally, $\pi(\cF_{r,\beta})=\frac{27r!}{6r^r}\pi(P_\beta)$. The maps are total and computable uniformly in $r,\beta$, and primitive recursive in $\beta$ for fixed $r$. The threshold is given by the rational formula \eqref{eq:v3-public-threshold}.
\end{theorem}
Thus the same finite template $K_{*,r}$ supplies the density lower bound for every $\beta$, and the decision problem asks whether this fixed construction is asymptotically optimal. Strict improvement has a finite certificate. An $\cF_{r,\beta}$-hom-free $r$-graph $G$ with $\lambda_r(G)>\tau_r$ certifies halting within $v(G)$ steps. Conversely, a halting trace gives such a certificate with uniform vertex weights.

We next turn the density dichotomy into a distinction between the possible structures of nearly extremal hypergraphs. The numerical construction supplies the inner forbidden family $\{P_\beta\}$, and the structural construction converts its density gap into one or two phases.

Limits are identified by their finite homomorphism densities. Let $\cE(\cF)$ be the compact space of $\cF$-free limits of density $\pi(\cF)$; definitions are given in Section~\ref{sec:prelim}. Write $\cE(P)\coloneqq\cE(\{P\})$. A \emph{finite profile} specifies finitely many vertex roles, their proportions, and allowed multisets of roles on edges. We say that $\cF$ is \emph{stable with respect to a single finite profile} if one fixed finite profile describes every asymptotically extremal free sequence within $o(n^r)$ edge changes.

\begin{theorem}\label{thm:structural}
There is a computable $r_{\rm str}$ such that, for every fixed $r\ge r_{\rm str}$, one can effectively construct finite families $\cG_{r,\beta}$ and two unlabeled $r$-graphs $T_r^{(1)},T_r^{(2)}$. The latter and a finite-profile limit $W_{0,r}$ depend only on $r$. Set
\[
 \Sep_r(W)\coloneqq t(T_r^{(1)},W)-t(T_r^{(2)},W).
\]
The following alternatives hold.
\begin{enumerate}
\item If $\mathsf U(\beta)$ does not halt, then $\cE(\cG_{r,\beta})=\{W_{0,r}\}$, the separator is zero, and every asymptotically extremal sequence is $o(n^r)$-close to the fixed profile.
\item If $\mathsf U(\beta)$ halts, then $\cE(\cG_{r,\beta})$ is the disjoint union of two nonempty compact phases, on which $\Sep_r$ is respectively $c_{r,\beta}$ and $-c_{r,\beta}$ for some $c_{r,\beta}>0$. Each phase is homeomorphic to $\cE(P_\beta)$, the entire extremal limit space of the inner problem.
\end{enumerate}
\end{theorem}
In each of the 2 phases, the inner extremal limit may vary over all of $\cE(P_\beta)$. At the finite level, every asymptotically extremal $\cG_{r,\beta}$-free sequence in the halting case is $o(n^r)$-close to models built from asymptotically extremal $P_\beta$-free $3$-graphs by adding neutral vertices and root classes. Section~\ref{sec:transfer} gives the precise approximation and class-size rounding conclusions.

Consequently, five properties are undecidable, namely uniqueness of the extremal limit; connectedness of the extremal space; the separator taking both signs; separation into two nonempty regions on which the separator has opposite signs and is bounded away from zero; and stability with respect to a single finite profile.

The finite certificate conversions above can be verified in ZFC, so ZFC proves the numerical equivalence uniformly in $\beta$. The structural characterization can likewise be proved there. For a consistent computably axiomatized extension of ZFC, the G\"odel--Rosser theorem~\cite{Rosser1936} supplies a universal arithmetic sentence independent of that theory. We apply the constructions to a program that searches for a counterexample to this sentence. Since ZFC proves that the program does not halt exactly when the sentence holds, independence transfers to the corresponding density and structural assertions.
\begin{theorem}\label{thm:independence}
For every consistent computably axiomatized $T\supseteq\mathrm{ZFC}$, the following statements hold.
\begin{enumerate}
\item For every fixed $r\ge3$, one can effectively construct a finite family $\cF_{r,T}$ such that $\pi(\cF_{r,T})=\tau_r$ is true but neither provable nor refutable in $T$. For $r=3$, this family consists of a single $3$-graph $P_T$. The theory proves that $K_{*,r}$ is $\cF_{r,T}$-hom-free and that $\lambda_r(K_{*,r})=\tau_r$, and hence proves the density lower bound.
\item For every fixed $r\ge r_{\rm str}$, there is a finite family for which all five structural properties above are independent of $T$. Uniqueness, connectedness and stability hold, while the two sign properties fail.
\end{enumerate}
\end{theorem}
The constructions also determine the complexity of exact comparison and the size of finite improvement certificates. \Cref{thm:finitecomparison} classifies strict comparison and equality, including comparisons between single forbidden $3$-graphs. \Cref{thm:witnessbounds} shows that the minimum order of an improving template is comparable to the area of the halting history, and is at most quadratic in the halting time and input length. Maximizing this minimum order over halting inputs of bounded length gives Busy-Beaver growth, while the positive density gaps admit no computable lower bound in terms of the encoded input.

For effectively enumerated forbidden families, we ask which values in the unrestricted density spectrum admit such a presentation. \Cref{thm:spectrum} characterizes them as the values with computable nonincreasing rational approximations. Each such value also has a presentation with decidable membership in which no member contains another.

\paragraph{Proof outline.}
The numerical construction must make a strict density improvement possible precisely when the encoded computation halts. We represent a finite computation by 2 vertex types and 3 partial matchings, encoded by root links in ordinary $3$-graphs. A single edge on 3 data vertices marks the computation and increases its weighted density. A halting history can be closed into a periodic tableau with every required matching incidence present. If the computation does not halt, each marked component has a missing incidence, and a path estimate shows that the resulting density loss absorbs the possible gain from the marked edge. A further inequality controls vertices whose root links do not determine a unique code. Together these estimates establish the fixed rational baseline and allow a halting history to be recovered from any improvement.

The next task is to realize this density dichotomy with a single forbidden $3$-graph. Root elimination first gives a finite forbidden basis. We enlarge the admissible class without changing its Lagrangian supremum, then form the categorical product of the canonical extensions of its new forbidden cores. The projections onto the factors give the upper bound. For the lower bound, we must exclude homomorphisms from the product into the constructions attaining or exceeding the baseline. A selection argument extracts compatible maps on one factor from any such homomorphism. Unique completion in the partial matchings and the single marked edge then reconstructs a legal rooting of a forbidden core, a contradiction. Effective suspension extends the resulting numerical theorem to every $r\ge3$.

The structural construction optimizes the proportion of active vertices as a function of the inner density gap. The optimal proportion is zero at zero gap and positive at positive gap. In the latter case, the 2 possible labels for the inner edges give 2 phases, each retaining the entire inner extremal space. To control partial models, we reconstruct complete legal objects from partial records and show that the score lost through missing records dominates the possible gain from reconstruction. Unlabeled homomorphism-density tests recover every inner homomorphism density, giving the homeomorphisms in \Cref{thm:structural}. Finally, vertex extendability and the stability theorem of Liu, Mubayi and Reiher~\cite{LMR2023} show that every ordinary near-extremal hypergraph is close to the corresponding model class.

Section~\ref{sec:prelim} gives the common transfer lemmas and limit-space definitions. Sections~\ref{ttm:sec:direct-tableau} and~\ref{sec:v3-three} establish the computational and weighted constructions. Section~\ref{sec:single} produces the single forbidden $3$-graph, and Section~\ref{sec:v3-suspension} proves the higher-uniformity and certificate statements. Sections~\ref{sec:limits}--\ref{sec:transfer} give the structural classification and stability. Section~\ref{sec:consequences} treats exact comparison, witness bounds, and logical independence. Section~\ref{sec:spectrum} treats effectively presented families, oracle variants, and finite schemes for exact values.

\section{Finite transfers and limit spaces}\label{sec:prelim}
Write $\Forb(\cF)$ for the finite $\cF$-free $r$-graphs and $\Lambda(\cD)\coloneqq\sup_{G\in\cD}\lambda(G)$ for a class $\cD$. A class is \emph{monotone} if it is closed under arbitrary subgraphs. An edge-preserving map is injective on each edge. Let $\cQ(\cF)$ be the family of all homomorphic images of members of $\cF$. For finite $\cF$, the family $\cQ(\cF)$ can be computed by enumerating partitions of the vertex sets of its members.

For a vertex set $A$ write $x(A)\coloneqq\sum_{v\in A}x_v$. For probability weights define
\[
 t(F,G;\x)\coloneqq\sum_{\phi\colon V(F)\to V(G)}
   \prod_{v\in V(F)}x_{\phi(v)}\prod_{e\in F}\mathbf1_{\phi(e)\in G}.
\]
When the weighting is uniform, we omit $\x$ from the notation. In particular $t(K_r^{(r)},G;\x)=r!\lambda(G;\x)$. Splitting weights among clones preserves every homomorphism density.

If $F\to G$, then $F$ embeds in a blowup of $G$. Given a homomorphism $\phi\colon F\to G$ and a probability weighting $\x$ on $F$, let $y_v\coloneqq\sum_{u\in\phi^{-1}(v)}x_u$ for $v\in V(G)$. Expanding the edge products gives $\lambda(F;\x)\le\lambda(G;\y)$. Maximizing over $\x$ yields $\lambda(F)\le\lambda(G)$.

For a vector $\x$, use $e_j(\x)$ for its elementary symmetric polynomial, with $e_0=1$; write $e_j(A)$ for its restriction to a vertex set $A$. For nonnegative $\x$, we have $e_j(\x)\le(\sum x_i)^j/j!$.

Let $\mathfrak W_r$ be the closure in $[0,1]^{\mathbb N}$ of the homomorphism-density vectors of finite $r$-graphs, with coordinates indexed by finite test $r$-graphs $F_1,F_2,\ldots$ and the topology of coordinatewise convergence. This is compact and agrees with the usual dense-limit description~\cite{Austin2008,ElekSzegedy2012,Zhao2015}.

The next lemma describes free limit vectors by the vanishing of the forbidden densities.
\begin{lemma}\label{lem:freerealization}
For any family $\cF$ of $r$-graphs, a vector $W\in\mathfrak W_r$ is a limit of $\cF$-free $r$-graphs of orders tending to infinity if and only if $t(F,W)=0$ for every $F\in\cF$.
\end{lemma}
\begin{proof}
For necessity, let $G$ be an $F$-free $r$-graph. Every homomorphism $F\to G$ identifies two vertices of $F$, so $t(F,G)\le\binom{v(F)}2/v(G)$. This tends to zero along any free sequence whose orders tend to infinity.

For sufficiency, suppose that $t(F,W)=0$ for every $F\in\cF$. We first construct a free $r$-graph approximating $W$ on any prescribed finite set of tests. Let $S_m(G)$ be the random $r$-graph on $[m]$ obtained by sampling $m$ vertices independently and uniformly from $V(G)$, with replacement, and putting an edge on an $r$-set of indices exactly when their images form an edge of $G$. For an $a$-vertex test $F$ and $m\ge a$, let $Z_F$ be the proportion of injections from $V(F)$ to $[m]$ that are homomorphisms from $F$ to $S_m(G)$. Then
\[
 \mathbb E[Z_F]=t(F,G),\quad
 \operatorname{Var}[Z_F]\le \frac{a^2}{m},\quad
 |t(F,S_m(G))-Z_F|\le\frac{\binom a2}{m}.
\]
The variance bound follows because disjoint index sets give independent indicators; the last bound counts noninjective index maps. Bounding the probability of an $F$-copy by the expected number of injective homomorphisms from $F$ to $S_m(G)$ gives $\Pr[F\subseteq S_m(G)]\le m^a t(F,G)$.

Fix finitely many tests and $\epsilon>0$. Choose $m$ so that, uniformly in $G$, all their densities in $S_m(G)$ and $G$ differ by less than $\epsilon/2$ with probability exceeding $3/4$. Choose $G$ approximating $W$ within $\epsilon/2$ on these tests and satisfying $\sum_{F\in\cF,\,v(F)\le m}m^{v(F)}t(F,G)<1/4$, with one representative of each isomorphism type. Some sample is free and approximates $W$ within $\epsilon$ on these tests.

At stage $j$, apply this construction to $F_1,\ldots,F_j$ with error $1/j$, choosing $m\ge j$. The resulting free $r$-graphs have orders tending to infinity and density vectors converging to $W$.
\end{proof}

The hom-free enumeration method for effective approximation is described in Keevash~\cite[Section~2]{Keevash2011}. We record the finite-family form, together with the density identities needed below.
\begin{lemma}\label{prop:approx}
For every family $\cF$,
\begin{equation}\label{eq:BS}
 \pi(\cQ(\cF))=\pi(\cF),\qquad
 \pi(\cF)=\inf_{\substack{\cF'\subseteq\cF\\\cF'\text{ finite}}}\pi(\cF'),
\end{equation}
and
\begin{equation}\label{eq:lagdensity}
 \pi(\cF)=\sup\{\lambda_r(G)\colon G\text{ finite and }\cF\text{-hom-free}\}.
\end{equation}
For finite $\cF$ and $m\ge r$, let $p_m\coloneqq\ex(m,\cQ(\cF))/\binom mr$. Then
\begin{equation}\label{eq:approx}
 \frac{r!\ex(m,\cQ(\cF))}{m^r}\le\pi(\cF)\le p_m,
 \qquad 0\le p_m-\pi(\cF)\le\frac{r(r-1)}{2m}.
\end{equation}
The lower and upper endpoints are computable uniformly from $r$, $m$, and $\cF$.
\end{lemma}
\begin{proof}
We first prove the classical identity $\pi(\cQ(\cF))=\pi(\cF)$~\cite{BrownSimonovits1984,Pikhurko2014,Keevash2011} by showing that vanishing homomorphism density for $F$ forces vanishing density for each homomorphic image $Q$ of $F$. Fix a finite target $r$-graph $G$ and a vertex $u$ of $Q$. We compare $t(Q,G)$ with the homomorphism density obtained by cloning $u$. First sample the images of all other vertices independently and uniformly from $V(G)$. Let $A\in\{0,1\}$ indicate that all edges avoiding $u$ are preserved, and let $B\in[0,1]$ be the conditional probability that an independent image of $u$ preserves all edges containing $u$. Then $t(Q,G)=\mathbb E[AB]$. Replace $u$ by $s\ge1$ clones and sample their images independently and uniformly from $V(G)$. The resulting homomorphism density is $\mathbb E[AB^s]$. Since $A\in\{0,1\}$, we have $AB^s=(AB)^s$. Jensen's inequality applied to $AB$ therefore gives $\mathbb E[AB^s]\ge(\mathbb E[AB])^s$.

Let $Q(s)$ be the blowup with $s$ clones at every vertex. If $F\to Q$ and $s\ge v(F)$, then $F\subseteq Q(s)$. Iterating the preceding inequality over the vertices of $Q$ gives
\[
 t(F,G)\ge t(Q(s),G)\ge t(Q,G)^{s^{v(Q)}}.
\]
Thus $t(F,G_n)\to0$ forces $t(Q,G_n)\to0$ for every homomorphic image $Q$ of $F$. Since $\cF\subseteq\cQ(\cF)$, the sets of limit vectors on which all forbidden homomorphism densities vanish are the same for these two families. By \Cref{lem:freerealization}, the maximum edge density on this common set equals both $\pi(\cF)$ and $\pi(\cQ(\cF))$.

For finite-subfamily approximation, let $\cF_{\le n}$ consist of one representative of each isomorphism class of members of $\cF$ with at most $n$ vertices. Then $\pi(\cF)\le\pi(\cF_{\le n})\le\ex(n,\cF)/\binom nr\to\pi(\cF)$.

We next prove \eqref{eq:lagdensity}. For every finite $\cF$-hom-free $G$, blowups with rational proportions approaching an optimal weighting give $r!\lambda(G)\le\pi(\cF)$. Conversely, an $n$-vertex extremal member $H_n$ of $\Forb(\cQ(\cF))$ is $\cF$-hom-free, so
\[
 \sup_{G\text{ finite }\cF\text{-hom-free}}\lambda_r(G)
 \ge\lambda_r(H_n)\ge\frac{r!|H_n|}{n^r}
 \longrightarrow\pi(\cF).
\]
These two inequalities prove \eqref{eq:lagdensity}.

For the effective bounds in \eqref{eq:approx}, now assume that $\cF$ is finite. Choose an $m$-vertex $\cQ(\cF)$-free extremal $H$. Its balanced blowups have limiting density $r!|H|/m^r$, giving the lower bound. Averaging an $n$-vertex $\cQ(\cF)$-free $r$-graph over its $m$-vertex subsets, for $n\ge m$, gives the upper bound $p_m$. Since $0\le p_m\le1$, these bounds imply
\[
 \begin{aligned}
 0\le p_m-\pi(\cF)
 &\le p_m\left[1-\frac{r!\binom mr}{m^r}\right]\\
 &\le1-\prod_{j=0}^{r-1}\left(1-\frac jm\right)
 \le\sum_{j=0}^{r-1}\frac jm=\frac{r(r-1)}{2m}.
 \end{aligned}
\]
Enumeration computes each finite extremal number. Deleting isolated vertices from a forbidden hypergraph does not change its homomorphism density, so all assertions also apply after deleting them.
\end{proof}
The extremal space introduced in Section~\ref{sec:intro} therefore has the coordinate description
\[
 \cE(\cF)=\bigl\{W\in\mathfrak W_r\colon
 t(F,W)=0\ \text{for every }F\in\cF,\quad
 t(K_r^{(r)},W)=\pi(\cF)\bigr\}.
\]
In particular, it is a nonempty compact subset of $\mathfrak W_r$.

For $r$-graphs $G\subseteq G'$ on the same vertex set, fix a probability weighting $\x$. Coupling the vertex samples and taking a union bound over the edges of $F$ gives $0\le t(F,G';\x)-t(F,G;\x)\le r!|F|[\lambda(G';\x)-\lambda(G;\x)]$.

For a fixed $r$-graph $G$ and probability weightings $\x,\y$ on its vertex set, couple draws from $\x$ and $\y$ with disagreement probability at most $\|\x-\y\|_1$. Independent copies of this coupling for the vertices of $F$ and a union bound give $|t(F,G;\x)-t(F,G;\y)|\le v(F)\|\x-\y\|_1$.

The numerical and structural constructions first describe legality using finitely many named roots. We use the bounded-search criterion of Li, Liu, and Liu~\cite[Lemmas 8.3--8.4]{LiLiuLiu2026} to remove those names when every failed rooting has a bounded witness. For completeness, we give the argument with the explicit order bound needed for the constructions.
\begin{lemma}\label{lem:unroot}
Fix a finite label set $M$. A \emph{partial rooting} is an injection $\mu\colon M\rightharpoonup V(G)$; write $\nu\supseteq\mu$ when $\nu$ extends $\mu$ as a partial map. Suppose that admissibility of a partial rooting $\mu$ of $G$ is specified by a decidable predicate on $(G,\mu)$, invariant under hypergraph isomorphisms that preserve every assigned root label. Let $\cD$ consist of the $r$-graphs having an admissible partial rooting. Suppose that $\cD$ is monotone and contains every finite edgeless $r$-graph. Assume every failed rooting has a witness $F_\mu\subseteq G$ of order at most $w$ such that, whenever $F_\mu\subseteq H\subseteq G$ and $\nu\supseteq\mu$ is admissible on $H$, some new root lies in $V(F_\mu)\setminus\operatorname{im}\mu$. Then $\cD=\Forb(\cA)$ for an effective finite family whose members have order at most
\[
 w\sum_{j=0}^{|M|}(|M|w)^j.
\]
\end{lemma}
\begin{proof}
In $G\notin\cD$, start at the empty rooting. At a failed node $\mu$, take one witness $F_\mu$ and create the children
\[
 \mu\cup\{a\mapsto x\},\qquad
 a\in M\setminus\operatorname{dom}\mu,\quad
 x\in V(F_\mu)\setminus\operatorname{im}\mu.
\]
Each step assigns exactly one new label. The depth is therefore at most $|M|$ and the branching at most $|M|w$. The union $F$ of all witnesses has the stated order; every previously assigned vertex lies in an ancestor witness and hence in $F$.

To show that $F\notin\cD$, suppose that $\rho$ is an admissible rooting of $F$. Starting at the empty rooting, consider a node $\mu$ that is a restriction of $\rho$. Apply the witness hypothesis with $H=F$ and $\nu=\rho$. It supplies a label $a\in\operatorname{dom}\rho\setminus\operatorname{dom}\mu$ with $\rho(a)\in V(F_\mu)\setminus\operatorname{im}\mu$, so the corresponding child is again a restriction of $\rho$. Repeating this step would require another label after all labels used by $\rho$ had been assigned, a contradiction.

Enumerate all $r$-graphs outside $\cD$ whose order is at most the stated bound. They all have an edge, since every edgeless $r$-graph belongs to $\cD$. The preceding construction shows that every $G\notin\cD$ contains a listed $r$-graph. Conversely, monotonicity prevents a member of $\cD$ from containing any listed $r$-graph. Thus the list is a forbidden basis for $\cD$.
\end{proof}

The pair-covering transfer of Li, Liu, and Liu~\cite{LiLiuLiu2026} converts the Lagrangian supremum of a finitely forbidden class into a finite-family Tur\'an density. Pair-covering extensions also appear in Mubayi~\cite{Mubayi2006}.

For an $r$-graph $A$ with an edge, let $\Ext(A)$ consist of all $r$-graphs without isolated vertices of order at most $v(A)+(r-2)\binom{v(A)}2$ containing a copy of $A$ whose vertex pairs are covered by edges of the extension. Every homomorphism from an extension is therefore injective on this copy. For completeness, we give the transfer proof and its support conclusion.
\begin{lemma}\label{lem:pairtransfer}
Suppose that $\cA$ is finite and $\Forb(\cA)$ contains an edge. Then the following statements hold.
\begin{enumerate}
\item We have $\pi(\bigcup_{A\in\cA}\Ext(A))=r!\Lambda(\Forb(\cA))$.
\item Every blowup of a member of $\Forb(\cA)$ avoids every member of $\bigcup_{A\in\cA}\Ext(A)$.
\item For an optimal weighting with minimum support size on an $r$-graph that has at least one edge and avoids every member of $\bigcup_{A\in\cA}\Ext(A)$, the subhypergraph induced by the positive-weight vertices belongs to $\Forb(\cA)$ and covers pairs.
\end{enumerate}
\end{lemma}
\begin{proof}
If a blowup of an $\cA$-free $r$-graph contained an extension of $A\in\cA$, the projection to the base would be injective on the pair-covered $A$-copy and would embed $A$ into the base, contradicting $\cA$-freeness. This proves blowup freeness and the density lower bound.

For the reverse bound, we use a standard minimum-support argument, as in Frankl and R\"odl~\cite{FranklRodl1984}. Let $H$ be an extension-free $r$-graph with at least one edge, choose an optimal weighting with minimum support size, and let $J$ be the $r$-graph induced by that support. If two distinct vertices $u,v\in V(J)$ lie in no common edge, then, with all other weights fixed, the edge polynomial has the form $C+x_uL_u+x_vL_v$. Move their combined weight to the vertex with the larger coefficient. The resulting weighting is still optimal and has smaller support, a contradiction.

Since $J$ covers pairs, a copy of $A\in\cA$ in $J$, together with a covering edge for each pair, would give a forbidden extension satisfying the prescribed order bound. Thus $J\in\Forb(\cA)$. For $n\coloneqq v(H)$, uniform weights on $H$ give
\[
 \frac{|H|}{n^r}\le\lambda(H)=\lambda(J)\le\Lambda(\Forb(\cA)).
\]
Taking $n\to\infty$ gives the reverse density bound.
\end{proof}

\paragraph{Rational blowups.}
Let $H$ be $\cF$-hom-free, and let $\x$ be a rational probability weighting. If $Dx_v$ is an integer for every vertex, replacing $v$ by $Dx_v$ clones gives a hom-free $r$-graph $B$ of order $D$ with
\begin{equation}\label{eq:rationalblowup}
 \frac{r!|B|}{D^r}=r!\lambda(H;\x).
\end{equation}
Indeed, the uniform weights on each clone class sum to $x_v$, so the identity follows by summing over edge transversals.

Now suppose that $(H,\x)$ is independent of $\beta$, that $H$ is $\cF_\beta$-hom-free for every $\beta$, and that $r!\lambda(H;\x)=\tau$. If $\pi(\cF_{\beta_\infty})=\tau$ for one input, then $\tau\le\lambda_r(B)\le\pi(\cF_{\beta_\infty})=\tau$. Thus $B$ is a fixed baseline template whose uniform weighting is optimal.

If $\pi(\cF)>q$, then \eqref{eq:lagdensity} gives a finite $\cF$-hom-free template $H$ with a probability weighting $\x$ such that $r!\lambda(H;\x)>q$. Approximate its weights rationally while preserving the strict inequality, and take the corresponding rational blowup. The resulting finite $\cF$-hom-free $r$-graph $G$, with uniform weights, satisfies
\begin{equation}\label{eq:finitecertificate}
 r!|G|>qv(G)^r.
\end{equation}
Conversely, the balanced blowups of any such $G$ give $\pi(\cF)\ge r!|G|/v(G)^r>q$. This equivalence will turn density comparisons into searches for finite certificates.

\section{Computation tableaux and partial matchings}\label{ttm:sec:direct-tableau}
We encode a halting computation by a finite periodic tableau. A locally certified reset closes its time direction, and 3 partial matchings represent its horizontal and vertical successors. The size of a closed marked component will also control the order of an improvement certificate.

\subsection{A locally certified reset}
Fix a deterministic universal Turing machine $\mathsf U$ on a one-sided tape, normalized to run forever on the empty input and on inputs beginning with $\mathtt0$, and to simulate a fixed universal machine on $\gamma$ when its input is $\mathtt1\gamma$. Its word-halting problem is undecidable \cite{Turing1936,Rogers1967}. For a halting input $\beta$, write $T(\beta)$ for the first halting time and $S(\beta)$ for the largest visited tape index. The machine has a read-only track with alphabet $\mathcal I=\{\mathtt0,\mathtt1,\mathtt\$,\mathtt{pad}\}$, a work alphabet $\Gamma$ with blank $\square$, a state set $\mathcal Q$, an initial state $q_0\ne q_{\rm h}$, and a halting state $q_{\rm h}$. Its fixed transition table is
\[
 \delta(q,i,a)=(q',b,\varepsilon),\qquad\varepsilon\in\{-1,+1\},
\]
for every nonhalting state $q$. A left move at cell $1$ leaves the head there. On the empty word and words beginning with $\mathtt0$, use a fixed nonhalting state that moves left forever without changing the work symbol. These conventions preserve universality, make $\beta_\infty=\mathtt0$ provably nonhalting, and ensure $T(\beta)\ge1$ on every halting input.

A row lies between fixed walls $\mathtt L,\mathtt R$. To make the time direction periodic, we certify locally that a row contains a halting head and then reset the entire row in one step. Let
\[
 D=\Gamma\sqcup(\mathcal Q\times\Gamma)\sqcup\{\mathtt{Crash}\},
 \qquad
 \Sigma=\{\mathtt L,\mathtt R\}\sqcup(\mathcal I\times\mathcal D),
\]
where the dynamic alphabet is the disjoint union
\[
 \mathcal D=\{\mathtt I\}\sqcup\{\mathrm{run}(d)\colon d\in D\}
 \sqcup\{\mathrm{pre}(a),\mathrm{halt}(a),\mathrm{post}(a)\colon a\in\Gamma\}.
\]
The decorations $\mathrm{pre}(a)$, $\mathrm{halt}(a)$, and $\mathrm{post}(a)$ all have underlying work symbol $a$; $\mathrm{halt}(a)$ additionally represents the head $(q_{\rm h},a)$. The symbol $\mathtt I$ represents an initial blank cell.

The horizontal grammar permits exactly the following successors, when the read-only track is omitted.
\[
\begin{array}{c|l}
\text{symbol}&\text{allowed successor}\\ \hline
\mathtt L&\mathtt I,\ \mathrm{run}(d),\ \mathrm{pre}(a),\ \mathrm{halt}(a)\\
\mathtt I&\mathtt I,\ \mathtt R\\
\mathrm{run}(d)&\mathrm{run}(d'),\ \mathtt R\\
\mathrm{pre}(a)&\mathrm{pre}(b),\ \mathrm{halt}(b)\\
\mathrm{halt}(a)&\mathrm{post}(b),\ \mathtt R\\
\mathrm{post}(a)&\mathrm{post}(b),\ \mathtt R\\
\mathtt R&\mathtt L.
\end{array}
\]
Thus every finite wall segment is entirely initial, entirely running, or has dynamic word
\[
 \mathrm{pre}(\Gamma)^*\mathrm{halt}(\Gamma)\mathrm{post}(\Gamma)^*.
\]
In particular, the last case certifies exactly one halting head by adjacent-pair conditions alone. Add the read-only rules $\mathtt\$\to\mathtt{pad}$ and $\mathtt{pad}\to\mathtt{pad}$ or $\mathtt R$. Let $\mathcal G\subseteq\Sigma^2$ denote this finite pair relation.

We next define a finite transition relation $\mathcal R\subseteq\Sigma^4$; its first 3 entries are a radius-one window, and its last entry is the next symbol at its middle position. Walls stay fixed regardless of their neighbors, and every interior update preserves its read-only symbol. The dynamic update is as follows.
\begin{enumerate}
\item A middle $\mathtt I$ becomes $\mathrm{run}(q_0,\square)$
if its left neighbor is $\mathtt L$, and $\mathrm{run}(\square)$ otherwise.
\item A middle reset decoration becomes $\mathtt I$.
\item At a running middle position, perform the usual Turing update
on the underlying running symbols, obtaining $d\in D$. The output may be $\mathrm{run}(d)$; if $d=a\in\Gamma$, it may instead be $\mathrm{pre}(a)$ or $\mathrm{post}(a)$; and if $d=(q_{\rm h},a)$, it may instead be $\mathrm{halt}(a)$.
\end{enumerate}
For completeness the underlying running update is determined locally as follows. A window with a nonrunning interior neighbor, a crash, a halting head, or at least 2 heads outputs $\mathtt{Crash}$. Otherwise a window with no head leaves its middle work symbol unchanged. A unique nonhalting head executes $\delta$. The old head position receives the written work symbol. The destination receives the new head and retains its existing work symbol. At the left wall a requested left move leaves the head in place with the written symbol; a requested right move through the right wall crashes the old head position. A neighbor not receiving the head stays unchanged. Any remaining malformed case outputs $\mathtt{Crash}$. These rules give an explicit finite table. The horizontal relation $\mathcal G$ requires every reset row to have exactly one certified halting head.

Write $\Sigma_{\rm ini}=\{(i,\mathtt I)\colon i\in\mathcal I\}$ and $\Sigma_{\rm inp}(i)=\{(i,d)\colon d\in\mathcal D\}$. The constant background symbol $b_0=(\mathtt{pad},\mathrm{run}(\square))$ satisfies
\begin{equation}
 (b_0,b_0)\in\mathcal G,\qquad
 (b_0,b_0,b_0,b_0)\in\mathcal R.
 \label{ttm:eq:constant-background-update}
\end{equation}

For $W\ge2$, an admissible row on $0,\ldots,W+1$ has fixed endpoint walls and satisfies $\mathcal G$ on adjacent positions. Write $c\Rightarrow_W c'$ when both rows are admissible and their radius-one windows and new middle symbols satisfy $\mathcal R$ at every interior position. For an input $\beta$ of length $\ell$ and $W\ge\ell+2$, let $c_{\beta,W}$ have read-only word $\beta\mathtt\$\mathtt{pad}^{W-\ell-1}$ and dynamic word $\mathtt I^W$.

\begin{lemma}\label{ttm:lem:reset-ca}
For every input $\beta$, the machine $\mathsf U(\beta)$ halts if and only if there are integers $W\ge|\beta|+2$ and $p\ge1$ such that $c_{\beta,W}\Rightarrow_W^p c_{\beta,W}$.
\end{lemma}
\begin{proof}
The first step launches precisely the blank-work-tape computation, with its head in cell $1$. Until a reset phase or a right-boundary crossing occurs, the underlying running configurations are the genuine deterministic computation. Indeed, the horizontal grammar forces a whole wall segment to use one phase, and the transition relation preserves the ordinary update under every permitted decoration. A reset phase is possible only when this update has produced an actual halting head; its local prefix--head--suffix certificate excludes a reset without one. A right-boundary crossing produces a persistent crash, which admits no reset decoration and therefore precludes return to an initial row.

If the machine halts, choose $W\ge\max\{S(\beta),|\beta|+2\}$. Decorate the resulting halting row by $\mathrm{pre}$ before its head, $\mathrm{halt}$ at its head, and $\mathrm{post}$ afterwards. The next step resets every interior cell to $\mathtt I$, preserving the read-only word. This gives the asserted finite return.

Conversely, a first return to the initial phase must be preceded by a certified reset phase. The preceding running history is genuine, so its halting head proves that $\mathsf U(\beta)$ halts.

In any returning sequence, the first halting update must receive reset decorations. Leaving it in the running phase would produce a persistent crash at the next step. Hence the first return occurs after exactly $T(\beta)+2$ transitions.
\end{proof}

For a halting input $\beta$, the reset construction uses an initial row, $T(\beta)$ running rows, and a certified halting row. Define its \emph{history area} by
\begin{equation}\label{eq:historyarea}
 W_\beta\coloneqq\max\{S(\beta),|\beta|+2\},\qquad
 A_{\rm hist}(\beta)\coloneqq(T(\beta)+2)(W_\beta+2),
\end{equation}
and set $A_{\rm hist}(\beta)=\infty$ when the machine does not halt. Here $W_\beta$ counts the interior cells, so each row has $W_\beta+2$ vertices including the walls.

The space--time description uses the usual local encoding of computation~\cite{Wang1961,Berger1966,Robinson1971}. The reset grammar makes the return condition a finite system of local relations.

\subsection{Two computation types and 3 matchings}\label{sec:two-type-computation-module}

Write $A_{\rm CA}=|\Sigma|$ for the size of the fixed alphabet above. We use its horizontal relation $\mathcal G$, transition relation $\mathcal R$, and the reset equivalence of \Cref{ttm:lem:reset-ca}.

There are just 2 computation types, $\mathsf C$ and $\mathsf B$. A vertex of type $\mathsf C$ is a cell. A vertex of type $\mathsf B$ carries one color of $\Sigma$. The signature has 3
\emph{distinct indexed relation symbols}
\[
 J\colon\mathsf C\longrightarrow\mathsf B,
 \qquad R\colon\mathsf B\longrightarrow\mathsf C,
 \qquad S\colon\mathsf B\longrightarrow\mathsf C.
\]
Each relation is a partial matching. The union of the 3 relations need not be a matching. Later we take unions only over copies indexed by the same relation symbol. The underlying type multigraph has 3 parallel indexed edges; both computation types have incidence degree 3.

A realization consists of finite classes $V_{\mathsf C},V_{\mathsf B}$, a color map $\chi\colon V_{\mathsf B}\to\Sigma$, the 3 partial matchings, and a subset $Z\subseteq V_{\mathsf C}$ of marked cells. Define $H(x,y)$ to hold when some $b$ satisfies $J(x,b)$ and $R(b,y)$, and define $V(x,y)$ to hold when some $b$ satisfies $J(x,b)$ and $S(b,y)$. Thus the 2 successors use the paths
\[
 x\mathrel{\mathop{\longrightarrow}\limits^J}b\mathrel{\mathop{\longrightarrow}\limits^R}y\qquad\text{and}\qquad x\mathrel{\mathop{\longrightarrow}\limits^J}b\mathrel{\mathop{\longrightarrow}\limits^S}y.
\]
A present label $\operatorname{lab}(x)=a$ means that $J(x,b)$ is present for a vertex $b$ with $\chi(b)=a$. It is unique whenever it is present, by the matching property of $J$.

The same $J$-edge can serve simultaneously in a horizontal path, a temporal path, and a label. These witnesses are not required to be disjoint. This sharing imposes no additional restriction on the horizontal and temporal successors. Indeed, given any 2 permutations $H,V$ of a finite cell set, take a copy of that set for $V_{\mathsf B}$, choose $J$ to be the identity pairing, and choose $R=H$ and $S=V$ under that identification. Assign the $\mathsf B$-colors to be the desired cell labels.

\paragraph{Local conditions.}

All paths and labels mentioned in a condition must be present for that condition to apply. Powers of $H$ below mean concatenated paths, with each step expanded into its $J,R$ pair.
\begin{enumerate}
\item[(Q1)] \emph{Commutation.} If
$H(x,x_E),V(x,x_N),V(x_E,y),H(x_N,z)$ are present, then $y=z$.
\item[(Q2)] \emph{Evolution.} If $H(x_{\rm left},x),H(x,x_{\rm right}),V(x,x')$
and the 4 labels $a_{\rm left},a_0,a_{\rm right},b$ at $x_{\rm left},x,x_{\rm right},x'$ are present, then $(a_{\rm left},a_0,a_{\rm right},b)\in\mathcal R$.
\item[(Q3)] \emph{Row grammar.} On a present $H(x,y)$ with both
endpoint labels present, require
\[
 (\operatorname{lab}(x),\operatorname{lab}(y))\in\mathcal G.
\]
\item[(Q4)] \emph{Input at a marked cell.} Write
$\beta=\beta_1\cdots\beta_\ell$. For each $a\in Z$, impose the following conditions.
\begin{enumerate}
\item If $H(y,a)$ and both labels are present, then
$\operatorname{lab}(y)=\mathtt L$ and $\operatorname{lab}(a)\in\Sigma_{\rm ini}$.
\item For $1\le j\le\ell$, a present path $H^{j-1}(a,z)$
and a present label at $z$ require $\operatorname{lab}(z)\in\Sigma_{\rm inp}(\beta_j)$.
\item A present path $H^\ell(a,z)$ and a present label at $z$
require $\operatorname{lab}(z)\in\Sigma_{\rm inp}(\mathtt{\$})$.
\end{enumerate}
\end{enumerate}
A path of length zero has its single endpoint; its label is still separately required in the prefix conditions. Thus a marked cell with no incident relations activates none of Q4. Marks merely select which local conditions are imposed. No fixed cardinality or separation condition on $Z$ belongs to the relational module. The halting construction below uses 3 marked cells as the endpoints of one data-only edge.

\begin{lemma}
\label{lem:two-type-positive-obstructions}
For each input $\beta$, the matching conditions and Q1--Q4 are equivalent to avoidance of an effective finite list of typed positive patterns. Their sizes are effectively bounded, with all $\beta$-dependent sizes $O(|\beta|+1)$. Legality is preserved by deleting vertices, relation instances, or marks. Every forbidden positive pattern is connected in the relation multigraph. Thus every occurrence lies in a single relation component, containing its marked cell when the pattern is marked.
\end{lemma}
\begin{proof}
A matching conflict is a fork of 2 instances of the same relation with distinct partners. Expand each displayed path into matching instances, and expand each label into its $J$-edge and the specified color of its $\mathsf B$-endpoint. For an implication, enumerate its finitely many forbidden color alternatives, or retain the inequality between the 2 endpoints of a failed equality. Each marked pattern also includes its unary mark.

These formal patterns can have overlapping realizations. Enumerate all their quotients that preserve types, all prescribed colors, and all declared inequalities. In particular, 2 formal $\mathsf B$ vertices carrying different prescribed colors cannot be identified. A vertex with no prescribed color can be identified with a colored one. Every actual violation factors through one of these quotients followed by an injective copy; conversely every included copy is a violation. This explicitly allows the shared $J$-edges in successor and label paths. No disjointness assumption is introduced.

Only Q4 depends on $\beta$. Its powers have at most $2|\beta|$ matching edges, and the predecessor condition and label witnesses add only a bounded number of vertices, giving the stated bound. The patterns require only present instances, colors, and marks. Deletions cannot create them. The relation edges in each pattern connect all its vertices, so every occurrence lies in a single relation component.
\end{proof}

\subsection{Closed components and recovery}

The relation multigraph has one undirected edge for every instance of $J,R,S$, retaining their different labels even when endpoints coincide. Call a component \emph{closed} when every vertex has all incidences prescribed by its type, and \emph{incomplete} otherwise. Equivalently a closed computation component is 3-regular in this indexed multigraph.

\begin{lemma}
\label{lem:two-type-closed-recovery}
Let $a$ be a marked cell in a finite legal realization. If its relation component is closed, then $\mathsf U(\beta)$ halts. More precisely, if the component has $N$ cell vertices, there exist integers $W,p$ with
\[
 W\ge|\beta|+2,\qquad W+2\le N,\qquad1\le p\le N,
 \qquad c_{\beta,W}\Rightarrow_W^p c_{\beta,W}.
\]
Moreover $A_{\rm hist}(\beta)\le N$ and $T(\beta)\le N$. Consequently, if $\mathsf U(\beta)$ does not halt, every marked cell in a legal realization for $\beta$ lies in an incomplete component.
\end{lemma}
\begin{proof}
Closedness makes $J,R,S$ bijections between the component's cell and color classes, each of size $N$. Thus $H=R\circ J$ and $V=S\circ J$ are permutations, every cell has a label, and Q1--Q3 give
\[
 HV=VH,\qquad
 (\operatorname{lab}(H^{-1}x),\operatorname{lab}(x),
 \operatorname{lab}(Hx),\operatorname{lab}(Vx))\in\mathcal R,
\]
as well as the horizontal grammar $\mathcal G$.

Q4 gives a left wall at $H^{-1}a$, initial phase at $a$, the input $\beta$ at $a,\ldots,H^{\ell-1}a$, and $\mathtt\$$ at $H^\ell a$. Successive cells after $a$ are initial until the first right wall. Such a wall must occur before the horizontal orbit closes. Otherwise the initial phase would reach the predecessor $H^{-1}a$, whose label is the left wall, contrary to $\mathcal G$. Let $W$ be the exponent of this first right wall. The input prefix, delimiter, and required following padding cell give $W\ge\ell+2$. The segment
\[
 H^{-1}a,a,Ha,\ldots,H^W a
\]
is therefore $c_{\beta,W}$. Its vertices are distinct up to the first right wall, since an earlier repetition would close the orbit before reaching the predecessor left wall. Hence $W+2\le N$.

Let $p\le N$ be the period of $a$ under $V$. Commutation implies that $V^p$ fixes this entire segment. Its labels under the successive images $V^t$ form an admissible row sequence. Q3 gives $\mathcal G$, Q2 gives $\mathcal R$, and wall labels remain fixed in time. Boundary updates read no cell beyond the fixed walls, so this does not require the whole component to be a rectangular torus. At time $p$ the segment returns pointwise, giving the asserted return and hence halting by \Cref{ttm:lem:reset-ca}.

To obtain the area bound, write $T=T(\beta)$ and consider the first $T+2$ row segments, at times $0,\ldots,T+1$. They are pairwise disjoint as vertex sets. Indeed, each segment has a left wall, then $W$ interior vertices, then a right wall. Suppose $V^uH^ia=V^vH^ja$ with $0\le u<v\le T+1$ and $-1\le i,j\le W$. A common interior vertex has the same distance to the preceding left wall in both segments, so $i=j$; an earlier wall in one segment would otherwise be an interior vertex of the other. For a common wall, its left or right label also forces $i=j$. Commutation and invertibility now give $V^{v-u}a=a$, so the initial segment returns pointwise after $v-u\le T+1$ transitions. This contradicts the first-return time $T+2$ in the reset construction. The computation cannot cross the right wall before halting, so $W\ge S(\beta)$. Hence
\[
 N\ge(T+2)(W+2)\ge(T+2)(W_\beta+2)=A_{\rm hist}(\beta).
\]
In particular $T(\beta)\le N$.

\end{proof}

\begin{lemma}
\label{lem:two-type-closed-constructions}
For every $N\ge1$ there is an unmarked legal realization with $N$ vertices of each computation type and all 3 matchings perfect. If $\mathsf U(\beta)$ halts, there is a closed legal realization containing a marked cell and having $A_{\rm hist}(\beta)$ vertices of each computation type. Three disjoint copies supply 3 distinct marked cells of the same type in a legal realization.
\end{lemma}
\begin{proof}
For the unmarked realization take cells indexed by $\mathbb Z/N\mathbb Z$, a copied color class, $J$ the identity pairing, $R$ the cyclic successor, and $S$ the identity. Give every color vertex the background label $b_0$. Then $H$ and $V$ commute;
\eqref{ttm:eq:constant-background-update} verifies Q2 and Q3,
and Q4 is vacuous.

For a halting input, take $W=W_\beta$ and the returning row sequence of length $p=T(\beta)+2$ constructed in \Cref{ttm:lem:reset-ca}. Thus $p(W+2)=A_{\rm hist}(\beta)$. Write this sequence as $c_0\Rightarrow_W c_1\Rightarrow_W\cdots\Rightarrow_W c_p=c_0$ with $c_0=c_{\beta,W}$. Use cells indexed by $\mathbb Z/p\mathbb Z\times\mathbb Z/(W+2)\mathbb Z$ and a copied color class. Let $J$ be the identity pairing, let $R$ increment the horizontal coordinate, and let $S$ increment time. Give cell $(t,z)$ the label of position $z$ in $c_t$, and mark cell $(0,1)$.

The 2 coordinate successors commute, so Q1 holds. The chosen row sequence satisfies Q2 and Q3 away from the seams. At the horizontal seam, the right wall is followed by the left wall; wall updates ignore exterior neighbors. The returning sequence makes the temporal seam valid. Since $c_0=c_{\beta,W}$, the marked cell satisfies Q4. Every matching is perfect, and disjoint copies preserve all conditions because every forbidden positive pattern is connected in the relation multigraph.
\end{proof}

\paragraph{Additional types.}
Let $m$ be the least integer satisfying
\[
 m\equiv2\pmod4,\qquad m\ge\max\{10^{16},A_{\rm CA}+1\}.
\]
Adjoin $m-2$ additional physical types with no prescribed relations and no marks. Each vertex of an additional type has no prescribed incidence and therefore forms an isolated closed component. Use one code for the cell type and each additional type, and one code for every color of type $\mathsf B$. Thus
\[
 C_{\rm code}=m+A_{\rm CA}-1\le2m.
\]
The enlarged type multigraph has 3 indexed edges and maximum degree 3. Adding equally sized additional classes extends every closed construction without introducing a computational constraint.

\section{A three-uniform weighted construction}\label{sec:v3-three}
The matching system is represented by root links in an ordinary $3$-graph. At most one data-only edge is allowed; its weight supplies the reward for a halting computation. We first define the legal rootings, then compare this reward with the loss from incomplete matchings.

\subsection{Root codes and admissible rootings}\label{sec:compiler}
Fix the $m$ chosen above, write $C=C_{\rm code}$, and let
\[
 q=\frac{m^4-10m^3+13m^2+12m+4}{8},
 \qquad L=\lceil2048\sqrt q\rceil.
\]
The number $q$ of roots is chosen to make the maximum of the unmarked density polynomial rational. Each base relation will be represented at $L$ distinct roots; this multiplicity makes the loss from a matching defect large enough to absorb the possible gain from the marked edge.

\begin{lemma}\label{lem:block-codes}
Let $C\ge1$ and let $q\ge2C+2$ be even. There are $C$ pairwise edge-disjoint perfect matchings $C_1,\ldots,C_C$ on $[q]$ and a partition $\mathcal P$ of $[q]$ into blocks of size at most $2C+2$ such that every edge of every $C_c$ lies inside a block. For these matchings, every family of pairwise disjoint lists of root labels, each of size at least $2C+2$, has a transversal containing no edge of any $C_c$.
\end{lemma}
\begin{proof}
Let $b$ be the smallest even integer at least $C+1$, and write $q=kb+r$, where $0\le r<b$. Since $q$ and $b$ are even, $r$ is even. Partition $[q]$ into $k-1$ blocks of size $b$ and one block of size $b+r$. Every block has even size at least $C+1$ and at most $2b-2\le2C+2$. Choose $C$ distinct factors in a $1$-factorization of the complete graph on each block. For an even block of size $v$, a factorization is obtained by labeling its vertices $\{\infty\}\cup\mathbb Z_{v-1}$. The factor at $a$ consists of $\{\infty,a\}$ and $\{a+j,a-j\}$ for $1\le j\le(v-2)/2$. Define $C_c$ as the union of the $c$th chosen factors over all blocks. These are the required global perfect matchings.

Form a bipartite graph between the lists and the blocks, joining a list to every block it meets. Any $j$ lists have a union of size at least $j(2C+2)$, since the lists are disjoint. Since each block has size at most $2C+2$, their union meets at least $j$ blocks. Hall's theorem gives distinct representative blocks. Choose from each list a label in its representative block. Selected labels lie in distinct blocks, so none of the matchings $C_c$ contains a pair of selected labels.
\end{proof}

Use the complete root $3$-graph on $Q=[q]$. Choose the $C$ missing-code perfect matchings $C_c$ by \Cref{lem:block-codes}, and let $E_c=\binom Q2\setminus C_c$. The bound $q\ge2C+2$ is verified below. Code $c$ records its physical type and, at type $\mathsf B$, its automaton color. Each of the 3 base relations $e\colon S\to T$ receives $L$ distinct root labels $r_{e,\ell}$. The remaining labels are additional roots.

A partial rooting is an injection $\mu$ from a subset of $Q$ to the vertices of a $3$-graph. Its image consists of roots; the other vertices are data. For data vertex $v$ define
\[
 L_G(v)=\{\{i,j\}\colon\mu(i)\mu(j)v\in G\},\qquad
 \Cand_G(v)=\{c\colon L_G(v)\subseteq E_c\},
\]
using only present roots. A data vertex $v$ is \emph{declared} if $\Cand_G(v)=\{c\}$. Its code $c$ determines its physical type and, when applicable, its color. The vertex is \emph{ambiguous} if $|\Cand_G(v)|\ge2$. Write $V_S$ for the set of declared data vertices of physical type $S$.

The 4 edge shapes have distinct roles. Triples of roots provide the background density. Edges with 2 roots determine the candidate codes of data vertices. Edges with one root encode matching instances and supply the terms that balance the physical types. A data-only edge marks the computation and contributes the possible density gain. A $3$-graph belongs to $\mathcal C_\beta$ if some partial rooting has the following properties.
\begin{enumerate}[label=\textnormal{(H\arabic*)},leftmargin=*]
\item Every data vertex has nonempty candidates. All-root triples are
permitted.
\item At a present root $\mu(r_{e,\ell})$, where $e\colon S\to T$, a pair of declared vertices of the same physical type is permitted exactly when both have type $S$ or both have type $T$. Every pair of declared vertices of distinct types other than $\{S,T\}$ is permitted. For declared vertices $u\in V_S$ and $v\in V_T$, a present edge $\mu(r_{e,\ell})uv$ records the oriented pair $(u,v)\in D_{e,\ell}$. Additional roots permit declared pairs exactly when their types differ. A pair with an ambiguous endpoint is always permitted.
\item The union $D_e=\bigcup_\ell D_{e,\ell}$ is a partial matching
and the union relations satisfy Q1--Q3. The color of a declared $\mathsf B$-vertex is determined by its code.
\item There is at most one data-only edge, including edges with
ambiguous endpoints. Every declared endpoint of that edge has cell type $\mathsf C$ and satisfies Q4. Ambiguous endpoints impose no marked condition.
\end{enumerate}
A fully declared legal target is a $3$-graph equipped with a fixed witnessing rooting for which every data vertex is declared. All target roots and actual codes refer to that rooting. It need not use all $q$ labels.

For each relation instance $(u,v)\in D_e$ in a forbidden positive pattern, choose a copy index $\ell$ and represent the instance by the edge $\mu(r_{e,\ell})uv$. Include every choice of these indices. This also forbids matching conflicts between different copies.

\begin{lemma}
\label{lem:ternary-marker-monotone}
The class $\mathcal C_\beta$ is monotone under taking subgraphs. For a fixed partial rooting, violations of H1--H4 have positive witnesses of effectively bounded order.
\end{lemma}
\begin{proof}
Restrict a witnessing rooting to the surviving vertices. A surviving root is still a root, so every data-only edge of the subgraph was already data-only. Deleting roots or edges enlarges candidate sets. Thus a surviving declaration keeps its code; a declared vertex may instead become ambiguous. Declared binary relations and the paths and labels required to apply Q4 at declared endpoints can only disappear. The one-root rules, matching conditions, and Q1--Q3 are therefore preserved by the positive-pattern construction in \Cref{lem:two-type-positive-obstructions}. The bound of one data-only edge is preserved. Each surviving declared endpoint of that edge already had the cell type and satisfied Q4 at that endpoint. Every remaining path and label to which Q4 applies at that endpoint was already present, regardless of any other endpoint becoming ambiguous. This proves monotonicity.

Two different data-only edges witness a cardinality violation, using at most 6 vertices. An incorrectly typed declared endpoint is witnessed by the data-only edge and its declaration certificate. Such a certificate consists of at most one present root pair excluding each competing code. If the certificate and the existing root assignments are retained and the vertex remains data, every competing code is still excluded. Any admissible extension therefore retains the declared code. To witness a Q4 violation, adjoin the data-only edge, the tested endpoint's declaration certificate, and the forbidden positive colored pattern with its required declaration certificates. The existing finite bounds for Q1--Q4 and for declaration certificates therefore supply an effective bound. Expanding copy indices and allowed identifications is finite, exactly as for the unmarked constraints.
\end{proof}
\subsection{The rational baseline and weight estimates}\label{sec:numerical-weighted}
For a probability weighting $\x$ on $G$, write $P\coloneqq\lambda(G;\x)$. We compare the contribution of a marked triple with the losses from matching defects and ambiguous codes.

For the baseline model, take the unmarked realization with one data vertex of each physical type. Use all $q$ roots and include every permitted edge containing at least one root. Give each root weight $(1-t)/q$ and each data vertex weight $t/m$, so $t$ is the total data weight. Its edge polynomial is $B(t)$ below, whose 3 terms correspond to edges with 0, 1, and 2 data vertices, respectively.

\begin{lemma}
\label{lem:wi-rational-baseline}
The integer $q$ is positive and even, and
\[
 m^4/12\le q\le m^4/8,\qquad
 \sqrt q\ge m^2/4\ge32768m.
\]
Set
\[
 a=\frac{(q-1)(q-2)}{6q^2},\quad
 b=\frac{q-2}{2q},\quad c=\frac{m-1}{2m},\quad
 B(t)=a(1-t)^3+b(1-t)^2t+c(1-t)t^2.
\]
Let
\[
 \kappa=\frac m2-\frac1{m-1},\qquad
 x_0=\frac\kappa q,\qquad t_0=\frac\kappa{q+\kappa}.
\]
Then $B$ has its unique maximum at $t_0\in(0,1)$. In particular $B_*\coloneqq B(t_0)$ and $\tau_3=6B_*$ are rational, and $a<B_*<1/6$.
\end{lemma}
\begin{proof}
Since $m\equiv2\pmod4$, the numerator defining $q$ is divisible by $16$, so $q$ is an even integer. Its displayed polynomial gives $m^4/12\le q\le m^4/8$ for $m\ge30$; the lower bound follows already by retaining $m^4-10m^3$. The bounds for $\sqrt q$ follow from $m\ge10^6$. Since $C\le2m$, we have $q\ge2C+2$, as required in \Cref{lem:block-codes}. Also $L\le2049\sqrt q$ gives $3L\le6147\sqrt q<q/2$, so all relation labels fit.

For $x=t/(1-t)$, write
\[
 B(t)=\frac{a+bx+cx^2}{(1+x)^3}.
\]
Its derivative with respect to $x$ has numerator
\[
 \ell+(2/q-1/m)x-cx^2,\qquad \ell=\frac{q-2}{2q^2}>0.
\]
This quadratic has exactly one positive zero. The identity
\[
 1-2\kappa+c\kappa^2=\frac{q}{m(m-1)}
\]
shows by substitution that this zero is $x_0=\kappa/q$. Hence $t_0$ is the unique maximizer and $B_*>B(0)=a$. Finally $a<1/6$, $b<1/2$, and $c<1/2$ give
\[
 B(t)<(1-t)^3/6+(1-t)^2t/2+(1-t)t^2/2=(1-t^3)/6
\]
for $0<t<1$, proving the upper bound.
\end{proof}

For a partial matching $D\subseteq V_S\times V_T$ define
\[
 \Def(D)=\sum_{uv\in D}(x_u-x_v)^2+
 \sum_{u\in V_S\setminus\operatorname{dom}D}x_u^2+
 \sum_{v\in V_T\setminus\operatorname{im}D}x_v^2.
\]
Equivalently $\Def(D)=\sum_{V_S\cup V_T}x_v^2-2\sum_Dx_ux_v$. Thus deleting matching pairs increases the defect. In particular $\Def(D_{e,\ell})\ge\Def(D_e)$ whenever their union is a matching. Represent absent root labels by zero weights. Write $s$ for total root weight, $t=1-s$ for total data weight, and $V=\sum_i(z_i-s/q)^2$ for root variance.

\begin{lemma}
\label{lem:wi-coarse-blocks}
Fix $G\in\mathcal C_\beta$, an admissible rooting of $G$, and a probability weighting $\x$. Let $s$ be total root weight and $t=1-s$ total data weight, and represent absent roots by zero weights $z_i$. Let $V=\sum_i(z_i-s/q)^2$. Assume that $P\ge a$. Then
\[
 t<4/\sqrt q,\qquad s>1/2,\qquad V\le96q^{-3/2}.
\]
For a base relation $e$, let $Z_e$ be the total weight of its $L$ indexed roots. Let $\Delta\coloneqq\max_S\sum_{e\text{ incident with }S}Z_e$, where $S$ ranges over the physical types. These quantities also satisfy
\[
 Z_e\ge\frac{Ls}{2q},\qquad \Delta\le\frac s4,
 \qquad
 \frac12\sum_{e,j}z_{e,j}\Def(D_{e,j})
 \ge\frac{Ls}{4q}\sum_e\Def(D_e).
\]
\end{lemma}
\begin{proof}
Writing $z_i=s/q+\delta_i$ gives
\[
 e_3(\boldsymbol z)
 =as^3-(s/2-s/q)V+\frac13\sum_i\delta_i^3
 \le as^3-\frac{s}{12}V,
\]
since $\delta_i\le s(1-1/q)$ and $q\ge8$. The graph $E_c$ is complete multipartite with $q/2$ parts of size 2, so its weighted root-pair sum is at most $bs^2$. Every one-root link has value at most $t^2/2$, and the data-only edge contributes at most $t^3/27$. Hence
\[
 P-a\le\frac{t}{2q}+\frac{t^2}{q}-\frac{7t^3}{54}
                -\frac{s}{12}V.
\]
The sum of the first 3 terms is negative for $t\ge4/\sqrt q$ and $q\ge64$. For smaller $t$ their 2 positive terms sum to at most $4q^{-3/2}$. This proves the coarse bounds.

Cauchy--Schwarz now gives
\[
 Z_e\ge Ls/q-\sqrt{LV},\qquad
 \frac{\sqrt{LV}}{Ls/q}
 \le\sqrt{\frac{96\sqrt q}{Ls^2}}
 \le\sqrt{96/512}<1/2.
\]
A type has at most $3L\le8L$ indexed incidences, so
\[
 \Delta\le8Ls/q+\sqrt{8LV}
 \le\frac{16392s+1280}{\sqrt q}
 \le\frac{18952s}{\sqrt q}\le s/4.
\]
Finally apply defect monotonicity in each block and the lower bound for $Z_e$.
\end{proof}
\begin{lemma}
\label{lem:no-ddd-unmarked-weighted-estimate}
Let $G\in\mathcal C_\beta$ have an admissible partial rooting and no data-only edge. Fix a probability weighting $\x$ on $V(G)$. Represent absent roots by zero weights, and write $z_i$ for root weights, $s=\sum_i z_i$, $t=1-s$, and $V=\sum_i(z_i-s/q)^2$. Let $w$ be the total weight of ambiguous data vertices, and let $y_S$ be the total weight of declared data vertices of physical type $S$. Let $v=t-w$ and $W=\sum_S(y_S-v/m)^2$, and let $\Delta$ be the maximum total weight of indexed roots incident with a physical type. For $q\ge8$,
\begin{align*}
 P\le{}& B(t)+\frac{s}{m}\left(tw-\frac{w^2}{2}\right)
              -\frac{s^2w}{2q}\\
 &-\left(\frac{s}{12}-\frac w2\right)V
 -(s/2-\Delta)W
 -\frac12\sum_{e,j}z_{e,j}\Def(D_{e,j}).
\end{align*}
This inequality holds for every probability weighting; it does not require $P\ge a$.
\end{lemma}
\begin{proof}
The complete root contribution is at most $as^3-sV/12$. For a missing-code perfect matching $C_c$, write $z_i=s/q+\delta_i$. Its weighted pair sum satisfies
\[
 \sum_{ij\in C_c}z_iz_j
 =\frac{s^2}{2q}+\sum_{ij\in C_c}\delta_i\delta_j
 \ge\frac{s^2}{2q}-\frac V2.
\]
The graph $E_c$ is complete multipartite with $q/2$ parts of size 2, so its weighted pair sum is at most $bs^2$. For an ambiguous data vertex $u$, choose distinct candidate codes $c,d$. The matchings $C_c$ and $C_d$ are edge-disjoint, so $C_d\subseteq E_c$ and $L_G(u)\subseteq E_c\cap E_d=E_c\setminus C_d$. Thus the total 2-root contribution is at most
\[
 bs^2t-\frac{s^2w}{2q}+\frac w2V.
\]
The lower bound for the missing matching may be negative; the resulting upper bound remains valid in that case.

At an indexed root for $e\colon S\to T$, complete its permitted unmarked declared pairs. The resulting link polynomial is
\[
\begin{aligned}
 &e_2(V_S)+e_2(V_T)+\sum_{uv\in D_{e,j}}x_ux_v
       +\sum_{\{A,B\}\ne\{S,T\}}y_Ay_B\\
 &\qquad=e_2(\boldsymbol y)+\frac12(y_S-y_T)^2
                  -\frac12\Def(D_{e,j}).
\end{aligned}
\]
Here $e_2(V_S)$ is the sum of products on unordered pairs of vertices in $V_S$; the formula follows by expanding the defect. An additional root contributes at most $e_2(\boldsymbol y)$. Applying $\tfrac12(\eta_S-\eta_T)^2\le\eta_S^2+\eta_T^2$, with $\eta_S=y_S-v/m$, and summing the incidence weights gives
\[
 \frac12\sum_{e,j}z_{e,j}(y_S-y_T)^2\le\Delta W.
\]
Since $e_2(\boldsymbol y)=cv^2-W/2$, summing these root links bounds the declared one-root contribution by
\[
 csv^2-(s/2-\Delta)W
       -\frac12\sum_{e,j}z_{e,j}\Def(D_{e,j}).
\]

At a fixed root, the weight of pairs with an ambiguous endpoint is at most $vw+w^2/2=tw-w^2/2$. Summing root weights gives at most $s(tw-w^2/2)$. Finally,
\[
 c(t-w)^2+tw-w^2/2
 =ct^2+\frac1m(tw-w^2/2).
\]
Adding these bounds proves the assertion, including absent roots and zero data weights. There is no data-only contribution.
\end{proof}
\begin{lemma}
\label{lem:ternary-ambiguity-short}
For the fixed parameters of \Cref{lem:wi-rational-baseline}, let $s+t=1$, $s,t\ge0$, and $0\le w\le t$. Then
\[
 B(t)+\frac{s}{m}(tw-w^2/2)-\frac{s^2w}{2q}
       +\frac{w^3}{27}\le B_*.
\]
The same holds after deleting the nonnegative last term. More precisely, the gap in the displayed inequality is at least
\[
 \frac{s}{2m}(t-w-x_0s)^2+\frac{7x_0^2}{27}s^2w.
\]
Equality in the displayed upper bound holds exactly when $w=0$ and $t=t_0$.
\end{lemma}
\begin{proof}
We first factor $B_*-B(t)$, then incorporate the terms involving $w$. Let
\[
 d=3B_*-c+2B_*x_0,\qquad
 \varepsilon=\frac1{2q}-\frac{x_0}{m}
             =\frac1{m(m-1)q}.
\]
Since $B(t_0)=B_*$, $B'(t_0)=0$, and $t_0=x_0/(1+x_0)$, comparing coefficients gives the factorization
\[
 B_*-B(t)=(t-x_0s)^2(B_*t+ds),\qquad s+t=1.
\]
Consequently the difference between $B_*$ and the claimed left side is
\[
 \frac{s}{2m}(t-w-x_0s)^2
 +\left[B_*t+s\left(d-\frac1{2m}\right)\right](t-x_0s)^2
 +\varepsilon s^2w-\frac{w^3}{27}.
\]

To bound this expression from below, we first estimate its coefficients. Since $q\ge m^4/12$ and $\kappa<m/2$, we have $3q>m(m-1)\kappa^2$, so $\varepsilon\ge x_0^2/3$. Also $B_*\ge a\ge1/8$ and $\kappa\ge6$ give
\[
 d-\frac1{2m}\ge3a-\frac12+\frac{x_0}{4}
 =\frac{\kappa/4-3/2}{q}+\frac1{q^2}\ge0.
\]
Finally, $0\le w\le t$ gives
\[
 w^3\le wt^2
 \le2t(t-x_0s)^2+2x_0^2s^2w.
\]
Since $B_*\ge1/8>2/27$ and $\varepsilon\ge x_0^2/3>2x_0^2/27$, every term left after this substitution is nonnegative, and retaining the first square and $(\varepsilon-2x_0^2/27)s^2w$ gives the quantitative bound. If equality holds, the positive coefficient of $t(t-x_0s)^2$ forces $t=x_0s$ whenever $t>0$; $t=0$ is impossible since $B_*>a$. Thus $s>0$, and the positive coefficient of $s^2w$ forces $w=0$. Conversely $w=0$, $t=t_0$ gives equality.
\end{proof}

\begin{lemma}
\label{lem:ternary-grounded-eight}
Let a finite relation component have nonnegative vertex weights, total weight $\mu$, and at least one missing incidence. Let $D_C$ be its matching energy, namely the sum of squared differences along relation edges, counted separately for each matching, and squared weights at missing incidences. Every vertex of weight $a$ satisfies
\[
 a^3\le\frac94\mu D_C.
\]
\end{lemma}
\begin{proof}
Take a simple path from this vertex to a missing incidence, append a formal zero-weight leaf, and write its weights $a_0=a,a_1,\ldots,a_N=0$. For $u,v\ge0$, Cauchy--Schwarz applied to $u^{3/2}-v^{3/2}=\tfrac32\int_v^u\sqrt z\,dz$ gives
\[
 (u^{3/2}-v^{3/2})^2\le\frac98(u+v)(u-v)^2.
\]
Telescoping and Cauchy--Schwarz therefore give
\[
 a^3\le\frac98
 \left(\sum_{i=0}^{N-1}(a_i+a_{i+1})\right)
 \left(\sum_{i=0}^{N-1}(a_i-a_{i+1})^2\right)
 \le\frac94\mu D_C.
\]
Each actual vertex occurs at most once. The final term $a_{N-1}^2$ is included in $D_C$ at the chosen missing incidence.
\end{proof}
\begin{proposition}
\label{prop:ternary-marker-nonhalting-bound}
Let $G\in\mathcal C_\beta$ have an admissible rooting. If every declared vertex belonging to a data-only edge lies in a component with a missing incidence, then every probability weighting satisfies $P\le B_*$. In particular, this holds for every $G\in\mathcal C_\beta$ when $\mathsf U(\beta)$ does not halt. Moreover, if $\lambda(G)>B_*$, then $A_{\rm hist}(\beta)\le v(G)$.
\end{proposition}
\begin{proof}
If $P<a$, there is nothing to prove. Otherwise
\Cref{lem:wi-coarse-blocks} gives $t<4/\sqrt q$, $s>1/2$,
$\Delta\le s/4$, and
\[
 \frac12\sum_{e,\ell}z_{e,\ell}\Def(D_{e,\ell})\ge\alpha D,
 \qquad \alpha=Ls/(4q),\quad D=\sum_e\Def(D_e).
\]
Let $w$ be ambiguous weight, $W$ the variance of the total declared weight in each physical type, and $M$ the data-only edge weight, or zero if there is no such edge. Apply \Cref{lem:no-ddd-unmarked-weighted-estimate} after deleting that edge, and then restore its weight. This gives
\begin{align*}
 P\le{}&B(t)+\frac{s}{m}(tw-w^2/2)-\frac{s^2w}{2q}+M\\
 &-\left(\frac{s}{12}-\frac w2\right)V
 -\frac{s}{4}W-\alpha D.
\end{align*}
The coarse bounds imply $w\le t<s/6$.

Let $\eta$ be the largest weight of a declared endpoint of the data-only edge, or zero if none exists. Every such endpoint is a marked cell in a component with a missing incidence, by the hypothesis. The estimate in \Cref{lem:ternary-grounded-eight} therefore gives $\eta^3\le\tfrac94tD$. We claim
\[
 M\le\frac{w^3}{27}+2\eta^3.
\]
For 3 ambiguous endpoints this is AM--GM. For one ambiguous endpoint of weight $c$, AM--GM gives $M\le c\eta^2\le c^3/27+2\eta^3$. For 2 ambiguous endpoints their product is at most $w^2/4$, and weighted AM--GM gives
\[
 M\le\frac{w^2\eta}{4}
 \le\frac{w^3}{27}+\frac{27\eta^3}{16}
 \le\frac{w^3}{27}+2\eta^3.
\]
With 3 declared endpoints, $M\le\eta^3$; if the edge is absent, $M=0$. Consequently $M\le w^3/27+\tfrac92tD$. Since
\[
 \alpha=Ls/(4q)\ge256/\sqrt q>\tfrac92t,
\]
the matching-defect penalty absorbs the second term. Discarding the nonpositive variance penalties now leaves exactly the expression in
\Cref{lem:ternary-ambiguity-short}, proving $P\le B_*$.

On a nonhalting input, \Cref{lem:two-type-closed-recovery} supplies the hypothesis. If instead $\lambda(G)>B_*$, a maximizing weighting violates the upper bound, so some declared endpoint lies in a closed component. The same recovery lemma gives $A_{\rm hist}(\beta)\le v(G)$.
\end{proof}
\begin{proposition}
\label{prop:ternary-marker-lower-witnesses}
The unmarked closed construction attains $B_*$ on every input. If an input has a legal closed marked two-type realization with $N$ vertices in each of its 2 computation types, there is a fully declared legal target for the three-uniform construction with a weighting of value
\[
 B_*+\left(\frac{t_0}{3mN}\right)^3>B_*.
\]
\end{proposition}
\begin{proof}
For an unmarked closed two-type realization with equal size in each of its 2 computation types, supply the same number of vertices in every additional physical type. Give all roots weight $(1-t_0)/q$ and give each physical type total weight $t_0/m$, distributed equally among its vertices. Include all root triples and, for every data vertex $v$ of prescribed code $c$, all edges $\mu(i)\mu(j)v$ with $\{i,j\}\in E_c$. For each base relation $e$, use its prescribed matching as $D_{e,\ell}$ at every copy root $r_{e,\ell}$. Include every other one-root pair permitted by H2, including the permitted same-type pairs.

This construction declares every data vertex with its prescribed code. Indeed, all $q$ roots are present, and if $v$ has prescribed code $c$ and $d\ne c$, a pair in $C_d\subseteq E_c$ supplies a present edge that excludes candidate $d$. The link identity in \Cref{lem:no-ddd-unmarked-weighted-estimate} now gives $P=B(t_0)=B_*$.

For a closed legal two-type realization containing a marked cell, take 3 disjoint copies. There are $3N$ vertices in each computation type; add $3N$ vertices in each additional physical type. Use the same completion and assign each data vertex weight $t_0/(3mN)$. Choose the marked cell in each of the 3 copies and add their triple as the sole data-only edge. All 3 marked cells satisfy Q4 in their respective components, so the new data-only edge satisfies H4. Root-pair candidates are unchanged, and the 3 marker vertices are distinct. The total contribution of edges containing a root remains $B_*$, and the new data-only edge contributes $(t_0/(3mN))^3$.
\end{proof}
\begin{corollary}\label{thm:weighted-compiler}
There is a finite fully declared unmarked target $T_0$, independent of the input, with $\lambda(T_0)=B_*$. On a nonhalting input, $\Lambda(\mathcal C_\beta)=B_*$. On a halting input, there is a fully declared target $T_\beta\in\mathcal C_\beta$ with $q+3mA_{\rm hist}(\beta)$ vertices and a rational weighting of value
\[
 B_*+\left(\frac{t_0}{3mA_{\rm hist}(\beta)}\right)^3.
\]
The number $\tau_3\coloneqq6B_*$ lies in $\mathbb Q\cap(0,1)$.
\end{corollary}

\begin{proof}
Use the one-vertex-per-type unmarked realization for $T_0$. Its specified weighting attains $B_*$, and its legality on the fixed nonhalting word $\mathtt0$ gives the reverse Lagrangian bound. The 2 preceding propositions and the closed realizations in
\Cref{lem:two-type-closed-constructions} give all remaining
claims. Rationality and the range follow from \Cref{lem:wi-rational-baseline}.
\end{proof}

\subsection{An effective finite forbidden basis}\label{sec:direct-finite-root-removal}

\begin{proposition}\label{prop:direct-finite-basis}
For every $\beta$, the class $\mathcal C_\beta$ has an effective finite forbidden basis.
\end{proposition}
\begin{proof}
Membership in $\mathcal C_\beta$ is decidable by testing all partial rootings. Monotonicity follows from \Cref{lem:ternary-marker-monotone}. All violations for a fixed rooting have bounded positive witnesses. Include every currently assigned root vertex in the witness vertex set, adding at most $q$ vertices, and take the induced subhypergraph on this enlarged set. This preserves the link of every retained data vertex on the assigned roots, and hence its candidate set. The finite positive computational list and the declaration certificates give an effective total witness bound $b(\beta)\le b_0(|\beta|+1)$ for an effective fixed integer $b_0$.

To check persistence, keep all current root assignments and add roots outside the witness. Its old data vertices remain data, and their candidates can only shrink. An empty candidate set stays empty; a singleton stays the same or becomes empty, which is itself a violation. Thus, if all candidates remain nonempty, every declaration used in the witness is unchanged. Its positive relation instances, marked forbidden patterns, and violations of the type requirements in H2 and H4 persist. Two old data-only edges remain data-only. Adding roots outside the witness therefore cannot repair it.

Every edgeless $3$-graph is admissible with the empty rooting. Applying \Cref{lem:unroot} therefore gives an effective finite forbidden basis whose members have order at most $b(\beta)\sum_{j=0}^q(qb(\beta))^j$.
\end{proof}

\section{A single forbidden three-uniform hypergraph}\label{sec:single}
We now convert the weighted construction into a single forbidden $3$-graph. For a finite $3$-graph $A$, its \emph{canonical extension} $E(A)$ adds one new vertex $z_{xy}$ and the edge $xyz_{xy}$ for every pair $\{x,y\}$ of vertices of $A$ that lies in no edge of $A$. We call $A$ the \emph{core} of this extension. Every core pair is covered by an edge of $E(A)$.

The \emph{categorical product} $\prod_i G_i$ has vertex set $\prod_i V(G_i)$; a triple is an edge when its projection in every coordinate is an edge of that factor. Each coordinate projection is a homomorphism. We will take such a product over a finite forbidden basis. Its projections supply an upper bound on its Tur\'an density. To obtain the lower bound, we first arrange that every forbidden core contains enough disjoint cliques to select compatible maps on one factor.

For $v\ge3$, write $p(v)\coloneqq1-3/v+2/v^2=6\lambda(K_v^{(3)})$.

\subsection{A clique-packing guard}\label{sec:safe-guard}
We will choose $H$ so that adjoining $\Forb(H)$ to $\mathcal C_\beta$ leaves its Lagrangian supremum unchanged. Every forbidden hypergraph for the enlarged class must then contain $H$.

\begin{lemma}\label{lem:safe-simple-packing}
Let $t\ge4$ and $q\ge10^6t^3$, and let
\[
 s=\left\lfloor\frac{24q}{25t}\right\rfloor,\qquad
 k=t(s-1),\qquad n=st,\qquad H=sK_t^{(3)}.
\]
Then $s\ge2$, $0.95q<k<n\le0.96q$, and
\[
 6\Lambda(\Forb(H))\le1-\frac3q<p(q).
\]
\end{lemma}
\begin{proof}
The floor gives $0.96q-2t<k<n\le0.96q$. Since $q>200t$, this implies $0.95q<k$ and $k\ge100$. Every $t$-set in a $K_t^{(3)}$-free $3$-graph contains a missing edge; double counting gives $\pi(K_t^{(3)})\le1-\delta$, where $\delta=\binom t3^{-1}>6/t^3$. A clique is pair-covered, so its avoidance class is closed under blowups. Rational weighted blowups therefore give $6\Lambda(\Forb(K_t^{(3)}))\le1-\delta$.

Fix an $H$-free $3$-graph $G$ and a probability weighting $\x$, and write $P\coloneqq\lambda(G;\x)$. Let $S$ be the union of a maximal collection of vertex-disjoint copies of $K_t^{(3)}$. Then $|S|\le k$, and the subhypergraph induced by $V(G)\setminus S$ is $K_t^{(3)}$-free. Let $x$ be the weight on $S$ and $u=1-x$. The contributions with 3, 2, one, and zero vertices in $S$ are bounded by
\[
 6P\le p(k)x^3+3(1-1/k)x^2u+3xu^2+(1-\delta)u^3
      =1-\frac{3x^2}{k}+\frac{2x^3}{k^2}-\delta u^3.
\]
If $u\ge1/100$, then $\delta u^3>6/q$, while $3x^2/k-2x^3/k^2\ge0$. If $u<1/100$, then $x>0.99$ and
\[
 \frac{3x^2}{k}-\frac{2x^3}{k^2}
 \ge\frac{(3-2/k)x^2}{k}
 >\frac{2.9}{k}\ge\frac{2.9}{0.96q}>\frac3q.
\]
Thus every weighting satisfies $6P<1-3/q$, and taking suprema proves the assertion.
\end{proof}

\subsection{Selecting compatible maps}\label{sec:selector}

Starting from a homomorphism of the whole product, we seek a fixed number of maps on one factor. No selected map may send 2 core vertices to the same target root. For each edge, taking its 3 vertex images under 3 distinct selected maps must still give a target edge. The number of maps will be independent of the number of factors. Repeated root images along coordinate lines lead to intersecting families, which the following estimates control.

Two words are intersecting if they agree in some coordinate. A family is intersecting if every 2 of its words intersect; 2 families are cross-intersecting if every word in one intersects every word in the other. In this section $K_n$ denotes the ordinary complete graph.

For $A\subseteq[n]^a$, let $\mu(A)\coloneqq |A|/n^a$. On $[n]^a$ use uniform probability measure, and let $P$ be the normalized adjacency operator of the categorical product $K_n^a$. Thus $(Pf)(x)$ averages $f$ over the points differing from $x$ in every coordinate. On one coordinate the operator has eigenvalue $1$ on constants and $-1/(n-1)$ on mean-zero functions. Decomposing each coordinate space into constants and mean-zero functions gives an orthogonal tensor decomposition. Level $k$ is the sum of the tensor subspaces with mean-zero factors in exactly $k$ coordinates and constant factors in the others. The product operator has eigenvalue $(-1/(n-1))^k$ on this level.

\begin{lemma}\label{lem:spectral-cross}
Let $a\ge1$, $n\ge2$, and let $A,B\subseteq[n]^a$ be cross-intersecting. Writing $\alpha=\mu(A)>0$ and $\beta=\mu(B)$, one has
\[
 \alpha\beta\le\frac{(1-\alpha)(1-\beta)}{(n-1)^2},
 \qquad
 \beta\le\frac1{(n-1)^2\alpha}.
\]
\end{lemma}
\begin{proof}
There are no product-graph edges between $A$ and $B$, so $\langle 1_A,P1_B\rangle=0$. On mean-zero functions the operator norm of $P$ is $1/(n-1)$. Hence
\[
 \alpha\beta
 =|\langle1_A-\alpha,P(1_B-\beta)\rangle|
 \le\frac{\sqrt{\alpha(1-\alpha)\beta(1-\beta)}}{n-1}.
\]
Squaring and cancelling proves the first inequality when $\beta>0$; the zero case is immediate. The second follows.
\end{proof}

A \emph{coordinate slice} is a set $J=\{x\in[n]^a\colon x_i=b\}$. The section of a family $F$ on $J$ is the subset of $[n]^{a-1}$ obtained from $F\cap J$ by deleting coordinate $i$.

\begin{lemma}
\label{lem:spectral-fourth}\label{cor:spectral-slice}
Let $n\ge2^{20}$, let $a\ge1$, and let $\theta=1/(64n)$. Every intersecting $F\subseteq[n]^a$ either has measure at most $\theta$, or admits a coordinate slice $J$ such that $\mu(F\setminus J)\le\theta$. In the latter case, if $a\ge2$, every family in $[n]^{a-1}$ cross-intersecting with the section of $F$ on this slice has measure at most $\theta$.
\end{lemma}
\begin{proof}
Suppose $\alpha\coloneqq\mu(F)>\theta$. If $a=1$, then $F$ is a singleton and the assertion is immediate, so assume $a\ge2$. Let $f=1_F$, and write its level-one projection as
\[
 g=\sum_{i=1}^a g_i(x_i),\qquad
 g_i=\mathbb E[f\mid x_i]-\alpha,\qquad W=\mathbb E[g^2].
\]
Since $F$ is independent, $\langle f,Pf\rangle=0$. The constant and level-one contributions are $\alpha^2$ and $-W/(n-1)$. The remaining eigenspaces have eigenvalues of absolute value at most $(n-1)^{-2}$, and the squared norms of the projections of $f$ onto them sum to at most $\mathbb E[f^2]=\alpha$. Hence
\[
 0\ge\alpha^2-\frac{W}{n-1}-\frac{\alpha}{(n-1)^2},
 \qquad
 W\ge(n-1)\alpha^2-\frac{\alpha}{n-1}
   \ge\frac{n\alpha^2}{4}.
\]
The last inequality follows from $\alpha(n-1)^2\ge n^2\alpha/4>n/256\ge4096$.

Let $M=\max_{i,b}\mu(F\mid x_i=b)$. Orthogonality gives $\mathbb E[fg]=W$, and $|g_i|\le M$. H\"older's inequality and independence of the mean-zero $g_i$ therefore imply
\[
 \frac{W^4}{\alpha^3}\le\mathbb E[g^4]
 \le M^2W+3W^2.
\]
Since $W^2/\alpha^3\ge n^2\alpha/16>n/1024\ge1024$,
\[
 M^2\ge\frac{W^3}{\alpha^3}-3W
 \ge\frac{W^3}{2\alpha^3}
 \ge\frac{(n\alpha)^3}{128}>2^{-25}.
\]
In particular $M>2^{-13}$. Choose a maximizing slice $J=\{x_i=b\}$. By \Cref{lem:spectral-cross}, every family cross-intersecting with its section has measure at most
\[
 \frac1{(n-1)^2M}
 <\frac{2^{13}}{(n-1)^2}
 \le\frac{2^{14}}{n^2}
 \le\frac1{64n}=\theta.
\]
This applies in particular to every section at $x_i=c\ne b$; averaging yields $\mu(F\setminus J)\le\theta$.
\end{proof}

\begin{lemma}
\label{lem:guard-clique-normalized-sampling}
Let $h\ge3$, $t\ge4h$, and $H=sK_t^{(3)}$, where $0.95q<n=st\le q$. In each of finitely many copies of $H$, mark $a_i\le q/2$ vertices, with $\sum_i a_i\le q$. There exists a probability distribution on ordered $h$-tuples in each copy with the following properties.
\begin{enumerate}
\item the $h$ entries are distinct unmarked vertices of a single
constituent clique;
\item when the tuples are sampled independently in the different copies, for each $k\in[h]$ the vector formed by the $k$th selected vertex in every copy has a product distribution whose density relative to uniform product measure on the copies of $V(H)$ is at most $12$.
\end{enumerate}
\end{lemma}
\begin{proof}
In coordinate $i$, delete the marked vertices. Discard every constituent clique with fewer than $h$ remaining vertices. If such a discarded clique has $r<h$ remaining vertices, then $t-r$ vertices were marked, and
\[
 r\le\frac{t-r}{3},
\]
since $t\ge4h>4r$. Therefore the discarded unmarked vertices number at most $a_i/3$. The union $Y_i$ of the remaining unmarked constituent cliques has size
\[
 N_i\coloneqq|Y_i|\ge n-\frac43a_i>0.
\]
Choose a remaining constituent clique with probability equal to its number of remaining vertices divided by $N_i$, and then choose a uniformly random ordered $h$-tuple of distinct vertices from it. Every position of the chosen tuple is uniform on $Y_i$. This proves the first assertion.

The vector formed by any fixed position across the copies has density at most $\prod_i n/N_i$ relative to uniform product measure. The function
\[
 F(u)=-\log\left(1-\frac{4u}{3n}\right)
\]
is convex on $[0,q/2]$ and vanishes at zero, so $F(a_i)\le(2a_i/q)F(q/2)$. Hence
\[
 \prod_i\frac n{N_i}
 \le\exp\left(\sum_iF(a_i)\right)
 \le\exp\bigl(2F(q/2)\bigr)
 =\left(\frac{3n}{3n-2q}\right)^2
 <\left(\frac{57}{17}\right)^2<12.
\]
This proves the asserted density bound.
\end{proof}
\begin{proposition}
\label{prop:elementary-profiles}
Let $h\ge3$, $q\ge2^{24}$, and let $H=sK_t^{(3)}$ have $t\ge8h$ and $0.95q<n=st\le q$. Suppose every core $A_1,\ldots,A_d$, $d\ge2$, contains $H$. For any homomorphism
\[
 \phi\colon\prod_{i=1}^dE(A_i)\longrightarrow T
\]
and any designated set of at most $q$ target roots, some direction $j$ has maps $f_1,\ldots,f_h\colon V(E(A_j))\to V(T)$ such that no $f_a$ repeats a root on $V(A_j)$, and for every ordered edge $(x,y,z)$ of $E(A_j)$ and distinct $a,b,c$,
\[
 \{f_a(x),f_b(y),f_c(z)\}\in T.
\]
\end{proposition}
The number $h$ of selected maps does not depend on the number $d$ of product factors.
\begin{proof}
We first select a coordinate in which, after excluding certain coordinate slices, each target root is repeated on only a small proportion of lines. We then sample tuples in the remaining coordinates to obtain $h$ maps with no repeated root images and the required edge condition.

Fix a copy $X_i$ of $H$ in each core. Two vertices of $\prod_i V(A_i)$ that differ in every coordinate belong together to a product edge, by completing their pair in every canonical factor. Thus their images under $\phi$ are distinct.

For a coordinate $j$ and a tuple $\boldsymbol w\in\prod_{i\ne j}V(X_i)$, write $\phi(x,\boldsymbol w)$ for the image of the product vertex obtained by inserting $x$ in coordinate $j$. A line in direction $j$ is obtained by fixing $\boldsymbol w$ and letting $x$ range over $V(A_j)$. For a target root $r$, define
\[
 Z_j(r)\coloneqq\left\{\boldsymbol w\in\prod_{i\ne j}V(X_i)\colon
 \bigl|\{x\in V(A_j)\colon\phi(x,\boldsymbol w)=r\}\bigr|\ge2\right\}.
\]
Thus $Z_j(r)$ consists of the tuples whose lines have at least 2 vertices with image $r$. Each $Z_j(r)$ is intersecting. Indeed, suppose that $\boldsymbol w,\boldsymbol w'\in Z_j(r)$ differ in every coordinate. Each of their lines has at least 2 vertices with image $r$, so we can choose distinct $x,x'\in V(A_j)$ with $\phi(x,\boldsymbol w)=\phi(x',\boldsymbol w')=r$. The corresponding vertices of $\prod_iV(A_i)$ differ in every coordinate, contradicting the distinctness of their images proved above.

For $j\ne k$, project $Z_j(r)$ and $Z_k(r)$ to their common coordinates. These projections are cross-intersecting. Otherwise choose a bad $j$-line and bad $k$-line that differ in every common coordinate. On the former choose an occurrence of $r$ differing from the latter's fixed $j$-entry; on the latter choose one differing from the former's fixed $k$-entry. The resulting 2 vertices of $\prod_iV(A_i)$ differ in every coordinate, a contradiction. When $d=2$ this means one of the 2 bad-line sets is empty.

Let $\theta=1/(64n)$. \Cref{cor:spectral-slice} applies since $n>0.95q>2^{20}$. Start in direction $1$ and assign a slice to each $Z_1(r)$ of measure greater than $\theta$. If no coordinate other than $1$ receives more than $q/2$ assignments, keep this direction. Otherwise let $k$ receive more than $q/2$, and let $R_0$ be the roots assigned slices in coordinate $k$.

For each $r\in R_0$, the assigned section is cross-intersecting with the projection of $Z_k(r)$ to the common coordinates of directions $1$ and $k$. If $d\ge3$, \Cref{cor:spectral-slice} shows that the projection has measure at most $\theta$, and hence $\mu(Z_k(r))\le\theta$. If $d=2$, then $Z_k(r)$ is empty. We therefore switch to direction $k$ and discard the old assignments. For each root $r\notin R_0$ with $\mu(Z_k(r))>\theta$, choose a slice supplied by \Cref{cor:spectral-slice}. Since $|R_0|>q/2$, fewer than $q/2$ assignments are needed in total, so every coordinate satisfies the required bound.

Fix the resulting direction $j$ and slice assignments. There are at most $q$ assignments in total and at most $q/2$ in each coordinate other than $j$. Every root has either a bad set of measure at most $\theta$, or an assigned slice outside which its bad set has measure at most $\theta$. Mark the vertices defining the slices. Apply \Cref{lem:guard-clique-normalized-sampling} with $2h$ in place of $h$ in each coordinate other than $j$, using $t\ge8h$. For each position $a\in[2h]$, let $w_a$ be the vector of the $a$th sampled vertices in those coordinates. Every $w_a$ avoids the assigned slices and has marginal density at most $12$ relative to uniform product measure. A union bound over roots gives, for each $a$,
\[
 \Pr\bigl[w_a\in\textstyle\bigcup_r Z_j(r)\bigr]
 \le12q\theta<\frac{12}{64(0.95)}<\frac14.
\]
The expected number of sampled vectors in $\bigcup_r Z_j(r)$ is less than $h/2$. Choose a sample with fewer than $h$ vectors in this union and retain any $h$ outside it, relabeling them $w_1,\ldots,w_h$. Only linearity of expectation is used; the vectors $w_a$ need not be independent. Let $f_a(x)\coloneqq\phi(x,w_a)$. None repeats a root on $V(A_j)$. Every 3 distinct retained tuples form an edge, because their entries in each coordinate are distinct vertices of one guard clique. This proves the universal ordered-edge condition.
\end{proof}

\subsection{Reconstruction by unique completion}
To reconstruct a legal rooting, we compare the images $f_i(x)$ of a fixed source vertex $x$ as the map index $i$ varies. In both a partial matching and the data-only edge, when present, fixing all but one entry leaves at most one possible completion. The next lemma will force these images to agree on suitable sets of indices.

\begin{lemma}\label{lem:transversal-rigidity}
Let $r\ge2$ and let $R\subseteq X_1\times\cdots\times X_r$ be a relation in which fixing any $r-1$ entries determines at most one value of the remaining entry. For each $j$, let $S_j$ be an index set of size at least $r+1$ and let $a_j\colon S_j\to X_j$. If
\[
 (a_1(i_1),\ldots,a_r(i_r))\in R
 \quad\text{whenever the indices }i_1,\ldots,i_r\text{ are distinct},
\]
then each $a_j$ is constant, and the tuple of constants belongs to $R$.
\end{lemma}
\begin{proof}
Fix $j$ and 2 indices in $S_j$. In each of the other $r-1$ index sets, choose an index avoiding both and all previously chosen indices. At most $r$ choices are excluded in any index set, so this is possible. The resulting 2 tuples in $R$ agree in all but entry $j$; uniqueness makes that entry equal as well. Thus each $a_j$ is constant. Choosing any transversal of the index sets proves the final assertion.
\end{proof}

\begin{proposition}
\label{thm:many-profile-ternary-reconstruction}
Let $C$ be the number of codes, let $T$ be a fully declared legal target in $\mathcal C_\beta$, and let $A$ be a finite core. Suppose $C\ge2$, $h\ge6C$, and maps
\[
 f_1,\ldots,f_h\colon V(E(A))\longrightarrow V(T)
\]
satisfy the following conditions.
\begin{enumerate}
\item No $f_i$ repeats a target root on $V(A)$.
\item Every ordered edge $(x,y,z)$ of $E(A)$ satisfies
\[
 \{f_i(x),f_j(y),f_k(z)\}\in T
 \qquad(i,j,k\text{ distinct}).
\]
\end{enumerate}
Then $A\in\mathcal C_\beta$.
\end{proposition}
\begin{proof}
Arrange the values $f_1(x),\ldots,f_h(x)$ as the row indexed by $x\in V(A)$. Pair coverage of the core gives
\begin{equation}
 f_i(x)\ne f_j(y)\qquad(x\ne y,\ i\ne j).
 \label{eq:ternary-cross-privacy}
\end{equation}
Together with condition (i), which excludes repeated root images within each map, this shows that the sets of root labels occurring in different core rows are disjoint.

We construct a partial rooting of $A$ in 3 steps.
\begin{enumerate}
\item If a root label occurs at least twice in the row of $x$, assign $x$ one such root label.
\item Among the remaining vertices, keep $x$ as data if at least 4 entries of its row are target data vertices with a common code $c_x$. Choose one such code and let $S_x$ be the set of all indices $i$ for which $f_i(x)$ is a data vertex of code $c_x$. Thus $|S_x|\ge4$.
\item Each remaining row has at most $3C$ data entries and no repeated root label, so it contains at least $h-3C\ge2C+2$ distinct root labels. By \Cref{lem:block-codes}, choose one label for each such row, with the chosen labels in distinct blocks, and assign these labels to the corresponding vertices.
\end{enumerate}
All selected root labels are distinct, because different rows have disjoint sets of root labels. Pairs of labels selected in step \textup{(iii)} are permitted by every code.

We first show that $c_x\in\Cand_A(x)$ for every data vertex $x$. Indeed, consider a source edge with data row $x$ and 2 root rows. If both root labels were selected in step \textup{(iii)}, their pair is permitted by every code. Otherwise, at least one selected label was chosen in step \textup{(i)} and occurs at least twice. We can therefore choose distinct occurrence indices for the two labels and then an index in $S_x$ avoiding both. The resulting target edge shows that the selected root pair belongs to $E_{c_x}$. Thus candidates are nonempty, and any declared source row has code $c_x$. This proves H1.

For H2, consider a source edge with one root and 2 declared data vertices. Choose an occurrence of the root label and then occurrences of the 2 codes at distinct indices avoiding it. The resulting target edge shows that the source edge satisfies the one-root rule H2. Ambiguous endpoints impose no such rule, so H2 holds.

It remains to recover the constrained data rows. For a declared source instance of an indexed relation, fix an index $k$ realizing its selected root label. All cross pairs on $S_x\setminus\{k\}$ and $S_y\setminus\{k\}$ lie in the target partial matching. These sets have size at least 3, so \Cref{lem:transversal-rigidity} with $r=2$ makes both rows constant there. For a source data-only edge, all cross triples on its 3 chosen occurrence sets are data-only edges of $T$. Their ordered triples form a relation with unique completion, since $T$ has at most one such edge. \Cref{lem:transversal-rigidity} with $r=3$ makes all 3 rows constant on their full occurrence sets.

These constants define one map $g$ on every row incident with a declared relation or a data-only edge. They agree whenever 2 constraints meet, since each relevant occurrence set is $S_x$ with at most one index removed, so their intersection is nonempty. Every value of $g$ occurs at least 3 times in its row. Equation~\eqref{eq:ternary-cross-privacy} therefore makes $g$ injective. By construction it preserves the chosen codes, every indexed relation, and every data-only edge.

All matching conditions and Q1--Q3 pull back through this injection, because each forbidden positive pattern uses only the preserved relations and colors, and equality conclusions pull back by injectivity. Every source data-only edge maps to the unique target data-only edge; injectivity therefore allows at most one source edge of this kind, including those with ambiguous endpoints. Its declared endpoints have cell code and inherit Q4. The marked cell and every vertex of a forbidden Q4 pattern lie in the domain of $g$, so such a source pattern would give the same forbidden marked pattern in the target. Ambiguous endpoints impose no marked condition. Thus H3 and H4 also hold, completing the reconstruction.
\end{proof}

\subsection{The forbidden product}\label{sec:final-singleton-reduction}

We now prove \Cref{thm:singleton}. The root-code parameters have already been fixed independently of $\beta$. We choose the guard parameters by
\[
 h=6C,\qquad t=8h,\qquad
 s=\left\lfloor\frac{24q}{25t}\right\rfloor,
 \qquad H=sK_t^{(3)}.
\]
The choice $m\ge10^{16}$ gives $h\le12m$ and
\[
 \frac q{h^3}\ge\frac m{12\cdot12^3}>2^{32},
 \qquad q\ge2^{32}h^3>10^6(8h)^3=10^6t^3.
\]
Consequently \Cref{lem:safe-simple-packing} and
\Cref{prop:elementary-profiles} apply, and
\begin{equation}\label{eq:final-singleton-safe-guard}
 6\Lambda(\Forb(H))<p(q)<6B_*=\tau_3.
\end{equation}

For each input $\beta$, define
\[
 \mathcal C^\dagger_\beta
 \coloneqq\mathcal C_\beta\cup\Forb(H).
\]
Let $\mathcal B^\dagger_\beta$ be the finite forbidden core family, and define the ordinary finite $3$-graph
\begin{equation}
 P_\beta
 \coloneqq\prod_{B\in\mathcal B^\dagger_\beta}E(B),
 \label{eq:final-singleton-product}
\end{equation}
where the product is the categorical product of simple $3$-graphs. The next lemma verifies that the basis and product are effective.

\begin{lemma}
\label{lem:final-singleton-effective-basis}
The class $\mathcal C^\dagger_\beta$ has an effectively computable, nonempty finite forbidden core family. Every member $B$ of that basis contains $H$ as an ordinary subgraph and satisfies $B\notin\mathcal C_\beta$. Moreover
\[
 \Lambda(\mathcal C^\dagger_\beta)
 =\Lambda(\mathcal C_\beta).
\]
\end{lemma}
\begin{proof}
\Cref{prop:direct-finite-basis} supplies an
effective finite family $\widehat{\mathcal A}_\beta$ with $\mathcal C_\beta=\Forb(\widehat{\mathcal A}_\beta)$. For each $A\in\widehat{\mathcal A}_\beta$, enumerate all overlap unions of one copy of $A$ and one copy of $H$. Choose a partial bijection between their vertex sets, identify the chosen vertices, and take the union of the 2 edge sets. The resulting finite family $\mathcal O_\beta$ satisfies
\[
 \Forb(\mathcal O_\beta)
 =\Forb(\widehat{\mathcal A}_\beta)\cup\Forb(H)
 =\mathcal C^\dagger_\beta.
\]
Indeed, a host is outside the union exactly when it contains both an original forbidden core and $H$, and the actual overlap of these 2 copies is one of the enumerated unions. Take $\mathcal B^\dagger_\beta$ to consist of one canonical representative of each isomorphism type in $\mathcal O_\beta$, obtained by finite permutation search. Every member is outside both classes in the union, proving the 2 assertions about its subgraphs and admissibility.

The basis is nonempty because $K_{q+1}^{(3)}$ is not admissible. In a hypothetical rooting, at least one vertex is data. Four data vertices would give more than one data-only edge. With one, 2, or 3 data vertices, at least $q-2>q/2$ roots remain. Such a root set contains a pair of every perfect matching on $Q$, so each data vertex has an empty candidate set. This is a contradiction. Moreover $K_{q+1}^{(3)}$ contains $H$ because $|V(H)|\le q$. Thus $K_{q+1}^{(3)}\notin\mathcal C^\dagger_\beta$.

It remains to compare the Lagrangian suprema. A supremum over a union is the maximum of the 2 suprema. The unmarked closed construction gives $\Lambda(\mathcal C_\beta)\ge B_*$ on every input, whereas \eqref{eq:final-singleton-safe-guard} places the extra class strictly below $B_*$. This proves equality of the weighted suprema.
\end{proof}

\begin{proposition}
\label{prop:final-singleton-target-exclusion}
For every input $\beta$ and every finite fully declared legal target $T$ of the three-uniform construction for that input,
\[
 P_\beta\nrightarrow T.
\]
\end{proposition}
\begin{proof}
Suppose a homomorphism exists, and fix a legal rooting of $T$. It has at most $q$ roots. If the product has one factor $E(B)$, core pair coverage makes the homomorphism injective on $V(B)$. Consequently $B$ is isomorphic to a subgraph of the admissible $T$, and monotonicity gives $B\in\mathcal C_\beta$, contrary to \Cref{lem:final-singleton-effective-basis}.

If there are at least 2 factors, each core contains $H$.
\Cref{prop:elementary-profiles},
with the parameters fixed above, gives some factor core $B\in\mathcal B^\dagger_\beta$ and $h=6C$ maps on $E(B)$. No map sends two core vertices to the same target root, and the maps satisfy the ordered-edge condition for every 3 distinct map indices.
\Cref{thm:many-profile-ternary-reconstruction} then gives
$B\in\mathcal C_\beta$, the same contradiction.
\end{proof}

\begin{lemma}
\label{lem:final-singleton-upper-transfer}
For a nonempty finite family $\mathcal B$ of cores with at least 3 vertices each, set $P=\prod_{B\in\mathcal B}E(B)$. Then
\[
 \pi(P)\le6\Lambda(\Forb(\mathcal B)).
\]
If $G$ has an edge and $P\nrightarrow G$, choose an optimal weighting whose support $S$ has minimum cardinality. Then $G[S]\in\Forb(\mathcal B)$ and $\lambda(G[S])=\lambda(G)$.
\end{lemma}
\begin{proof}
$3$-graphs with no edges have Lagrangian zero. Take a finite $3$-graph $G$ with at least one edge and $P\nrightarrow G$, and choose an optimal weighting whose support $S$ has minimum cardinality. Every pair of vertices in $S$ lies in an edge of $G[S]$. Otherwise the polynomial is affine when the combined weight of that pair is redistributed, so moving all its weight to one endpoint would preserve optimality and shrink the support.

If $G[S]$ contained a copy of $B\in\mathcal B$, pair coverage would extend that embedding to a homomorphism $E(B)\to G[S]$. Send each private vertex to any completing vertex for its core pair. Composing with the factor projection would give $P\to E(B)\to G$, a contradiction. Hence $G[S]\in\Forb(\mathcal B)$ and $\lambda(G)=\lambda(G[S])\le\Lambda(\Forb(\mathcal B))$. Take the supremum and apply \eqref{eq:lagdensity}.
\end{proof}

\begin{proof}[Proof of \Cref{thm:singleton}]
The guarded forbidden family is computable by
\Cref{lem:final-singleton-effective-basis}; forming its canonical
extensions and their product is a finite operation. The root-code and guard parameters are independent of $\beta$. The bounded enumerations, overlap identifications, and finite products are primitive recursive; their size bounds are given in Section~\ref{sec:consequences}.

Let $\mathcal T_\beta$ be the fully declared legal targets for the input. Target exclusion, the homomorphism--Lagrangian identity, and the product upper bound give the following bounds
\[
 6\sup_{T\in\mathcal T_\beta}\lambda(T)
 \le\pi(P_\beta)
 \le6\Lambda(\mathcal C^\dagger_\beta)
 =6\Lambda(\mathcal C_\beta).
\]
The unmarked target $T_0$ has Lagrangian $B_*$ for every input. On a nonhalting input, the global weighted bound gives $\Lambda(\mathcal C_\beta)=B_*$, so both sides force $\pi(P_\beta)=6B_*$. On a halting input, the closed marked target $T_\beta\in\mathcal T_\beta$ has Lagrangian greater than $B_*$, so the left inequality gives $\pi(P_\beta)>6B_*$.
\Cref{lem:wi-rational-baseline} gives the fixed rational
threshold $\tau_3=6B_*\in(0,1)$.

An algorithm deciding either $\pi(F)>\tau_3$ or $\pi(F)=\tau_3$, applied to $F=P_\beta$, would therefore decide whether the fixed universal machine halts. This is impossible.
\end{proof}

\section{Effective suspension and numerical certificates}\label{sec:v3-suspension}

Write $q_3\coloneqq q$, $m_3\coloneqq3m$, and $t_3\coloneqq t_0$ for the parameters of the numerical construction. The factor $3$ in $m_3$ accounts for the 3 disjoint marked realizations.

For a $3$-graph $H$ and an $(r-3)$-set $S$ disjoint from $V(H)$, the \emph{suspension} $S*H$ has vertex set $S\sqcup V(H)$ and edges $S\cup e$ for $e\in H$. This operation scales the Lagrangian by a rational factor; compare Liu and Pikhurko~\cite[Proposition 2.3]{LiuPikhurko2023}. To obtain finite forbidden families, we also need a finite basis for the class of suspended hypergraphs.

\begin{lemma}\label{lem:v3-effective-suspension}
Let $\mathcal A$ be a finite family of $3$-graphs, each having an edge, and suppose that $\mathcal C=\Forb(\mathcal A)$ contains a $3$-graph with an edge. For every integer $r\ge3$, one can effectively construct a finite family $\mathcal A^{[r]}$ of $r$-graphs such that
\[
 \Lambda\bigl(\Forb(\mathcal A^{[r]})\bigr)
 =\frac{27}{r^r}\Lambda(\mathcal C).
\]
The construction is uniform in $r$ and $\mathcal A$. Every member of $\mathcal A^{[r]}$ has an edge.
\end{lemma}
\begin{proof}
For $r=3$, take $\mathcal A^{[3]}=\mathcal A$. Suppose that $r\ge4$, and let $s=r-3$. Let $\mathcal D$ consist of all edgeless $r$-graphs and all $r$-graphs $G$ admitting an $s$-set $S$ contained in every edge such that the $3$-graph
\[
 G_S\coloneqq\bigl\{e\setminus S\colon e\in G\bigr\},
 \qquad V(G_S)\coloneqq V(G)\setminus S,
\]
belongs to $\mathcal C$. Deleting a core vertex deletes every edge; other subgraphs have links contained in $G_S$. Thus $\mathcal D$ is monotone. Membership is decidable by testing each possible core and its link.

Let $a\coloneqq\max(\{0\}\cup\{v(J)\colon J\in\mathcal A\})$ and $N\coloneqq r+\binom r3\max\{r,a\}$. Every $G\notin\mathcal D$ has a subgraph outside $\mathcal D$ on at most $N$ vertices. To see this, fix an edge $e_0$. A common core must be one of the $\binom r3$ sets $S\in\binom{e_0}s$. For each candidate $S$, retain an edge that does not contain $S$, if one exists. Otherwise retain the vertices of a forbidden copy $J\in\mathcal A$ in $G_S$, together with its corresponding suspended edges. Each witness adds at most $\max\{r,a\}$ vertices to $e_0$. Their union rules out every candidate core, proving the bound.

Enumerating all nonmembers of $\mathcal D$ of order at most $N$, up to isomorphism, therefore gives the required basis $\mathcal A^{[r]}$. Every member has an edge.

To compute the Lagrangian scaling factor, fix $H\in\mathcal C$ and an $s$-set $S$ disjoint from $V(H)$. For total weight $t$ on $V(H)$ and core weights $y_1,\ldots,y_s$ with $t+\sum_i y_i=1$, AM--GM on $t/3,t/3,t/3,y_1,\ldots,y_s$ gives
\[
 \lambda(S*H;\x)\le t^3\lambda(H)\prod_{i=1}^s y_i
 \le\frac{27}{r^r}\lambda(H).
\]
Equality is attained by giving $V(H)$ total weight $3/r$, in proportions attaining $\lambda(H)$, and every core vertex weight $1/r$. Every $r$-graph in $\mathcal D$ with an edge has the form $S*H$ for a link $H\in\mathcal C$, so taking suprema proves the identity.
\end{proof}

For $r\ge4$, apply \Cref{lem:v3-effective-suspension} to $\mathcal A=\cQ(\{P_\beta\})$, and denote its output by $\mathcal A_{r,\beta}$. Thus $\mathcal D_{r,\beta}\coloneqq\Forb(\mathcal A_{r,\beta})$ consists of the edgeless $r$-graphs and the suspensions whose $3$-graph links are $P_\beta$-hom-free. Define
\[
 \mathcal F_{3,\beta}\coloneqq\{P_\beta\},\qquad
 \mathcal F_{r,\beta}\coloneqq\bigcup_{J\in\mathcal A_{r,\beta}}\Ext(J)
 \quad(r\ge4).
\]
By \cref{lem:pairtransfer,prop:approx},
\begin{equation}\label{eq:v3-numericalexact}
 \pi(\mathcal F_{r,\beta})
 =\frac{27r!}{r^r}\Lambda\bigl(\Forb(\cQ(\{P_\beta\}))\bigr)
 =\frac{27r!}{6r^r}\pi(P_\beta)
 \quad(r\ge3).
\end{equation}
The fixed threshold is therefore
\begin{equation}\label{eq:v3-public-threshold}
 \tau_r\coloneqq\frac{27r!}{6r^r}\tau_3
 =\frac{27r!}{r^r}B_*.
\end{equation}
Since $\tau_{r+1}/\tau_r=(r/(r+1))^r<1$ and $\tau_3\in(0,1)\cap\mathbb Q$, every $\tau_r$ lies in that interval. The numerical dichotomy follows from \Cref{thm:singleton}. We next prove the fixed-template assertion and the certificate bounds.

A halting history first gives a small template with a rational weighting above the threshold. A rational blowup replaces this weighting by uniform vertex weights, allowing improvement to be verified by an integer inequality involving only the edge count and order.

\begin{proposition}\label{prop:v3-certificates}
For every fixed $r\ge3$, there are an effectively specified integer $D_r$ and a $D_r$-vertex $r$-graph $K_{*,r}$ with the following properties.
\begin{enumerate}
\item The $r$-graph $K_{*,r}$ is $\mathcal F_{r,\beta}$-hom-free for every input. Its uniform weighting is optimal and $\lambda_r(K_{*,r})=\tau_r$.
\item A halting history of area $n=A_{\rm hist}(\beta)$ gives an $\mathcal F_{r,\beta}$-hom-free improving template of order $q_3+m_3n+r-3$ and a hom-free certificate $G$ on $D_rn$ vertices satisfying
\begin{equation}\label{eq:v3-certificate-gap}
 \frac{r!|G|}{(D_rn)^r}
 =\tau_r+\frac{27r!t_3^3}{r^rm_3^3n^3}.
\end{equation}
\item Every $\mathcal F_{r,\beta}$-hom-free $r$-graph $G$ with $\lambda_r(G)>\tau_r$ satisfies $A_{\rm hist}(\beta)\le v(G)$ and $T(\beta)\le v(G)$.
\end{enumerate}
The finite constructions and both certificate conversions are primitive recursive for fixed $r$ and uniformly computable in $r$.
\end{proposition}
\begin{proof}
We first construct the baseline and improving templates, then recover a history from an arbitrary improving template.

Choose $D_r$ to clear the denominators of $1/r$, $3(1-t_3)/(rq_3)$, and $3t_3/(rm_3)$. Use an unmarked fully declared realization with 3 vertices in each physical type. Give every root weight $(1-t_3)/q_3$ and every data vertex weight $t_3/m_3$. By \Cref{prop:ternary-marker-lower-witnesses,prop:final-singleton-target-exclusion}, this weighting attains $B_*$ and the target is $P_\beta$-hom-free for every input.

For $r\ge4$, suspend this target on $r-3$ new vertices, each of weight $1/r$, and scale every old weight by $3/r$. The suspended hypergraph belongs to $\mathcal D_{r,\beta}$, so \Cref{lem:pairtransfer} makes it $\mathcal F_{r,\beta}$-hom-free. At $r=3$ use the target itself. In both cases the normalized value is $\tau_r$, and the chosen denominators give a $D_r$-vertex rational blowup $K_{*,r}$. Applying \eqref{eq:rationalblowup} and the upper bound for the fixed nonhalting input shows that its uniform weighting is optimal.

For a halting history of area $n$, \Cref{thm:weighted-compiler} gives a fully declared target on $q_3+m_3n$ vertices, with root weights $(1-t_3)/q_3$ and data weights $t_3/(m_3n)$. Its value is $B_*+(t_3/(m_3n))^3$, and target exclusion again gives $P_\beta$-hom-freeness. Suspension and the same weight assignment give normalized value
\[
 r!\left(\frac3r\right)^3\left(\frac1r\right)^{r-3}
 \left[B_*+\left(\frac{t_3}{m_3n}\right)^3\right]
 =\tau_r+\frac{27r!t_3^3}{r^rm_3^3n^3}.
\]
The weights have common denominator $D_rn$, proving part \textup{(ii)} by rational blowup.

For part \textup{(iii)}, first suppose $r=3$ and choose an optimal weighting of $G$ whose support $S$ has minimum cardinality. By \Cref{lem:final-singleton-upper-transfer}, the induced $3$-graph $J\coloneqq G[S]$ belongs to $\mathcal C_\beta\cup\Forb(H)$ and satisfies $\lambda(J)=\lambda(G)>B_*$. The strict bound $\Lambda(\Forb(H))<B_*$ in \eqref{eq:final-singleton-safe-guard} forces $J\in\mathcal C_\beta$. The weighted bound in \Cref{prop:ternary-marker-nonhalting-bound} now gives
\[
 A_{\rm hist}(\beta)\le v(J)\le v(G).
\]

For $r\ge4$, choose an optimal weighting of $G$ whose support has minimum cardinality. By \Cref{lem:pairtransfer}, the subhypergraph induced by this support belongs to $\mathcal D_{r,\beta}$. Write it as $S*J$, where $J$ is $P_\beta$-hom-free. The suspension identity gives
\[
 6\lambda(J)=\frac{6r^r}{27r!}\lambda_r(G)
 >\frac{6r^r}{27r!}\tau_r=\tau_3.
\]
Apply the $r=3$ case to $J$. Since $v(J)\le v(G)$, it gives the required area and time bounds.

All forward constructions use bounded finite operations on the history. For the reverse conversion from a certificate with uniform vertex weights, check that $G$ is $\mathcal F_{r,\beta}$-hom-free and that $r!|G|>\tau_r v(G)^r$, then simulate $\mathsf U(\beta)$ for $v(G)$ steps. Part \textup{(iii)} guarantees a halt whenever the checks pass. The finite bases, products, homomorphic images, suspensions, and rational blowups have primitive recursive size bounds for fixed $r$; explicit output bounds are recorded in Section~\ref{sec:consequences}.
\end{proof}

Write $\tau_r=p_r/b_r$ in lowest terms. A certificate with uniform weights satisfies
\begin{equation}\label{eq:balancedimprovement}
 b_rr!|G|>p_rv(G)^r.
\end{equation}
The numerical dichotomy and \Cref{prop:v3-certificates} prove \Cref{thm:numerical}.

\section{From a density gap to two phases}\label{sec:limits}
The next three sections turn the inner density gap into a structural dichotomy, with a fixed neutral limit at zero gap and two copies of the whole inner extremal space at positive gap. The weighted calculation comes first; the passage to ordinary free hypergraphs is proved separately in Section~\ref{sec:transfer}.

Let $g\coloneqq4$ and $\tau\coloneqq\tau_3=p/q\in(0,1)$ in lowest terms. Let $\cB_\beta$ be the family of homomorphic images of $P_\beta$ after its isolated vertices have been deleted. Set
\[
 z\coloneqq\pi(P_\beta),\qquad \delta\coloneqq z-\tau\in[0,1].
\]
At each fixed target uniformity, the construction uses only $\cB_\beta$ and the fixed rational $\tau$; the value $z$ is used only to analyze its extremal configurations.

The class $\Forb(\cB_\beta)$ is blowup closed and has Lagrangian supremum $z/6$. Moreover, the common vanishing conditions proved in \Cref{prop:approx} give $\cE(\cB_\beta)=\cE(P_\beta)$. The arguments in this and the next 2 sections use only a finite inner forbidden family, blowup closure of its free class, and $\tau\le z\le1$. They therefore apply to the homomorphic closure of any finite family of $3$-graphs after deletion of isolated vertices, whenever its density is at least $\tau$. The choice $\cB_\beta$ specializes the construction to the numerical problem.

Choose the background payoff as $L\coloneqq8g(g+1)q=160q$, independently of the input.
\paragraph{Two types and witnessed edges.}
A legal object consists of a partition $C=U\dotcup W$ into active vertices $U$ and neutral vertices $W$, together with two $\cB_\beta$-free $3$-graphs $J^{(1)},J^{(2)}$ on $U$, at most one of which has an edge. A neutral vertex supplies the extra position needed to report an inner $3$-edge on a $g$-set.

Give the data vertices a probability weighting $\x$, and let $\alpha\coloneqq x(U)$, so the neutral weight is $1-\alpha$. We choose the payoffs to give the upper bound $L/g!+(q/6)f_\delta(\alpha)$, where $f_\delta(\alpha)\coloneqq\delta\alpha^3(1-\alpha)-\alpha^4$. The factor $\delta$ makes the reward depend on the excess of the inner density over $\tau$, while the negative quartic term limits the total active weight.

For each $g$-set $X$, let $j(X)\coloneqq|X\cap U|$. Let $\sigma(X)=a$ if $j(X)=3$ and $X\cap U\in J^{(a)}$, for $a\in\{1,2\}$, and let $\sigma(X)=0$ otherwise. The record on $X$ is $\ell(X)\coloneqq(j(X),\sigma(X))$. We call $1,2$ the phase labels. Thus an all-active set reports no inner-edge information. There are $g+3$ possible labels, independent of the input. The label $(0,0)$ is neutral and is denoted by $0$. For each label $\ell=(j,\sigma)$, define the payoff $n_\ell$ and correction $s_\ell$ by
\begin{equation}\label{eq:payoff}
 n_\ell\coloneqq L+s_\ell,\qquad
 s_\ell\coloneqq-p\mathbf1_{j=3}-gq\mathbf1_{j=g}
                       +q\mathbf1_{\sigma\ne0}.
\end{equation}
In particular $|s_\ell|\le gq$ and $L-gq\le n_\ell\le L+q$. In Section~\ref{sec:weighted}, $n_\ell$ will be realized as the number of code words assigned to label $\ell$. The background contribution lost when records are missing will dominate any possible increase in the correction terms.

For a $4$-set $X=\{u_1,u_2,u_3,w\}$ with $u_1,u_2,u_3\in U$ and $w\in W$, the label is $(3,a)$ if $\{u_1,u_2,u_3\}\in J^{(a)}$, and $(3,0)$ otherwise. The corresponding payoffs are $L-p+q$ and $L-p$, with $p/q=\tau$ setting the density threshold. An all-active $4$-set has payoff $L-4q$; the deduction of $4q$ supplies the quartic penalty.

Write $J\coloneqq J^{(1)}\cup J^{(2)}$ for the \emph{inner projection} and $\Psi(J;\x)$ for the score of the complete legal object. Then
\[
 \Psi(J;\x)=\sum_{j=0}^g n_{(j,0)} e_j(\x_U)e_{g-j}(\x_W)
                      +q(1-\alpha)\lambda(J;\x_U).
\]
For $\alpha>0$ let $z_J\coloneqq 6\lambda(J;\x_U)/\alpha^3\le z$; at $\alpha=0$ let $z_J\coloneqq0$. Applying $e_j(\y)\le(\sum y_i)^j/j!$ to the nonnegative terms and collecting the background and penalty contributions gives
\[
 \Psi(J;\x)\le\frac L{g!}-\frac p6\alpha^3(1-\alpha)
       -\frac{gq}{g!}\alpha^g+q(1-\alpha)\lambda(J;\x_U).
\]
Substituting $g=4$, $p=q\tau$, and $\lambda(J;\x_U)=z_J\alpha^3/6$ gives the relaxed upper bound
\begin{equation}\label{eq:colourupper}
 \begin{aligned}
 \Psi(J;\x)&\le\overline S(J;\x)\coloneqq\frac L{g!}+\frac q{6}
       [(z_J-\tau)\alpha^3(1-\alpha)-\alpha^g],\\
 0\le\overline S(J;\x)-\Psi(J;\x)
       &\le\frac{L+q}{g!}\binom g2\sum_{v\in C}x_v^2.
 \end{aligned}
\end{equation}
The inner-edge reward is unchanged by the relaxation. Since $0\le n_{(j,0)}\le L+q$, the loss satisfies
\[
 \begin{aligned}
 \overline S(J;\x)-\Psi(J;\x)
 &=\sum_{j=0}^g n_{(j,0)}
   \left[\frac{\alpha^j(1-\alpha)^{g-j}}{j!(g-j)!}
         -e_j(\x_U)e_{g-j}(\x_W)\right]\\
 &\le(L+q)\left[\frac1{g!}-e_g(\x_C)\right].
 \end{aligned}
\]
The last bracket is $1/g!$ times the probability that an ordered $g$-sample contains a repeated vertex. A union bound over pairs of positions gives the second inequality in \eqref{eq:colourupper}.

Since $z_J\le z$, the relaxed score is at most $L/g!+(q/6)f_\delta(\alpha)$, as intended. To maximize $f_\delta$ on $[0,1]$, write $f_\delta(\alpha)=\alpha^3[\delta-(1+\delta)\alpha]$. At zero gap its maximum is attained only at $\alpha=0$. For $\delta>0$, it is nonpositive outside $[0,\delta/(1+\delta)]$. On this interval, let $y\coloneqq(1+\delta)\alpha/\delta$. Then
\[
 f_\delta(\alpha)=\frac{\delta^4}{(1+\delta)^3}y^3(1-y)
 \le\frac{27\delta^4}{4^4(1+\delta)^3},
\]
with equality exactly when $y=3/4$. Since $g=4$, adding $L/g!$ to $q/6$ times this maximum gives an upper bound for the scores of complete legal objects. This bound and the corresponding optimal active proportion are
\begin{equation}\label{eq:phasefunction}
 S(z)\coloneqq\frac L{g!}+\frac{9q}{2g^g}
                       \frac{\delta^g}{(1+\delta)^3},\qquad
 \alpha_*\coloneqq\frac{3\delta}{g(1+\delta)}.
\end{equation}
In particular, $1-\alpha_*>1/2$.

We now show that complete legal objects can approach this bound. For $\delta>0$, choose a sequence of $\cB_\beta$-free $3$-graphs $J$ with probability weightings $\y$ such that $6\lambda(J;\y)\to z$. Give their active vertices weights $\alpha_*\y$, add neutral vertices of total weight $1-\alpha_*$, and split the active and neutral weights among sufficiently many clones that $\sum_vx_v^2\to0$. By \eqref{eq:colourupper}, the actual scores then approach $S(z)$. At $\delta=0$, take $n$ neutral vertices of weight $1/n$ and give every $g$-set the neutral label. The score is $L\binom ng/n^g$, which tends to $L/g!=S(\tau)$. Thus $S(z)$ is the supremum of the scores of complete legal objects.

To analyze complete legal objects whose scores approach $S(z)$, decompose the deficit as
\begin{equation}\label{eq:phasedeficit}
 \begin{aligned}
 S(z)-\Psi(J;\x)
 ={}&\bigl[\overline S(J;\x)-\Psi(J;\x)\bigr]\\
 &+\frac q{6}\bigl[f_\delta(\alpha_*)-f_\delta(\alpha)\bigr]\\
 &+\frac q{6}(z-z_J)\alpha^3(1-\alpha).
 \end{aligned}
\end{equation}
All 3 terms are nonnegative. The unique scalar maximizer therefore forces $\alpha\to\alpha_*$ along any sequence of complete legal objects whose scores approach $S(z)$. At zero gap, the total active weight tends to zero, so the phase labels become indistinguishable in the limit. At positive gap, $0<\alpha_*<1$, and the last term forces $z_J\to z$. Both phases have the same optimal active weight $\alpha_*$; the two copies of the inner extremal space arise from the choice of which of $J^{(1)},J^{(2)}$ is nonempty.

\paragraph{Partial records.}
We next extend the score bound $S(z)$ to systems with missing records that satisfy bounded local legality tests. Let
\[
 Q\coloneqq g+1+\max\{2g,\max_{F\in\cB_\beta}(v(F)+|F|)\}.
\]
An empty maximum is zero. The bound $v(F)+|F|$ allows a witness to retain each vertex of a forbidden inner hypergraph and one neutral vertex witnessing each inner edge.

A partial system assigns a label to each member of a chosen family of $g$-sets. It is locally legal if, on each set of at most $Q$ vertices, there is a complete legal object on that same set whose labels agree with all available records. Conflicting labels are rejected. The test is finite and effective, and deletion preserves it. Data cloning also preserves it. Indeed, realize the records on the images of a tested clone set and clone that realization. The inner class is blowup closed, and the vertices supporting each cloned record have distinct images in the original data set.

For a partial system on $C$ and a $(g-1)$-set $B\subseteq C$, define the neighborhood $N(B)$ and the reported active counts by
\[
 N(B)\coloneqq\{x\in C\setminus B\colon B\cup\{x\}\text{ has a record}\},
 \qquad j_x\coloneqq j(B\cup\{x\})\quad(x\in N(B)).
\]
\begin{samepage}
The next lemma reconstructs a complete legal object and bounds the disagreement probability and the score of the partial system. An ordered sample is erased if it repeats a vertex or consists of distinct vertices whose set has no record.

\begin{lemma}\label{lem:recordrepair}\label{prop:recordpenalty}
Let $\mathcal R$ be locally legal and let $\x$ be a probability weighting on its data vertices. Let
\begin{equation}\label{eq:erasure}
 \begin{aligned}
 \epsilon&\coloneqq1-g!\sum_{X\text{ recorded}}\prod_{v\in X}x_v\\
 &=\bigl[1-g!e_g(\x_C)\bigr]
   +g!\sum_{\substack{X\in\binom Cg\\X\text{ not recorded}}}\prod_{v\in X}x_v.
 \end{aligned}
\end{equation}
The two terms in the second line are the probabilities that a sample has repeated vertices and that it consists of distinct vertices whose underlying set has no record, respectively. The following conclusions hold.
\begin{enumerate}
\item There is a set $V'\subseteq C$ with $x(V')\ge1-\epsilon$ whose available records have a legal realization.
\item Keep this realization on $V'$, declare every vertex of $C\setminus V'$ neutral, and use the resulting legal object to label every $g$-set of $C$, without changing $\x$. The disagreement probability on an ordered $g$-sample is at most $(g+1)\epsilon$; repetitions and missing records count as disagreements.
\item The score of the partial system satisfies
\begin{equation}\label{eq:recordpenalty}
 \Psi(\mathcal R;\x)\coloneqq\sum_{X\text{ recorded}}n_{\ell(X)}\prod_{v\in X}x_v
 \le S(z)-\frac{L}{2g!}\epsilon.
\end{equation}
\end{enumerate}
\end{lemma}
\end{samepage}
\begin{proof}
We first find a set of weight at least $1-\epsilon$ on which all available records have a legal realization. Making the remaining vertices neutral then gives the disagreement bound, and comparing the lost background contribution with the change in corrections gives the score bound.

The local step is to realize all records supported in a common neighborhood $N(B)$. When the counts $j_x$ take two values, the records through $B$ determine the types in $N(B)$. When all counts are equal, we instead construct a realization using a single type.

\begin{claim}\label{lem:commonlink}
For every $(g-1)$-set $B\subseteq C$, the restriction of the partial record system to $N(B)$ admits a complete legal realization on $N(B)$.
\end{claim}
\begin{proof}
In a local realization on $B\cup\{x,y\}$, with active set $U$, we have $j_x-j_y=\mathbf1_{x\in U}-\mathbf1_{y\in U}$. Hence $|j_x-j_y|\le1$ for all $x,y\in N(B)$.

Suppose first that 2 values $s,s+1$ occur, and choose $u,v\in N(B)$ with these respective counts. Every local realization containing $B,u,v$ and $x$ forces
\begin{equation}\label{eq:typerecovery}
 |B\cap U|=s,\qquad \mathbf1_{x\in U}=j_x-s.
\end{equation}
Use this equation to assign a type to every vertex in $N(B)$. For each recorded $g$-set $X\subseteq N(B)$, the set $B\cup\{u,v\}\cup X$ has at most $2g+1\le Q$ vertices. By local legality, it admits a complete legal realization, which verifies the active count reported by $X$.

For each record with phase label $\sigma\in\{1,2\}$, add its set of $3$ active vertices as an edge of $J^{(\sigma)}$. It remains to check that these witnessed edges give a legal realization of the available records. A failure has one of the following forms.
\begin{enumerate}[beginpenalty=10000]
\item Two records have different nonzero phase labels. They use at most $2g$ vertices.
\item A record with a nonzero phase label and a record labeled $(3,0)$ have the same $3$ active vertices. They use at most $g+1$ vertices.
\item The witnessed inner edges contain a forbidden $F\in\cB_\beta$. Retain the vertices of $F$ and, for each of its edges, one neutral vertex supplying a record with a nonzero phase label. This gives a witness on at most $v(F)+|F|$ vertices.
\end{enumerate}
Every possible failure therefore has a witness on at most $Q-g-1$ vertices. Adjoining $B,u,v$ gives at most $Q$ vertices. Their local realization must have the types fixed by \eqref{eq:typerecovery} and exhibit the same failure, a contradiction. Thus the witnessed inner hypergraphs realize all available records in the neighborhood.

Suppose instead that all $j_x$ are equal. For a recorded $g$-set $X\subseteq N(B)$, a realization on $B\cup X$ gives every vertex of $X$ the same type, because $j_x-|B\cap U|$ is constant. Its label is therefore $(0,0)$ or $(g,0)$. If 2 recorded sets $X,Y$ had different such labels, a complete legal realization on $B\cup X\cup Y$, which has at most $3g-1\le Q$ vertices, would assign the same type to every vertex of $X\cup Y$, a contradiction. Assign all of $N(B)$ the indicated type and use empty inner hypergraphs. If no set in the neighborhood is recorded, use the neutral type.
\end{proof}

To find a neighborhood of weight at least $1-\epsilon$, sample $g-1$ vertices independently according to $\x$. If the chosen vertices are distinct, let $B$ be their set; otherwise let $N(B)=\varnothing$ for this sample. A further independent vertex belongs to $N(B)$ exactly when the resulting ordered $g$-sample consists of distinct vertices whose underlying set has a record. Consequently
\[
 \mathbb E_B[x(N(B))]=1-\epsilon.
\]
Some neighborhood has total weight at least $1-\epsilon$; if $\epsilon=1$, take the empty set. Apply \Cref{lem:commonlink} to realize its available records and make all remaining data vertices neutral. If a sample consists of distinct vertices in the chosen neighborhood and its underlying set has a record, the reconstructed label agrees with the original label. Thus disagreement requires an erased sample or a vertex outside the neighborhood, and has probability at most $\epsilon+g\epsilon$.

We compare the partial score with the relaxed score $\overline S(J;\x)$ of the reconstructed object. For this comparison, both corrections must be defined on all ordered $g$-samples, including repetitions. Let $s_{\rm obs}\coloneqq s_{\ell(X)}$ when the sample has distinct entries and its underlying $g$-set $X$ has a record, and let $s_{\rm obs}\coloneqq0$ otherwise. Then
\[
 g!\Psi(\mathcal R;\x)=L(1-\epsilon)+\mathbb E[s_{\rm obs}].
\]
For the reconstructed legal object, extend its correction to samples with repeated vertices as follows. Let $a$ be the number of positions occupied by active vertices, and let $I_J$ indicate that $a=3$ and that the vertices occupying those positions are distinct and form an edge of $J$. Set
\[
 s_J\coloneqq-p\mathbf1_{a=3}-gq\mathbf1_{a=g}+qI_J.
\]
This agrees with the reconstructed correction on samples with distinct vertices. On a sample with repeated vertices, $a$ counts active vertices with multiplicity, but $I_J=0$, since $3$ distinct active vertices and one neutral vertex would make all $g$ entries distinct. In particular $|s_J|,|s_{\rm obs}|\le gq$. Since $g=4$, independent sampling gives $\Pr[a=3]=g\alpha^3(1-\alpha)$ and $\Pr[a=g]=\alpha^g$. There are $g$ choices for the neutral position and $3!$ orders of the vertices of an inner edge, so
\[
 \mathbb E[I_J]=g\cdot3!(1-\alpha)\lambda(J;\x_U)
              =g!(1-\alpha)\lambda(J;\x_U).
\]
Substituting these expectations into the definition of $s_J$ gives
\[
 \begin{aligned}
 L+\mathbb E[s_J]
 &=L-gp\alpha^3(1-\alpha)-gq\alpha^g
       +qg!(1-\alpha)\lambda(J;\x_U)\\
 &=g!\overline S(J;\x)\le g!S(z).
 \end{aligned}
\]
Outside the disagreement event, $s_{\rm obs}=s_J$. Since both corrections have absolute value at most $gq$, we have $\mathbb E[s_{\rm obs}]-\mathbb E[s_J]\le2gq(g+1)\epsilon$. The background contribution is smaller by exactly $L\epsilon$. Since $L=8g(g+1)q$, it follows that
\[
 g!\Psi(\mathcal R;\x)
 \le g!S(z)-L\epsilon+2gq(g+1)\epsilon
 =g!S(z)-\frac34L\epsilon.
\]
This implies \eqref{eq:recordpenalty}.
\end{proof}
For any sequence of locally legal record systems whose scores tend to $S(z)$, we have $\epsilon\to0$. If the first two sampled vertices coincide, the sample is erased. Hence $\sum_vx_v^2\le\epsilon\to0$.

For uniform weights on $|C|=n\ge g$ data vertices, let $h$ be the proportion of $g$-sets with no record. Each recorded $g$-set contributes once to $|N(B)|$ for each of its $g$ subsets $B$ of size $g-1$. Hence
\[
 \frac{1}{\binom n{g-1}}\sum_{B\in\binom C{g-1}}|N(B)|
 =\frac{g(1-h)\binom ng}{\binom n{g-1}}=(n-g+1)(1-h).
\]
Some $B$ attains at least this average, and \Cref{lem:commonlink} gives a legal realization on $N(B)$ preserving every available record.

\section{Code supports and weighted extremal limits}\label{sec:weighted}
The encoding makes the edge-polynomial value a fixed multiple of the record score at prescribed root weights. To classify all extremal limits, we must show that, along every sequence of nearly optimal weighted models, the root weights converge to the prescribed values. We must also recover the inner homomorphism densities from the unlabeled hypergraph. The code inequalities control the weights, while the neutral frame identifies the root roles and supplies tests that recover the inner densities.

\subsection{Encoding the record scores}
Since $(g+3)(L+q)=1127q$, choose the code parameters, independently of the input, as
\[
 b\coloneqq2+2\left\lceil\log_2(1127q)\right\rceil,
 \qquad r_{\rm str}\coloneqq3b+4.
\]
Thus $L$, $b$, and $r_{\rm str}$ do not depend on the orders of the inner obstructions.

The encoding will use $b$ blocks of 5 roots. Each ordinary edge will consist of 4 data vertices, 3 roots from each block, and all suspension roots. A data record will determine which combinations of root triples are allowed. Within each block, the two possible triple families must have the same maximizing weights. On $\mathbb Z_5$ let
\[
 P_0\coloneqq\{\{i,i+1,i+2\}\colon i\in\mathbb Z_5\},\qquad
 P_1\coloneqq\binom{\mathbb Z_5}{3}\setminus P_0.
\]
Multiplication by $2$ interchanges these two $3$-graphs.
\begin{lemma}\label{lem:pentagon}
For every nonnegative vertex weighting $\x$ of total weight $s$, we have $\lambda(P_\eta;\x)\le s^3/25$ for $\eta\in\{0,1\}$. If $s>0$, equality holds exactly when all five vertex weights equal $s/5$.
\end{lemma}
\begin{proof}
With indices modulo five, expansion gives
\begin{equation}\label{eq:pentagonsquares}
 \begin{aligned}
 400\left(\frac{s^3}{25}-\lambda(P_0;\x)\right)
 =\sum_i x_i\Bigl[&(5x_i-s)^2
 +5(x_{i+1}+x_{i-1}-x_{i+2}-x_{i-2})^2\\
 &+50(x_{i+1}-x_{i-1}-x_{i+2}+x_{i-2})^2\Bigr].
 \end{aligned}
\end{equation}
Every term is nonnegative. At equality each positive $x_i$ equals $s/5$, so, when $s>0$, all five do. The permutation $i\mapsto2i$ of $\mathbb Z_5$ maps $P_0$ to $P_1$ and transfers the same bound and equality condition to $P_1$.
\end{proof}
On disjoint five-sets $B_1,\ldots,B_b$, a word $\omega\in\{0,1\}^b$ defines $\cS_\omega\coloneqq\{E_1\cup\cdots\cup E_b\colon E_j\in P_{\omega_j}(B_j)\}$. These support families are disjoint. If block $j$ has total weight $s_j$, then
\begin{equation}\label{eq:binarybound}
 \lambda(\cS_\omega;\y)=\prod_j\lambda(P_{\omega_j};\y|_{B_j})
 \le25^{-b}\prod_j s_j^3,
\end{equation}
with simultaneous equality at blockwise uniform weights. Replacing the weights in each block by their average does not decrease the value of any nonnegative linear combination of these polynomials.

For $r\ge r_{\rm str}$ take $b$ 5-vertex blocks and $t\coloneqq r-r_{\rm str}$ suspension roots $S_*$. Let $M$ be the set of the $5b+t$ root labels. The neutral word pool will contain a path that identifies the block order up to reversal. The remaining words will be fixed by reversing the block order and complementing every bit, so this symmetry preserves each label's pool. Identify each word in $\{0,1\}^b$ with the set of its 1-coordinates, and let
\[
 R(\omega)_i\coloneqq1-\omega_{b+1-i},\qquad
 \mathcal V\coloneqq\{\omega\colon R(\omega)=\omega,\ \omega_1=0\},\qquad
 \mathcal P\coloneqq\{\varnothing,[1],\ldots,[b]\}.
\]
Here $[i]\coloneqq\{1,\ldots,i\}$. Each of the $2^{b/2-1}$ words in $\mathcal V$ has exactly $b/2$ entries equal to $1$, begins with $0$, and ends with $1$. Thus no unit-distance cube edge joins two such words or joins one to $\mathcal P$, since each word in $\mathcal V$ differs from every intermediate prefix in its first and last bits and has distance $b/2\ge2$ from each endpoint.

Order the words of $\mathcal V$ lexicographically and define the neutral word pool by
\[
 \Omega_0\coloneqq\mathcal P\cup\{\text{first }L-b-1\text{ words of }\mathcal V\}.
\]
For the remaining labels, in a fixed order, assign the next unused words of $\mathcal V$ in pools $\Omega_\ell$ of size $n_\ell$. These choices are possible because $\sum_\ell n_\ell\le(g+3)(L+q)\le2^{b/2-1}$ and $b+1\le L$.

For the bound $b+1\le L$, we have $L\ge32$ and $(g+3)(L+q)<L^2<2^{L/2-2}$; the last inequality holds at $L=32$ and is preserved when $L$ is increased by $2$. Taking logarithms and using that $L$ is even gives $b+1=3+2\lceil\log_2((g+3)(L+q))\rceil\le L-1$.

Every word in a nonneutral pool is fixed by $R$ and has first bit zero. Write $\mathcal S(\Omega)\coloneqq\bigcup_{\omega\in\Omega}\cS_\omega$. By \eqref{eq:binarybound}, these support polynomials simultaneously satisfy
\begin{equation}\label{eq:allcommonmax}
 \lambda(\mathcal S(\Omega);\y)
       \le |\Omega|25^{-b}\prod_i s_i^3,
\end{equation}
where $s_i$ is the total weight of block $i$, with equality when the vertex weights are equal within each block.

Under a partial rooting, an allowed edge is $S_*\cup Z\cup X$, where $|X|=g$, $X$ consists of data vertices, and $Z\in\mathcal S(\Omega_\ell)$; all displayed roots are assigned. It declares the record $\ell$ on $X$. Figure~\ref{fig:recordedge} illustrates the construction for a witnessed inner edge. The records declared by these edges must be locally legal. Let $\cD_\beta$ be the class admitting such a rooting.

To complete the root supports, add a new vertex for each unassigned root label and all edges specified by the code pool of each existing record, without creating new data records. For weighted models, assign weight zero to new roots.

\begin{figure}[htbp]
\centering
\begin{tikzpicture}[x=1cm,y=1cm,>=Stealth,every node/.style={font=\small},
 root/.style={rectangle,draw,fill=white,inner sep=2pt},
 chosen/.style={root,fill=black},
 active/.style={circle,draw,fill=black,inner sep=1.7pt},
 neutral/.style={circle,draw,fill=white,inner sep=1.7pt}]
 \node at (1.7,1.35) {Record on $X$};
 \filldraw[fill=black!8,draw=black!50] (0.4,0.55)--(1.8,0.55)--(1.1,-0.55)--cycle;
 \node[active,label=above:$u_1$] at (0.4,0.55) {};
 \node[active,label=above:$u_2$] at (1.8,0.55) {};
 \node[active,label=below:$u_3$] at (1.1,-0.55) {};
 \node[neutral,label=above:$w$] at (2.8,0) {};
 \node at (1.6,-1.05) {$\{u_1,u_2,u_3\}\in J^{(\sigma)}$};
 \draw[->] (3.4,0)--(5.1,0) node[midway,above] {$(3,\sigma)$};
 \node[draw,rounded corners=2pt,inner sep=8pt] (pool) at (7.15,0) {$\omega\in\Omega_{(3,\sigma)}$};
 \node at (7.15,0.85) {Choose a code word};
 \draw[->] (pool.south)--(7.15,-1.65);
 \node[anchor=west] at (7.35,-1.25) {choose one triple per block};
 \begin{scope}[yshift=-3cm]
  \node at (0.8,0.95) {$S_*$};
  \node[chosen] at (0.3,0) {};
  \node at (0.8,0) {$\cdots$};
  \node[chosen] at (1.3,0) {};
  \node at (0.8,-0.95) {all $|S_*|$ roots};
  \begin{scope}[xshift=3.4cm]
   \node at (0,0.95) {$B_1$};
   \foreach \j/\angle in {0/90,1/18,2/-54,3/-126,4/162}
    \coordinate (a\j) at (\angle:0.58);
   \filldraw[fill=black!8,draw=black!50] (a0)--(a1)--(a2)--cycle;
   \foreach \j in {0,1,2} \node[chosen] at (a\j) {};
   \foreach \j in {3,4} \node[root] at (a\j) {};
   \node at (0,-0.95) {$E_1\in P_{\omega_1}(B_1)$};
  \end{scope}
  \node at (5.5,0) {$\cdots$};
  \begin{scope}[xshift=7.6cm]
   \node at (0,0.95) {$B_b$};
   \foreach \j/\angle in {0/90,1/18,2/-54,3/-126,4/162}
    \coordinate (b\j) at (\angle:0.58);
   \filldraw[fill=black!8,draw=black!50] (b0)--(b1)--(b3)--cycle;
   \foreach \j in {0,1,3} \node[chosen] at (b\j) {};
   \foreach \j in {2,4} \node[root] at (b\j) {};
   \node at (0,-0.95) {$E_b\in P_{\omega_b}(B_b)$};
  \end{scope}
  \begin{scope}[xshift=10.2cm]
   \node at (1,0.95) {$X$};
   \node[active,label=above:$u_1$] at (0,0.25) {};
   \node[active,label=above:$u_2$] at (0.8,0.25) {};
   \node[active,label=below:$u_3$] at (0.4,-0.45) {};
   \node[neutral,label=above:$w$] at (1.8,0) {};
   \node at (1,-0.95) {all $4$ data vertices};
  \end{scope}
 \end{scope}
 \node at (6.1,-4.7) {$e=S_*\cup E_1\cup\cdots\cup E_b\cup X,\qquad |e|=|S_*|+3b+4=r$};
\end{tikzpicture}
\caption{Encoding a witnessed inner edge. Filled and open circles are active and neutral data vertices. In the lower row, filled squares are the selected roots; the other roots in each block are omitted from the edge. The displayed first and last triples use $\omega_1=0$ and $\omega_b=1$, as required for every nonneutral pool. Each code word permits $5^b$ choices of root support; one is shown schematically, with intermediate blocks suppressed. Suspension roots are absent when $S_*=\varnothing$.}
\label{fig:recordedge}
\end{figure}

Fix $D\in\cD_\beta$ with an admissible rooting and a probability weighting $\x$ on all its vertices. Let $C$ be the data set and let $T\coloneqq x(C)$. For $T>0$, give the extracted record system the probability weighting $\x_C/T$, and let $\epsilon$ be its erasure probability. Let
\[
 K\coloneqq25^{-b}T^g\prod_i s_i^3\prod_{j\in S_*}x_j,
 \qquad \epsilon\coloneqq1-\frac{g!}{T^g}
             \sum_{X\text{ recorded}}\prod_{v\in X}x_v.
\]
Completing the root supports can only increase the edge polynomial. Grouping the completed edges by their data set $X$, and applying \eqref{eq:allcommonmax} with $|\Omega_\ell|=n_\ell$, gives
\[
 \begin{aligned}
 \lambda(D;\x)
 &\le\left(\prod_{j\in S_*}x_j\right)
       \sum_{X\text{ recorded}}\left(\prod_{v\in X}x_v\right)
       \lambda\bigl(\mathcal S(\Omega_{\ell(X)});\x\bigr)\\
 &\le K\sum_{X\text{ recorded}}n_{\ell(X)}
                         \prod_{v\in X}\frac{x_v}{T}.
 \end{aligned}
\]
The final sum is the score of the extracted records with probability weighting $\x_C/T$. Let $\chi_r\coloneqq3^{3b}g^g/(25^b r^r)$ and $\Lambda\coloneqq\chi_rS(z)$. Applying \eqref{eq:recordpenalty} gives
\begin{equation}\label{eq:defectglobal}
 \lambda(D;\x)\le K\left[S(z)-\frac L{2g!}\epsilon\right]
 \le KS(z)\le\Lambda.
\end{equation}
For the last step, apply AM--GM to $g$ copies of $T/g$, 3 copies of $s_i/3$ for each block $i\in[b]$, and the suspension weights. These $r$ numbers sum to one, giving $K\le\chi_r$, with equality exactly at $s_i=3/r$, $T=g/r$, and suspension weights $1/r$. At $T=0$ the edge polynomial is zero.

For the lower bound, use the complete legal objects constructed after \eqref{eq:phasefunction}, with scores tending to $S(z)$. Scale their data vertex weights by $g/r$, assign weight $3/(5r)$ to each root in a 5-vertex block and $1/r$ to each suspension root, and include all edges prescribed by the records. Since each $P_\eta$ contains 5 edges and $g+3b+t=r$, the edge-polynomial value is the record score multiplied by
\[
 \left(\frac gr\right)^g
 \left[5\left(\frac{3}{5r}\right)^3\right]^b r^{-t}
 =\frac{3^{3b}g^g}{25^b r^{g+3b+t}}=\chi_r.
\]
It therefore tends to $\Lambda$.

\subsection{A neutral frame and recovery tests}
To recover record labels after forgetting the rooting, we need a finite frame whose homomorphisms identify the root blocks and data vertices, up to root permutations preserving every code pool. The neutral pool $\Omega_0$ supplies such a frame. We first determine its word symmetries, then use the frame to construct unlabeled tests recovering the inner homomorphism densities.

The unit-distance graph of $\Omega_0$ is the prefix path $\mathcal P$ and isolated vertices. The path has exactly one edge in each coordinate direction. It is the only nontrivial connected component, so every cube isometry preserving $\Omega_0$ preserves this path and either fixes or reverses its order. If it fixes the order, it fixes $\varnothing$ and every successive coordinate direction, and hence is the identity. If it reverses the order, its composition with $R$ fixes the order and is the identity. Since $R$ fixes every added word, the cube stabilizer of $\Omega_0$ is exactly $\{1,R\}$. The geometry is illustrated in \Cref{fig:prefixframe}.

\begin{figure}[htbp]
\centering
\begin{tikzpicture}[x=1.75cm,y=0.9cm,every node/.style={font=\small}]
 \foreach \i/\word in {0/0000,1/1000,2/1100,3/1110,4/1111}{
  \fill (\i,0) circle (1.6pt);
  \node[below=4pt] at (\i,0) {$\word$};
 }
 \foreach \i/\j/\bit in {0/1/1,1/2/2,2/3/3,3/4/4}{
  \draw (\i,0)--(\j,0) node[midway,above=3pt] {$\bit$};
 }
 \fill (1.25,-1.2) circle (1.6pt);
 \node[right=4pt] at (1.25,-1.2) {$0011$};
 \fill (2.75,-1.2) circle (1.6pt);
 \node[right=4pt] at (2.75,-1.2) {$0101$};
\end{tikzpicture}
\caption{An illustration of the support geometry in dimension $4$. The prefix words form a path whose edge labels indicate the changed coordinates. The two other displayed words are fixed by reverse complementation and are isolated in the unit-distance graph; the construction in dimension $b$ adds a selected set of analogous isolated words.}
\label{fig:prefixframe}
\end{figure}

To represent the word symmetries on the root labels, let $\iota(j,a)\coloneqq(b+1-j,2a\bmod5)$ for $a\in\mathbb Z_5$, fixing the suspension roots, and let
\[
 \mathscr H_r\coloneqq\langle D_5^b\times\operatorname{Sym}(S_*),\iota\rangle.
\]
The subgroup $D_5^b\times\operatorname{Sym}(S_*)$ has order $10^b t!$. Reversing the blocks permutes its factors, and multiplication by $2$ conjugates each dihedral factor to itself. Hence $\iota$ normalizes this subgroup. Since $\iota$ reverses the block order and $\iota^2$ is a blockwise reflection, adjoining $\iota$ doubles its order, giving $|\mathscr H_r|=2\cdot10^b t!$. The dihedral factors fix every word, and $\iota$ induces $R$ on words, so every element of $\mathscr H_r$ preserves every word pool separately.

Let $D_0$ be $g+1$ formal data vertices. On $M\cup D_0$, let $\mathscr R$ be the $r$-graph whose edges are all sets $S_*\cup Z\cup X$ with $Z\in\mathcal S(\Omega_0)$ and $X\in\binom{D_0}{g}$. The following lemma recovers the root and data roles from $\mathscr R$, regarded as an ordinary unlabeled $r$-graph.
\begin{lemma}\label{lem:rigidity}
Fix $D\in\cD_\beta$ with an admissible partial rooting, and let $\phi\colon\mathscr R\to D$ be a homomorphism. Then $\phi$ maps $S_*$ onto the suspension roots and $D_0$ to data vertices. It induces a permutation of the root labels in $\mathscr H_r$, so the block order is preserved or reversed. The image of each $g$-subset of $D_0$ has the neutral label. Conversely, every element of $\mathscr H_r$ acts as an automorphism of $\mathscr R$ fixing $D_0$ pointwise.
\end{lemma}
\begin{proof}
We first identify the suspension, block, and data roles. We then show that each image data set has the neutral label and determine the allowed permutations of the root roles.

The edges corresponding to the word $\varnothing$ already cover every pair of frame vertices, since $P_0$ covers pairs in each block and every $g$-subset of $D_0$ is used. Thus the images of any two formal vertices are distinct. Every suspension root of the target belongs to every image edge, so its preimage lies in the common intersection $S_*$ of the source edges. There are exactly $t$ target suspension roots and $t$ formal suspension vertices. Thus $S_*$ maps onto the target suspension roots.

Fix a target 5-vertex block. Let $I(v)=1$ if $\phi(v)$ lies in that block, and $I(v)=0$ otherwise. For every $e\in\mathscr R$, $\sum_{v\in e}I(v)=3$. Among the frame edges associated with the word $\varnothing$, vary the $P_0$ triple in one source block while fixing the other choices. Comparing consecutive triples gives $I_i+I_{i+1}+I_{i+2}=I_{i+1}+I_{i+2}+I_{i+3}$, so $I_i=I_{i+3}$ modulo $5$. Since $\gcd(3,5)=1$, $I$ is constant on each source block. For distinct $u,v\in D_0$, keep the root choices fixed and compare the edges whose data sets are $D_0\setminus\{u\}$ and $D_0\setminus\{v\}$. Subtracting their equations $\sum_{x\in e}I(x)=3$ gives $I(u)=I(v)$. Thus $I$ is constant on $D_0$.

The constant on $D_0$ must be $0$, since otherwise an edge would contain $g>3$ vertices of the target block. Each source block on which $I=1$ contributes exactly $3$ vertices to every edge, so exactly one whole source block maps to the chosen target block. Its five images are distinct and fill that block. Applying this to all target blocks identifies the block roles bijectively, and $\phi(D_0)$ consists entirely of data vertices.

Fix $X\in\binom{D_0}{g}$. Among the source edges with data set $X$, all $10$ triples occur in each source block, because $\varnothing,[b]\in\Omega_0$. All these image edges have the same data set $\phi(X)$, so local legality forces them to declare the same record on $\phi(X)$. If this record were nonneutral, every word in its pool would have first bit $0$, so only the five $P_0$ triples would be permitted in the first target block. The identified block bijection sends the $10$ source triples to $10$ distinct target triples, which cannot all belong to that class. Thus every image data set has neutral label.

It remains to determine the permutations within and between the blocks.
\begin{claim}\label{clm:tripleclasses}
For each source root block, the induced bijection onto its target block maps $P_0$ onto either $P_0$ or $P_1$.
\end{claim}
\begin{proof}
Suppose a five-point bijection in one block sends $P_0$ partly into $P_0$ and partly into $P_1$. Since it permutes all $10$ triples and both classes have size $5$, the image of $P_1$ also meets both classes. Consider the source support families $\cS_{\varnothing}$ and $\cS_{[b]}$. Within either family, the target bit in the chosen block can take both values independently of all other blocks. In any other block, these two support families together use all $10$ triples, whose images meet both target classes. Consequently, these two support families give at least two distinct words on the remaining target coordinates, each accompanied by both values in the chosen coordinate.

All resulting words belong to $\Omega_0$, by the neutral-label conclusion. They give two distinct cube edges parallel to the same coordinate, whereas the unit-distance graph of $\Omega_0$ has exactly one edge in each coordinate direction.
\end{proof}

By \Cref{clm:tripleclasses}, each block bijection either preserves the two triple classes or exchanges them. Together with the block permutation it induces a cube isometry mapping $\Omega_0$ into itself, and hence onto itself. The stabilizer calculation shows that this isometry is either the identity or $R$. The five-point permutations preserving $P_0$ form $D_5$, because the complementary pairs of its triples form a five-cycle. The permutations interchanging $P_0$ and $P_1$ are obtained by composing an element of $D_5$ with multiplication by $2$. Suspension vertices can be permuted freely, so these actions give precisely $\mathscr H_r$. Conversely, every element of $\mathscr H_r$ preserves each word pool. Extending it by the identity on $D_0$ therefore gives an automorphism of $\mathscr R$.
\end{proof}

The role identification in \Cref{lem:rigidity} allows us to recover inner homomorphism densities from unlabeled outer tests. For $\sigma\in\{1,2\}$ and a $3$-graph $F$ with an edge and no isolated vertices, adjoin a disjoint copy of $V(F)$ to the vertex set of $\mathscr R$. Choose $y_1\in D_0$ and $Z\in\mathcal S(\Omega_{(3,\sigma)})$, and add $S_*\cup Z\cup e\cup\{y_1\}$ for every $e\in F$. Denote the ordinary unlabeled test by $T^{(\sigma)}(F)$.

We next compute the density of this test in a model compiled from a complete legal object. Let $D$ be such a model, with a probability weighting $\x$ on $V(D)$. By \Cref{lem:rigidity},
\begin{equation}\label{eq:readback}
 t(T^{(\sigma)}(F),D;\x)=|\mathscr H_r|\!\prod_{i\in M}x_i\,
       (g+1)!e_{g+1}(\x_W)\,x(U)^{v(F)}
       t\left(F,J^{(\sigma)};\frac{\x_U}{x(U)}\right).
\end{equation}
By rigidity, the root assignments contribute $|\mathscr H_r|\prod_{i\in M}x_i$. The $g+1$ frame data vertices must have distinct neutral images. All such assignments are allowed, with total weight $(g+1)!e_{g+1}(\x_W)$. Since the root images already meet the prescribed root counts, the $3$ vertices coming from $F$ in each added edge must map to data vertices. The record $(3,\sigma)$ forces their images to be active and to form an edge of $J^{(\sigma)}$. Since $F$ has no isolated vertices, its vertex assignments contribute $x(U)^{v(F)}t(F,J^{(\sigma)};\x_U/x(U))$, including noninjective inner homomorphisms. Multiplying the three contributions proves \eqref{eq:readback}.

If $x(U)=0$, the outer test density is $0$, and we interpret the right-hand side of \eqref{eq:readback} as $0$.

For a single inner edge write $T^{(\sigma)}$, and let $\Sep\coloneqq t(T^{(1)},\cdot)-t(T^{(2)},\cdot)$.

\subsection{Classification of weighted extremal limits}
For a class $\cD$ of $r$-graphs, take the closure of its weighted homomorphism-density vectors. Write $\mathfrak L(\cD)$ for the vectors in this closure with edge coordinate $r!\Lambda(\cD)$. By \Cref{lem:freerealization,prop:approx}, the weighted extremal limits of $\Forb(\cB_\beta)$ form precisely $\cE(P_\beta)$.

\begin{theorem}\label{thm:loss}
We have $\Lambda(\cD_\beta)=\chi_rS(z)$. The following statements hold.
\begin{enumerate}[beginpenalty=10000]
\item Let $D_i\in\cD_\beta$, each with a fixed admissible partial rooting and a probability weighting $\x_i$. Suppose that $\lambda(D_i;\x_i)\to\Lambda$ and that all code supports of existing records are present. Each root in a 5-vertex block then has weight $3/(5r)+o(1)$, each suspension root has weight $1/r+o(1)$, and the total data weight is $g/r+o(1)$.
\item For every sequence as in \textup{(i)}, one can replace the partial records by complete legal labelings on the same data vertex sets, keeping the weights fixed. Under independent sampling of $g$ data vertices with the normalized data weighting, the old and new labels disagree with probability $o(1)$, counting missing records and repetitions as disagreements. The resulting compiled models have limiting active and neutral weights
\begin{equation}\label{eq:optimalmasses}
 a_\delta\coloneqq\frac{3\delta}{r(1+\delta)},\qquad
 w_\delta\coloneqq\frac{g+\delta}{r(1+\delta)}.
\end{equation}
\item At zero gap, $\mathfrak L(\cD_\beta)$ consists of the single density vector determined by a fixed rational finite profile.
\item At positive gap the weighted extremal space $\mathfrak L(\cD_\beta)$ is the disjoint union of two nonempty compact phases, each homeomorphic to $\cE(P_\beta)$, with separator values
\begin{equation}\label{eq:sepconstant}
 \pm c_{r,\beta},\qquad
 c_{r,\beta}\coloneqq |\mathscr H_r|\left(\frac{3}{5r}\right)^{5b}r^{-t}
                     w_\delta^{g+1}a_\delta^3 z>0.
\end{equation}
\end{enumerate}
\end{theorem}
\begin{proof}
The upper bound in \eqref{eq:defectglobal} and the construction following it give $\Lambda(\cD_\beta)=\chi_rS(z)$. Fix a sequence $D_i\in\cD_\beta$, with admissible partial rootings and probability weightings $\x_i$, such that $\lambda(D_i;\x_i)\to\Lambda$.

\begin{claim}\label{clm:outerweights}
Completing the root supports changes each fixed homomorphism density by $o(1)$. For the completed sequence, $K_i\to\chi_r$, $\epsilon_i\to0$, and all root and total data weights have the limits stated in the theorem.
\end{claim}
\begin{proof}
Completing the root supports gives $D_i^{\rm cmp}\in\cD_\beta$; extend $\x_i$ by zero on any newly supplied roots. Since the completed edge polynomial is at most $\Lambda$, the added edges have total weight at most $\Lambda-\lambda(D_i;\x_i)$. For every fixed test $F$, coupling the vertex samples and taking a union bound over its edges gives
\[
 0\le t(F,D_i^{\rm cmp};\x_i)-t(F,D_i;\x_i)
       \le r!|F|\bigl[\Lambda-\lambda(D_i;\x_i)\bigr]\longrightarrow0.
\]
Replace the sequence by its completions, retaining the notation $D_i$. Its data weight $T_i$ is positive for all sufficiently large $i$, because $\lambda(D_i;\x_i)$ approaches $\Lambda>0$, whereas this value is zero when $T_i=0$. The separate losses in \eqref{eq:defectglobal} give
\[
 \Lambda-\lambda(D_i;\x_i)
 \ge S(z)(\chi_r-K_i)
       +\frac L{2g!}K_i\epsilon_i.
\]
Both terms on the right are nonnegative. Since $S(z)>0$, first $K_i\to\chi_r>0$, and then $\epsilon_i\to0$. The vector of block total weights, suspension weights, and total data weight lies in a compact simplex. The unique AM--GM equality case in \eqref{eq:defectglobal} therefore forces every block total weight to approach $3/r$, every suspension weight to approach $1/r$, and $T_i\to g/r$.

These limits determine the total weight of each block. We must still show that its five vertex weights become asymptotically equal. For a 5-vertex block $B$ of positive total weight $s$, normalize $p_a\coloneqq x_a/s$ for $a\in B$ and let $E_B\coloneqq\sum_a p_a(5p_a-1)^2$. This quantity is unchanged by a permutation of the five coordinates. Discarding the last two nonnegative squared terms in \eqref{eq:pentagonsquares} therefore gives, for either bit $\eta\in\{0,1\}$,
\[
 \frac{s^3}{25}-\lambda(P_\eta;\x|_B)\ge\frac{s^3}{400}E_B,
 \qquad
 \lambda(P_\eta;\x|_B)\le\frac{s^3}{25}\left(1-\frac{E_B}{16}\right).
\]
Let $\bar\x$ assign weight $s/5$ to each vertex of $B$, leaving every other weight unchanged. For $\eta=0,1$, we have $\lambda(P_\eta;\bar\x|_B)=s^3/25$. Applying the last bound to every summand of the completed edge polynomial, whose remaining factors are nonnegative, gives
\[
 \lambda(D;\x)\le\left(1-\frac{E_B}{16}\right)\lambda(D;\bar\x)
                   \le\left(1-\frac{E_B}{16}\right)\Lambda.
\]
Consequently $E_B\to0$ along the near-optimal sequence. Since $\sum_a p_a(5p_a-1)=5\sum_a(p_a-1/5)^2$, Cauchy--Schwarz with weights $p_a$ gives $25[\sum_a(p_a-1/5)^2]^2\le E_B$. Thus each normalized weight approaches $1/5$. Combining this with $s\to3/r$ proves the stated limit $3/(5r)$ for every root in a 5-vertex block.
\end{proof}

\begin{claim}\label{clm:decodedphase}
The complete legal labelings in the theorem can be chosen so that their actual record scores approach $S(z)$. For every fixed test $F$, the homomorphism densities in the original and recompiled models then differ by $o(1)$. The sum of their squared data vertex weights tends to zero. In the recompiled models, the ratio of the total active weight to the total data weight tends to $\alpha_*$. At positive gap their normalized inner densities approach $z$, while at zero gap their total active weight tends to zero.
\end{claim}
\begin{proof}
Let $C_i$ be the data set, $\mathcal R_i$ the extracted records, and $\widehat\x_i$ the normalized data weighting. Applying \eqref{eq:allcommonmax} to the code supports of each recorded $g$-set gives
\[
 \lambda(D_i;\x_i)\le K_i\Psi(\mathcal R_i;\widehat\x_i)
       \le K_iS(z)\le\chi_rS(z).
\]
\Cref{clm:outerweights} and near-optimality imply $\Psi(\mathcal R_i;\widehat\x_i)\to S(z)$. Apply \Cref{lem:recordrepair} and denote the decoded complete legal object by $\mathcal J_i$, with inner $3$-graph $J_i\coloneqq J_i^{(1)}\cup J_i^{(2)}$. Its disagreement probability is at most $(g+1)\epsilon_i$. Sample $g$ data vertices independently according to $\widehat\x_i$ and compare the observed and reconstructed payoffs on the same sample, using $n_\ell$ for a recorded label $\ell$. Assign payoff zero to repetitions and, for the observed system, to missing records. Both payoffs lie in $[0,L+q]$ and agree outside the disagreement event. Since each expected payoff is $g!$ times the corresponding score,
\[
 \left|\Psi(\mathcal J_i;\widehat\x_i)-\Psi(\mathcal R_i;\widehat\x_i)\right|
 \le\frac{(L+q)(g+1)}{g!}\epsilon_i\longrightarrow0.
\]
Thus the actual finite record score of $\mathcal J_i$ approaches $S(z)$. Since all three terms on the right-hand side of \eqref{eq:phasedeficit} are nonnegative, each tends to zero. The unique maximizer of $f_\delta$ forces the total weight of active vertices under $\widehat\x_i$ to approach $\alpha_*$. If $\delta>0$, then $\alpha_*\in(0,1)$, and the inner-deficit term also forces $z_{J_i}\to z$. If $\delta=0$, the second term is $(q/6)[f_0(0)-f_0(\alpha)]=(q/6)\alpha^4$, so $\alpha\to0$. Since $T_i\to g/r$ and $g=4$, the limits of the actual active and neutral weights are
\[
 \begin{aligned}
 \frac gr\alpha_*&=\frac{3\delta}{r(1+\delta)}=a_\delta,\\
 \frac gr(1-\alpha_*)&=\frac{g+(g-3)\delta}{r(1+\delta)}
                   =\frac{g+\delta}{r(1+\delta)}=w_\delta,
 \end{aligned}
\]
as asserted in \eqref{eq:optimalmasses}.

To show that decoding preserves the outer density vector, let $\widetilde D_i$ be the model obtained by replacing each record's code pool according to the decoded labels, with the same vertex weights. Couple the vertex samples in the two models. The sampled image of a test edge can be an edge in one model but not the other only when exactly $g$ positions map to data vertices and the original and decoded records on those vertices disagree. Conditional on those positions, the data samples are independent with weighting $\widehat\x_i$, so their disagreement probability is at most $(g+1)\epsilon_i$. A union bound over the edges of a fixed test $F$ therefore gives
\[
 \begin{aligned}
 |t(F,\widetilde D_i;\x_i)-t(F,D_i;\x_i)|
 &\le |F|\binom rg T_i^g(1-T_i)^{r-g}(g+1)\epsilon_i\\
 &\le |F|(g+1)\epsilon_i\longrightarrow0.
 \end{aligned}
\]
The erasure event includes the event that the first two sampled data positions receive the same vertex, so
\[
 \sum_{v\in C_i}x_{i,v}^2\le T_i^2\epsilon_i\longrightarrow0.\qedhere
\]
\end{proof}

At zero gap, take the decoded models from \Cref{clm:decodedphase} and write $C_i,U_i$ for their data and active sets. Then $x_i(U_i)\to0$. Relabel every active vertex as neutral and compile the resulting complete neutral records. For a fixed test $F$, this changes its density by at most $v(F)x_i(U_i)=o(1)$, since a test map can be affected only if it uses an active vertex. Repeated data images have probability at most $\binom{v(F)}2\sum_{v\in C_i}x_{i,v}^2=o(1)$. Up to this error, each test density is a polynomial in the root weights and the total data weight, obtained by summing over compatible assignments of test vertices to root roles and data. Their limits therefore give the density vector $W_{0,r}$ of the complete neutral model.

Its finite profile has one role of proportion $3/(5r)$ for each of the $5b$ code roots, one role of proportion $1/r$ for each suspension root, and one neutral data role of proportion $4/r$. The allowed edge multisets consist of all suspension roles, the root roles in a support $Z\in\mathcal S(\Omega_0)$, and 4 copies of the neutral role.

\begin{claim}\label{clm:phasemaps}
If $\delta>0$, there are continuous maps $\Phi_\sigma\colon\cE(P_\beta)\to\mathfrak L(\cD_\beta)$, for $\sigma\in\{1,2\}$, whose compact images exhaust the weighted extremal space. Each map is a homeomorphism onto its image.
\end{claim}
\begin{proof}
Fix $\sigma$. We construct outer models from an inner limit and show that each outer test density depends only on that limit, defining a continuous map $\Phi_\sigma$. For $W_{\rm in}\in\cE(P_\beta)$, choose a sequence of finite $\cB_\beta$-free $3$-graphs with probability weightings whose homomorphism-density vectors converge to $W_{\rm in}$. Split their vertex weights among clones so that the largest vertex weight tends to zero. Use these as the active part of a complete legal object using phase label $\sigma$, with total active weight $a_\delta$. Add neutral vertices of total weight $w_\delta$, split so that their maximum weight tends to zero, and give every root its prescribed limiting weight. Compile these objects using the complete code pools.

The normalized inner densities approach $z$. The data weight is $g/r$, so the normalized active proportion is $ra_\delta/g=\alpha_*$. By \eqref{eq:colourupper}, the normalized record scores approach $S(z)$. At the prescribed outer weights, the edge-polynomial value of each compiled model is $\chi_r$ times its record score, and hence tends to $\Lambda$. It remains to express each outer test density in terms of $W_{\rm in}$.

Let $H$ be a fixed outer test. Write $C_i$ for the data set of the $i$th constructed outer model and $\x_i$ for its vertex weighting. The total probability of maps that repeat a data vertex is at most
\[
 \binom{v(H)}2\sum_{v\in C_i}x_{i,v}^2=o(1).
\]
Discard these maps. Assign each vertex of $H$ a root role in $M$, the active data type, or the neutral data type. For such an assignment $\theta$, write $A_\theta$ and $N_\theta$ for its active and neutral vertices. Retain only assignments for which every edge has the prescribed root profile and a compatible word pool. In particular, an edge with $3$ active vertices and one neutral vertex requires an inner edge when its word belongs to $\Omega_{(3,\sigma)}$, and requires an inner nonedge when its word belongs to $\Omega_{(3,0)}$. The pool with the other nonzero phase label is incompatible. Edges with other active counts impose only the corresponding count label.

Let $\mathcal P_\theta$ and $\mathcal N_\theta$ be the resulting sets of required inner edges and nonedges on $A_\theta$. Assignments with conflicting requirements contribute zero. Repeated root roles remain in the sum over $\theta$, since their weights need not vanish.

In each finite model, inclusion--exclusion expresses the contribution as a finite linear combination of inner homomorphism densities, with an $o(1)$ error from dropping the requirement that data images be distinct. The coefficients depend only on the prescribed root, active, and neutral weights. Since the inner densities converge to those of $W_{\rm in}$, every outer test density converges. These limits define $\Phi_\sigma(W_{\rm in})\in\mathfrak L(\cD_\beta)$. With $x_j^*$ denoting the prescribed weight of root role $j$, the limit coordinates are
\[
\begin{split}
 t(H,\Phi_\sigma(W_{\rm in}))
 ={}&\sum_{\theta\text{ compatible}}
       \left(\prod_{j\in M}(x_j^*)^{|\theta^{-1}(j)|}\right)
       a_\delta^{|A_\theta|}w_\delta^{|N_\theta|}\\
 &\hspace{10mm}\cdot\sum_{\mathcal B\subseteq\mathcal N_\theta}
           (-1)^{|\mathcal B|}
           t(\mathcal P_\theta\cup\mathcal B,W_{\rm in}),
\end{split}
\]
where each inner test has vertex set $A_\theta$. The formula shows that $\Phi_\sigma$ is independent of the approximating sequence and continuous.

To prove that the two images exhaust $\mathfrak L(\cD_\beta)$, approximate any weighted extremal limit by finite weighted models. \Cref{clm:outerweights,clm:decodedphase} give complete legal models with the same limiting density vector, optimal limiting outer weights, and asymptotically extremal inner projections. Since $z>0$, these projections eventually have an edge. Pass to a subsequence with a fixed phase label $\sigma$ and an inner limit $W_{\rm in}\in\cE(P_\beta)$. Apply the preceding finite-model expansion with the actual root weights and total active and neutral weights. Its coefficients converge to those in the displayed formula, and repeated data images contribute $o(1)$ because $\sum_{v\in C_i}x_{i,v}^2\to0$. Passing to the limit term by term identifies the outer limit with $\Phi_\sigma(W_{\rm in})$, so the two images exhaust $\mathfrak L(\cD_\beta)$.

To prove injectivity and continuity of the inverse, we use the tests $T^{(\sigma)}(F)$ to recover the inner coordinates. Fix a $3$-graph $F$ with an edge and no isolated vertices. In the approximating models, a union bound over pairs of neutral positions gives
\[
 0\le x(W)^{g+1}-(g+1)!e_{g+1}(\x_W)
 \le\binom{g+1}{2}x(W)^{g-1}\sum_{v\in W}x_v^2.
\]
Since $x(W)\to w_\delta$ and $\sum_{v\in W}x_v^2\to0$, \eqref{eq:readback} yields
\[
 t(T^{(\sigma)}(F),\Phi_\sigma(W_{\rm in}))
 =|\mathscr H_r|\left(\frac{3}{5r}\right)^{5b}r^{-t}
       w_\delta^{g+1}a_\delta^{v(F)}t(F,W_{\rm in}).
\]
The coefficient is positive, so division recovers $t(F,W_{\rm in})$. Isolated vertices do not affect homomorphism densities, and edgeless tests have density $1$. Thus $\Phi_\sigma$ is injective with a continuous inverse. Its image is compact because the inner extremal space is compact.
\end{proof}

At positive gap, apply the single-edge case of \eqref{eq:readback} to the maps in \Cref{clm:phasemaps}. On the image of $\Phi_\sigma$, the test $T^{(\sigma)}$ has density $c_{r,\beta}$, because the normalized inner edge density is $z$. The inner hypergraph with the other phase label is empty, so its test has density $0$. Thus $\Sep=c_{r,\beta}$ on the image of $\Phi_1$ and $\Sep=-c_{r,\beta}$ on the image of $\Phi_2$, with the positive constant in \eqref{eq:sepconstant}. The two images are therefore disjoint and give the asserted compact phases.
\end{proof}

\section{Finite forbidden families and ordinary stability}\label{sec:transfer}
We now construct a finite forbidden family whose near-extremal free $r$-graphs are close to blowups of the compiled legal models. Root elimination and pair covering determine its density. To transfer the weighted structural conclusions, we prove vertex extendability and apply the degree-stability theorem of Liu, Mubayi and Reiher.

Every finite edgeless $r$-graph belongs to $\cD_\beta$, since the empty partial rooting extracts no records and is admissible. Failure of admissibility for a fixed partial rooting is witnessed by an edge with an incorrect profile, conflicting labels, or a locally illegal system of records on at most $Q$ data vertices. An edge with an incorrect profile has $r\le |M|+Q$ vertices. For conflicting labels or local illegality, retain an edge declaring each relevant record. Together these edges use at most $Q$ data vertices and $|M|$ roots. In any subgraph containing these retained edges, an extension of the rooting that assigns no new root to their data vertices retains the same edge profiles and extracted records. The failed constraint therefore persists. Thus \Cref{lem:unroot} gives an effective finite basis $\cA_\beta$ with witness bound $|M|+Q$.

To describe ordinary blowups of the rooted models, we now allow one vertex class for each root label. A \emph{canonical presentation} of an ordinary $r$-graph consists of a data set and $|M|$ possibly empty root classes that together partition its vertex set. No edge contains two vertices from the same root class. Identifying the vertices in each nonempty root class, while keeping the data vertices distinct, gives a quotient that is admissible with the prescribed root labels. Write $\mathfrak C_\beta$ for the class of $r$-graphs admitting such a presentation.

For a homomorphism from $\mathscr R$ into a canonical $r$-graph, first compose with the map to its quotient. Since $\mathscr R$ covers pairs, this composition cannot identify two of its vertices. Thus \Cref{lem:rigidity} applies to canonical targets as well.

Restricting a canonical presentation to a subgraph shows that $\mathfrak C_\beta$ is monotone. For a blowup, keep all root clones in their original root classes and clone the data records as in Section~\ref{sec:limits}, where local legality under cloning was verified. Thus the blowup is canonical. If $H\to G$ with $G\in\mathfrak C_\beta$, embed $H$ in a blowup of $G$ by mapping the vertices in each fiber to distinct clones of their common image. Monotonicity then gives $H\in\mathfrak C_\beta$, proving closure under inverse homomorphisms.

Given $H\in\mathfrak C_\beta$, form the quotient of a canonical presentation and complete its root supports to obtain $D\in\cD_\beta$. Replace each root vertex of $D$ by its original root class and keep the data vertices unchanged, obtaining $H^{\rm cmp}\supseteq H$ on the same vertex set. Fix a probability weighting $\x$ on $H$. Obtain $\bar\x$ by summing within each root class and retaining each data vertex's weight. Then
\begin{equation}\label{eq:completecompress}
 \lambda(H;\x)\le\lambda(H^{\rm cmp};\x)
       =\lambda(D;\bar\x)\le\Lambda,\qquad
 t(F,H^{\rm cmp};\x)=t(F,D;\bar\x).
\end{equation}
It follows that $\Lambda(\mathfrak C_\beta)=\Lambda(\cD_\beta)=\Lambda$.

Pair-covering extensions determine the density. For vertex extendability, we must also exclude vertices that cannot be assigned a canonical role. We forbid all noncanonical hypergraphs up to an order $s_\beta$ large enough to include both kinds of witnesses. Let $R_Q\coloneqq\max\{2,(g+3)\binom Qg\}$ and set
\begin{equation}\label{eq:cutoff}
 \begin{gathered}
 s_\beta\coloneqq\max\left\{\max_{A\in\cA_\beta}
 [v(A)+(r-2)\tbinom{v(A)}2],\ r^2(|M|+1)R_Q+2\right\},\\
 \cG_{r,\beta}\coloneqq\{F\colon v(F)\le s_\beta,\ F\notin\mathfrak C_\beta\}.
 \end{gathered}
\end{equation}
To construct $\cG_{r,\beta}$, enumerate all $r$-graphs of order at most $s_\beta$ and test all assignments of their vertices to the data role and the $|M|$ root roles. The explicit order bound and the bounded legality tests make these operations primitive recursive for fixed $r$.

If $F\in\cG_{r,\beta}$ and $F'$ is a homomorphic image of $F$, then $F'\notin\mathfrak C_\beta$; otherwise inverse homomorphism closure would give $F\in\mathfrak C_\beta$. Since $v(F')\le v(F)\le s_\beta$, we have $F'\in\cG_{r,\beta}$. Thus the output is homomorphism closed. All canonical $r$-graphs are free.

Every pair-covering extension of an $A\in\cA_\beta$ is noncanonical, since the quotient map would otherwise embed $A$ into the admissible quotient, contradicting $\cD_\beta=\Forb(\cA_\beta)$. By \eqref{eq:cutoff}, every such extension has order at most $s_\beta$, so $\bigcup_{A\in\cA_\beta}\Ext(A)\subseteq\cG_{r,\beta}$. Canonical blowups give the lower bound, and \Cref{lem:pairtransfer} gives the upper bound in
\[
 r!\Lambda\le\pi(\cG_{r,\beta})
 \le\pi\left(\bigcup_{A\in\cA_\beta}\Ext(A)\right)=r!\Lambda.
\]

Together with \eqref{eq:defectglobal}, this gives
\begin{equation}\label{eq:structrationalbaseline}
 \pi(\cG_{r,\beta})=r!\chi_r S(z).
\end{equation}
In particular, $\pi(\cG_{r,\beta})=r!\chi_r L/g!$ if and only if $z=\tau$.

We next show that every symmetrized $\cG_{r,\beta}$-free $r$-graph is canonical. An $r$-graph is \emph{symmetrized} if any two vertices that lie in no common edge have the same link. For a symmetrized $\cG_{r,\beta}$-free $r$-graph, identify vertices having the same link~\cite{LMR2023}. The quotient therefore covers pairs, and choosing one representative from each class realizes it as a subgraph of the original $r$-graph. A copy of $A\in\cA_\beta$ in the quotient, together with edges covering its pairs, would give a forbidden extension. Thus the quotient belongs to $\cD_\beta$, and the original $r$-graph is canonical because it is a blowup of this quotient.

\begin{proof}[Proof of \Cref{thm:structural}]
\Cref{thm:loss} identifies the weighted extremal space. We first establish vertex extendability for the canonical class, then transfer the weighted conclusions to ordinary $\cG_{r,\beta}$-free $r$-graphs, and finally round the class sizes.

The witnesses excluding the possible roles of a new vertex involve only boundedly many old vertices. We construct a common frame and replacement frames through their root vertices. Together with \Cref{lem:rigidity}, these frames will fix the old root roles up to one common allowed permutation.
\begin{claim}\label{clm:commonframe}
Fix $r$ and $\beta$. Let $H_N$ be canonical $N$-vertex $r$-graphs with specified canonical presentations and $\delta(H_N)\ge(r\Lambda-o(1))N^{r-1}$. For every sequence $S_N\subseteq V(H_N)$ with $\sup_N|S_N|<\infty$, $H_N-S_N$ contains a copy of $\mathscr R$ respecting the specified root and data roles for all sufficiently large $N$. For every root vertex $y\in S_N$, replacing the representative of its root class by $y$ still gives a copy of $\mathscr R$ in $H_N$.
\end{claim}
\begin{proof}
Write $H_0\coloneqq H_N$ and $S\coloneqq S_N$, and complete the root supports in the specified canonical presentation of $H_0$ to obtain $H_1$. Give every vertex of $H_1$ weight $1/N$, and let $(D,\bar\x)$ be the completed weighted quotient from \eqref{eq:completecompress}. Summing degrees in $H_0$ and applying the Lagrangian bound to $D$ gives
\[
 \frac{N\delta(H_0)}r\le |H_0|\le |H_1|
       =N^r\lambda(D;\bar\x)\le\Lambda N^r.
\]
The minimum-degree bound therefore gives $|H_i|/N^r\to\Lambda$ for $i=0,1$. Apply \Cref{thm:loss} to the sequence of weighted quotients. The code-root classes have sizes $(3/(5r)+o(1))N$, the suspension classes have sizes $(1/r+o(1))N$, and decoding yields a neutral set $W$ of linear size.

To count the missing or changed records, let $C$ be the data vertex set and let $\mathcal Z$ consist of its $g$-sets whose record is missing or whose original and decoded labels differ. We have $|C|=(g/r+o(1))N$. The normalized data weights are all $1/|C|$, with erasure probability $\epsilon\to0$, so the disagreement bound in \Cref{lem:recordrepair} gives
\[
 \frac{g!|\mathcal Z|}{|C|^g}\le(g+1)\epsilon\longrightarrow0.
\]
Thus $|\mathcal Z|=o(N^g)$, and every $g$-subset of $W$ outside $\mathcal Z$ has the neutral record in $H_1$.

To obtain replacement frames, we uniformly bound the number of edges in $H_1\setminus H_0$ incident to each root vertex. In $H_1$ each root class consists of twins. For each code word, the corresponding product support is invariant under cyclic permutation of the five root roles within each block. Cyclically permuting the five asymptotically equal class sizes within a block changes each nonnegative link monomial by a factor $1+o(1)$. Write $n_i$ for the size of class $i$ in one 5-vertex block and $d_i$ for its common per-vertex degree. Then
\[
 \sum_{i=1}^5 n_id_i=3|H_1|,\qquad
 n_i=\left(\frac{3}{5r}+o(1)\right)N,\qquad
 d_i=(1+o(1))d_1.
\]
Using $|H_1|=(\Lambda+o(1))N^r$ and summing the five class sizes gives
\[
 d_i=(1+o(1))\frac{3|H_1|}{\sum_{j=1}^5n_j}
 =\frac{(3\Lambda+o(1))N^r}{(3/r+o(1))N}
 =(r\Lambda+o(1))N^{r-1}.
\]
A suspension class $V_j$ has per-vertex degree $|H_1|/|V_j|$, with the same asymptotic value. Let $\mathcal M\coloneqq H_1\setminus H_0$. Since $H_0\subseteq H_1$, we have $|\mathcal M|=|H_1|-|H_0|$, and for every root vertex $y$,
\[
 \deg_{\mathcal M}(y)=\deg_{H_1}(y)-\deg_{H_0}(y)
       \le\deg_{H_1}(y)-\delta(H_0).
\]
The edge counts and degree bounds above therefore give
\begin{equation}\label{eq:rootmissing}
 |\mathcal M|=o(N^r),\qquad
 \max_{y\text{ in a root class}}\deg_{\mathcal M}(y)=o(N^{r-1}).
\end{equation}

We now sample the common frame. Choose one uniformly random representative from each root class and $g+1$ uniformly random vertices of $W$, all choices being independent. We require all chosen vertices to lie outside $S$ and the data choices to be distinct; failure has probability $O(1/N)$ because $S$ is bounded and all sampling classes have linear size. We estimate every other failure under this original sampling law, without conditioning on these requirements.

Use the chosen vertices for the neutral frame $\mathscr R$. For every root vertex $y\in S$, also require a copy of $\mathscr R$ obtained by replacing the representative of $y$'s root class by $y$, while leaving all other sampled vertices fixed. The number of required frames and edges is bounded in terms of the fixed input, $r$, and $\sup_N|S_N|$. Every data $g$-set used by these frames must avoid $\mathcal Z$; the union bound gives failure probability $O(|\mathcal Z|/N^g)$. When its record is neutral, each required edge belongs to the completion $H_1$. A given common-frame edge belongs to $\mathcal M$ with probability $O(|\mathcal M|/N^r)$, since each possible image has probability $O(N^{-r})$. For a replacement edge through a fixed root $y$, only $r-1$ vertices remain random, so the corresponding bound is $O(\deg_{\mathcal M}(y)/N^{r-1})$.

The probability that any requirement fails is therefore at most
\[
 O\left(\frac{|\mathcal Z|}{N^g}\right)
 +O\left(\frac{|\mathcal M|}{N^r}\right)
 +\sum_{\substack{y\in S\\y\text{ in a root class}}}O\left(\frac{\deg_{\mathcal M}(y)}{N^{r-1}}\right)
 +O\left(\frac1N\right)=o(1).
\]
Thus one choice gives all the required frames inside $H_0$.
\end{proof}

\begin{claim}\label{prop:vertexextend}
For some $\zeta>0$ and all large $n$, an $\cG_{r,\beta}$-free $n$-vertex $r$-graph $H$ with $\delta(H)\ge(r\Lambda-\zeta)n^{r-1}$ and $H-v$ canonical is canonical.
\end{claim}
\begin{proof}
Suppose there are counterexamples with $\zeta\to0$ and $n\to\infty$, with the input and $r$ fixed. For each possible role of $v$, we will retain edges witnessing why that role cannot extend the canonical presentation of $H-v$, and force all old vertices in these witnesses to keep their roles. The resulting bounded noncanonical subgraph will contradict $\cG_{r,\beta}$-freeness.

Let $H_0\coloneqq H-v$ and $N\coloneqq n-1$. Deleting $v$ decreases every remaining degree by $O(N^{r-2})$, so $\delta(H_0)\ge(r\Lambda-o(1))N^{r-1}$. Fix a canonical presentation of $H_0$.

Consider each of the $|M|$ root roles and the data role of $v$. An incorrect edge profile is witnessed by the edge itself. Conflicting labels are witnessed by 2 edges declaring different labels on the same $g$-set. A local-legality failure on at most $Q$ data vertices is preserved by retaining one edge declaring each available label on each $g$-set, so at most $(g+3)\binom Qg$ edges suffice. Thus every failed role has a witness with at most $R_Q$ edges. For each role, at least one edge in its witness contains $v$, since otherwise the same constraint would already fail in the canonical $H_0$. The union $S_0$ of their old vertices therefore satisfies $|S_0|\le(|M|+1)(rR_Q-1)$.

To keep the old data vertices in their data roles, choose, for each such vertex $y\in S_0$, an edge $E_y\in H_0$ containing $y$; the positive minimum degree of $H_0$ guarantees such an edge. Set $S\coloneqq S_0\cup\bigcup_y E_y$, so $|S|\le r|S_0|$. Apply \Cref{clm:commonframe} to $H_0$ and this bounded set $S$. It gives a common frame outside $S$ and a replacement frame through each root vertex of $S$, as illustrated in Figure~\ref{fig:replacementframe}.

\begin{figure}[htbp]
\centering
\begin{tikzpicture}[x=1cm,y=1cm,>=Stealth,every node/.style={font=\small},
 root/.style={rectangle,draw,fill=white,inner sep=2.2pt},
 neutral/.style={circle,draw,fill=white,inner sep=1.7pt}]
 \foreach \shift/\name/\vertex in {0/Common frame/u,7.3/Replacement frame/y} {
  \begin{scope}[xshift=\shift cm]
   \node at (2.1,1.5) {\name};
   \draw[rounded corners=3pt,black!45] (-0.35,-1.05) rectangle (4.55,1.05);
   \node[root] at (0,0.45) {};
   \node at (1.05,0.45) {$\cdots$};
   \node[root,fill=black!25,label=above:$\vertex$] at (2.1,0.45) {};
   \node at (3.15,0.45) {$\cdots$};
   \node[root] at (4.2,0.45) {};
   \foreach \xx in {0.5,1.3,2.1,2.9,3.7} \node[neutral] at (\xx,-0.55) {};
   \node at (2.1,-1.45) {$g+1=5$ neutral data vertices};
  \end{scope}
 }
 \draw[->] (4.8,0.45)--(7.05,0.45);
 \node[align=center] at (5.9,1.15) {replace $u$\\by $y$};
 \node[align=center] at (5.8,-2.1) {The other $|M|-1$ root representatives and all data vertices are shared.};
\end{tikzpicture}
\caption{The common frame and a replacement frame through a root vertex $y\in S$. Here $u$ is the chosen representative of the same root class. Corresponding unshaded vertices are identical in the two panels, and omitted root representatives are unchanged. Edges are suppressed to display the overlap. Rigidity forces the two root permutations to agree on the shared roles, and hence on the remaining role.}
\label{fig:replacementframe}
\end{figure}

Consider the subgraph formed by the retained failure witnesses, the edges $E_y$, and all these frames. In any canonical presentation of this subgraph, \Cref{lem:rigidity} fixes the root roles of the common frame modulo some $\pi\in\mathscr H_r$. A replacement frame agrees with it on all but one formal root, so their root permutations agree on $|M|-1$ roles. As both are permutations of the same $|M|$ roles, they agree on the last role as well. Hence every root vertex of $S$ has its original role modulo the same $\pi$.

The old root vertices of each $E_y$ occupy exactly 3 roles in each 5-vertex block and all suspension roles, modulo $\pi$. Giving any other vertex of $E_y$ a root role would repeat a role or exceed the prescribed number in its block, so $y$ remains a data vertex. Every old vertex in each failure witness therefore retains its role, up to the same permutation preserving all code pools, and each retained constraint still excludes the corresponding role for $v$. Thus this subgraph admits no canonical presentation. This contradicts \eqref{eq:cutoff} and $\cG_{r,\beta}$-freeness, since its order is at most
\[
 r|S_0|+|M|+g+2\le r^2(|M|+1)R_Q+2.\qedhere
\]
\end{proof}
Here $\zeta$ and the stability constants below may depend on $\beta$ and $r$.

We now transfer the weighted classification to ordinary free $r$-graphs.
\begin{claim}\label{thm:ordinarytransfer}
The following statements hold.
\begin{enumerate}[beginpenalty=10000]
\item Every asymptotically extremal $\cG_{r,\beta}$-free sequence is $o(n^r)$-close, on the same vertex sets, to blowups of $r$-graphs compiled from complete legal objects.
\item For every such sequence, deleting $o(n)$ vertices from each member makes it a subgraph of a blowup of an $r$-graph compiled from a complete legal object.
\item We have $\cE(\cG_{r,\beta})=\mathfrak L(\cD_\beta)$.
\end{enumerate}
\end{claim}
\begin{proof}
We first obtain a canonical presentation after deleting $o(n)$ vertices, and then repair its data records to a complete legal object. The first step uses degree stability; the second uses \Cref{lem:recordrepair}.

Apply the degree-stability theorem of Liu, Mubayi and Reiher~\cite[Definitions 1.5--1.6 and Theorem 1.7]{LMR2023} to the forbidden family $\cG_{r,\beta}$ and the class $\mathfrak C_\beta$ of canonical $r$-graphs. The hypotheses hold for the following reasons.
\begin{enumerate}
\item The family $\cG_{r,\beta}$ has positive Tur\'an density and is closed under homomorphic images, hence is blowup invariant.
\item The class $\mathfrak C_\beta$ is hereditary and consists of $\cG_{r,\beta}$-free $r$-graphs.
\item Every symmetrized $\cG_{r,\beta}$-free $r$-graph belongs to $\mathfrak C_\beta$.
\item \Cref{prop:vertexextend} gives vertex extendability for $\mathfrak C_\beta$.
\end{enumerate}
The degree threshold is $\pi(\cG_{r,\beta})/(r-1)!=r\Lambda$.

\emph{From degree stability to a canonical presentation.} Degree stability gives $\gamma>0$ and $N_0$ such that every $\cG_{r,\beta}$-free $r$-graph $H$ of order $N\ge N_0$ with $\delta(H)\ge(r\Lambda-\gamma)N^{r-1}$ belongs to $\mathfrak C_\beta$.

Start with an $n$-vertex $\cG_{r,\beta}$-free member of an asymptotically extremal sequence. At each step, let $H$ be the current $r$-graph and $N$ its order, and delete a vertex $y$ if $\deg_H(y)<(r\Lambda-\gamma)N^{r-1}$. To bound the number of deletions, we track the difference $|H|-\Lambda N^r$, which is always nonpositive by \eqref{eq:lagdensity}, since $\cG_{r,\beta}$-freeness is equivalent to $\cG_{r,\beta}$-hom-freeness. It starts at $-o(n^r)$. Since $N^r-(N-1)^r=rN^{r-1}+O(N^{r-2})$, the increase at each deletion satisfies
\[
 \Lambda[N^r-(N-1)^r]-\deg_H(y)
 \ge\gamma N^{r-1}-O(N^{r-2})\ge\frac\gamma2N^{r-1}
\]
for large $N$. While $N\ge n/2$, each deletion increases this difference by an amount of order $n^{r-1}$. Since the difference starts at $-o(n^r)$ and remains nonpositive, only $o(n)$ deletions are possible. In particular, reaching $N=n/2$ is impossible. The process therefore stops at order $n-o(n)$ with the minimum degree required for degree stability, and the remaining $r$-graph is canonical. This is the usual degree-to-vertex stability argument (cf.\ \cite[Fact 2.5]{LMR2023}).

To obtain the approximation in part \textup{(i)} on the original vertex set, restore the deleted vertices as isolated data vertices and complete the root supports in the canonical presentation. The restored $r$-graph still has $(\Lambda-o(1))n^r$ edges, and its completion is $\cG_{r,\beta}$-free, so the completion adds only $o(n^r)$ edges by \eqref{eq:lagdensity}. At this stage only code supports have been supplied; the extracted data records may still be partial.

\emph{Repairing the data records.} Compress the root classes, with weights equal to their proportions, and give data vertices weight $1/n$. By \eqref{eq:completecompress}, the edge-polynomial values tend to $\Lambda$. Thus \Cref{thm:loss} gives $|C|/n\to g/r$ and erasure probability $\epsilon\to0$ for the data set $C$. Applying \Cref{lem:recordrepair} with weights $1/|C|$ yields $V'\subseteq C$ with $|C\setminus V'|\le\epsilon|C|=o(n)$ and a legal realization preserving every available record on $V'$.

Make $C\setminus V'$ neutral and label every $g$-set using this legal realization. By the definition of $\epsilon$, at most $\epsilon|C|^g/g!$ records are missing, and at most $|C\setminus V'|\binom{|C|-1}{g-1}$ records have an underlying $g$-set meeting $C\setminus V'$. Both bounds are $o(n^g)$. Only these records can change, and each affects at most $\binom n{r-g}$ edges. Recompile the records using the existing root classes. The resulting blowup differs by $o(n^r)$ edges on the same vertex set.

For the subgraph conclusion in part \textup{(ii)}, omit both the vertices removed in the degree-pruning step and the data vertices discarded by record repair. Their union has size $o(n)$, and every edge on the remaining vertex set belongs to the reconstructed legal completion.

To prove the forward inclusion in \textup{(iii)}, take an asymptotically extremal $\cG_{r,\beta}$-free sequence with density-vector limit $W$, and apply \textup{(i)}. A union bound over the edges of each fixed test shows that the $o(n^r)$ edge changes preserve its limiting density. Compress the root classes of each resulting blowup, assigning each root its class proportion and each data vertex weight $1/n$. The weighted quotients belong to $\cD_\beta$ and, by \eqref{eq:completecompress}, have the same homomorphism densities as the blowups. Their edge-polynomial values tend to $\Lambda$, so $W\in\mathfrak L(\cD_\beta)$.

Conversely, fix $W\in\mathfrak L(\cD_\beta)$ and choose weighted $D_i\in\cD_\beta$ whose density vectors approach $W$. Approximate their weights rationally and take the corresponding blowups. Choose the approximations and blowup orders so that the first $i$ test densities change by at most $1/i$ and the orders tend to infinity. These blowups are canonical and therefore $\cG_{r,\beta}$-free. Their density vectors approach $W$, and their edge densities approach $r!\Lambda$, so $W\in\cE(\cG_{r,\beta})$. This proves the reverse inclusion and hence the equality of limit spaces.
\end{proof}

The equality of limit spaces in \Cref{thm:ordinarytransfer} transfers the classification in \Cref{thm:loss} to $\cE(\cG_{r,\beta})$. The readback identity \eqref{eq:readback} gives separator value $0$ at zero gap and values $\pm c_{r,\beta}$ on the two positive-gap phases.

We now round the class sizes in the sequence of blowups supplied by \Cref{thm:ordinarytransfer}. By \Cref{thm:loss}, each class size differs by $o(n)$ from its prescribed target. The target sizes are $3n/(5r)$ for each code-root class, $n/r$ for each suspension class, and $a_\delta n,w_\delta n$ for the active and neutral classes, respectively, with $a_\delta,w_\delta$ as in \eqref{eq:optimalmasses}. We choose integer roundings summing to $n$ while preserving the edge approximation and, at positive gap, inner extremality. At zero gap the active set has $o(n)$ vertices; make them neutral and adjust the root and data class sizes to the prescribed integer roundings.

At positive gap the model uses a single phase label and the inner $3$-graph is asymptotically optimal. First adjust the root class sizes to their prescribed integer roundings by moving $o(n)$ vertices between root classes and the neutral class, keeping the whole vertex set fixed. If the active class is too large, remove the excess vertices from the inner $3$-graph and reassign them to the neutral class. If it is too small, reassign the required number of neutral vertices as clones of a fixed surviving inner vertex. There are linearly many neutral vertices available. Restriction and cloning preserve inner freeness.

Write $U,J$ and $U',J'$ for the active sets and inner $3$-graphs before and after rounding. Only edges meeting $U\triangle U'$ can change, giving $\bigl||J'|-|J|\bigr|\le |U\triangle U'|\binom{n-1}{2}=o(n^3)$. Since $|U|/n,|U'|/n\to a_\delta>0$, we have $|U|/|U'|\to1$ and
\[
 \frac{6|J'|}{|U'|^3}
 =\frac{6|J|}{|U|^3}\left(\frac{|U|}{|U'|}\right)^3
  +\frac{6(|J'|-|J|)}{|U'|^3}=z+o(1),
\]
using $6|J|/|U|^3\to z$ and $|J'|-|J|=o(n^3)$ in the last equality. In both gap cases, rebuilding the affected ordinary edges changes only edges meeting the $o(n)$ reassigned vertices, so the cost is $o(n^r)$ edge changes.

The two disjoint nonempty compact phases exclude uniqueness and connectedness. A fixed finite profile has a unique homomorphism-density limit, and $o(n^r)$ edge changes preserve every fixed test density. Thus stability with respect to a single finite profile would force the extremal space to be a singleton, contrary to the two distinct separator values. At zero gap, the extremal space is a singleton and hence connected, and stability holds with respect to the fixed neutral profile.
\end{proof}

\section{Computational consequences and logical independence}\label{sec:consequences}
\subsection{Finite comparisons and witness bounds}
Effective approximation allows a search for a strict density inequality. We first classify the complexity of exact comparison for finite forbidden families.

A relation is \emph{$\Sigma^0_1$} if it has the form $\exists s\,R(e,s)$ for a decidable predicate $R$. Equivalently, there is a search that halts exactly when the relation holds. A \emph{$\Pi^0_1$} relation has the form $\forall s\,R(e,s)$, so its complement has such a search. A \emph{computable many-one reduction} from a relation $P$ to a relation $Q$ is a total computable map $f$ such that $P(e)$ holds if and only if $Q(f(e))$ holds. A relation is \emph{complete} for a class if it belongs to the class and every relation in that class admits such a reduction to it. Tuples of finite inputs are encoded by the conventions in Subsection~\ref{sec:logic}.

\begin{theorem}\label{thm:finitecomparison}
The following statements hold.
\begin{enumerate}[beginpenalty=10000]
\item Given $r\ge2$, a finite family $\cF$ of $r$-graphs, and a rational $\eps>0$, a single algorithm returns a rational $a$ satisfying $|a-\pi(\cF)|<\eps$.
\item For fixed $r\ge3$, the relations $>,<,\ne$ between two such densities are $\Sigma^0_1$-complete, and $\le,\ge,=$ are $\Pi^0_1$-complete. For $r=3$, all these classifications remain valid when both families have one member. For $>,\ne,=,\le$ the second family can be fixed to $\cF_{r,\beta_\infty}$; for the reverse order relations the first can be fixed.
\end{enumerate}
\end{theorem}
\begin{proof}
Equation \eqref{eq:approx} gives intervals of any prescribed rational width. To semidecide $\pi(\cF)>\pi(\cH)$, compute successively narrower certified intervals for the two densities and halt when the lower endpoint for $\cF$ exceeds the upper endpoint for $\cH$. This happens exactly when the strict inequality holds. Searching in the other direction handles $<$, and running both searches in parallel handles $\ne$. Taking complements shows that the relations $\le,\ge,=$ belong to $\Pi^0_1$.

For hardness, \Cref{thm:numerical} reduces halting to strict excess over $\tau_r$ and nonhalting to equality with $\tau_r$. Interchanging the two families gives hardness for the reversed inequalities. At $r=3$, \Cref{thm:singleton} gives the same reductions using only single forbidden $3$-graphs, with $P_{\beta_\infty}$ as the fixed comparison object.
\end{proof}
For a finite family $\cF$ such that every blowup of $K_{*,r}$ is $\cF$-free, \eqref{eq:finitecertificate} shows that $\pi(\cF)=\tau_r$ if and only if no finite $\cF$-hom-free $r$-graph $G$ satisfies \eqref{eq:balancedimprovement}. This is a $\Pi^0_1$ assertion.

If a positive rational $\eta$ is supplied with the promise that $\pi(\cF)=\tau_r$ or $\pi(\cF)\ge\tau_r+\eta$, these 2 cases can be distinguished effectively. To do so, approximate $\pi(\cF)$ to error less than $\eta/3$ and compare the result with $\tau_r+\eta/2$.

On halting inputs, the minimum order of an improvement witness is comparable to the halting-history area, and the density gain is bounded below by a constant times the inverse cube of that area. Let
\[
 \omega_r(\beta)\coloneqq\min\{v(G)\colon G\text{ is }\cF_{r,\beta}\text{-hom-free},\
                              \lambda_r(G)>\tau_r\},
\]
with empty minimum equal to infinity. We use a variant of the Busy-Beaver function introduced by Rad\'o~\cite{Rado1962}. Let $\BB_{\mathsf U}(m)$ and $\BB_r^{\rm Tur}(m)$ be the maxima of $T(\beta)$ and $\omega_r(\beta)$, respectively, over halting words of length at most $m$, with an empty maximum taken to be zero.
\begin{theorem}\label{thm:witnessbounds}
For fixed $r\ge3$, the following statements hold.
\begin{enumerate}[beginpenalty=10000]
\item The value $\omega_r(\beta)$ is finite exactly on halting inputs. There are effective constants $C_r,c_r>0$ such that, on every halting input, writing $A\coloneqq A_{\rm hist}(\beta)$, we have
\begin{align}
 A&\le\omega_r(\beta)\le q_3+m_3A+r-3
       \le C_r(T(\beta)+|\beta|+1)^2,\label{eq:witnessbounds}\\
 \pi(\cF_{r,\beta})-\tau_r
       &\ge\frac{27r!t_3^3}{r^rm_3^3A^3}
        \ge c_r(T(\beta)+|\beta|+1)^{-6}.\label{eq:densitygap}
\end{align}
\item With the same constant $C_r$, the two Busy-Beaver functions satisfy
\[
 \BB_{\mathsf U}(m)\le\BB_r^{\rm Tur}(m)
       \le C_r(\BB_{\mathsf U}(m)+m+1)^2.
\]
\item There is no total computable upper bound for the finite values $\omega_r(\beta)$ on halting inputs as a function of the input length or the length of the code for the output family $\cF_{r,\beta}$. There is also no total computable function assigning a positive rational number $\eta(\beta)$ to each input $\beta$ such that $\eta(\beta)$ is a lower bound for the density gap on every halting input.
\end{enumerate}
\end{theorem}
\begin{proof}
By \Cref{prop:v3-certificates}, every improving template has at least $A=A_{\rm hist}(\beta)$ vertices, and a halting history gives one on $q_3+m_3A+r-3$ vertices. The same proposition gives a certificate $G$ for which $r!|G|/v(G)^r$ exceeds $\tau_r$ by $27r!t_3^3/(r^rm_3^3A^3)$.

Since $S(\beta)\le T(\beta)+1$ and $T(\beta)\ge1$, the definition \eqref{eq:historyarea} gives
\[
 W_\beta+2\le\max\{T(\beta)+3,\ |\beta|+4\}
            \le2(T(\beta)+|\beta|+1).
\]
Since $T(\beta)+2\le\frac32(T(\beta)+|\beta|+1)$, multiplication gives $A\le3(T(\beta)+|\beta|+1)^2$, and hence $A^{-3}\ge\frac1{27}(T(\beta)+|\beta|+1)^{-6}$. Substituting these bounds into the template order and density gain proves \cref{eq:witnessbounds,eq:densitygap} with $C_r=3m_3+q_3+r-3$ and $c_r=r!t_3^3/(r^rm_3^3)$, respectively. Since $T(\beta)\le A$, maximizing \eqref{eq:witnessbounds} over halting words of length at most $m$ gives the Busy-Beaver comparison.

For the noncomputability assertions, a computable witness bound would bound halting time by \eqref{eq:witnessbounds}, deciding halting by bounded simulation. A bound depending on the length of the code for the output family has the same consequence because that length is computable from $\beta$. If a computable positive gap bound $\eta(\beta)$ existed, approximation to error less than $\eta(\beta)/3$ and comparison with $\tau_r+\eta(\beta)/2$ would likewise decide halting.
\end{proof}

\emph{Output size.} The preceding bounds concern hypergraphs witnessing a density improvement. We now bound the order of the forbidden $3$-graph $P_\beta$ produced from the input $\beta$. Let $b=b(\beta)\le b_0(|\beta|+1)$ be the local witness bound in \Cref{prop:direct-finite-basis}, and let $H$ be the fixed guard. Define
\[
 M\coloneqq b\sum_{j=0}^{q_3}(q_3b)^j+v(H),\qquad
 d\coloneqq(M+1)2^{\binom M3},\qquad
 N\coloneqq\left(M+\binom M2\right)^d.
\]
After removing duplicate isomorphism types, the guarded basis has at most $d$ members, each on at most $M$ vertices. A canonical extension has at most $M+\binom M2$ vertices. Consequently $v(P_\beta)\le N$. In particular, for an effective constant $C>0$,
\[
 v(P_\beta)\le 2^{\,2^{C(|\beta|+1)^{3(q_3+1)}}}.
\]
For $r\ge4$, every member of $\cQ(\{P_\beta\})$ has at most $N$ vertices. Let
\[
 N_r\coloneqq r+\binom r3\max\{r,N\},\qquad
 M_r\coloneqq N_r+(r-2)\binom{N_r}2.
\]
The suspension basis and pair-covering extensions then give
\[
 \max_{F\in\cF_{r,\beta}}v(F)\le M_r,\qquad
 |\cF_{r,\beta}|\le(M_r+1)2^{\binom{M_r}r}.
\]
These are primitive recursive bounds in $\beta$ for every fixed $r$, and are computable uniformly in $r$.

\subsection{Uniform formalization and logical independence}\label{sec:logic}
The incompleteness argument uses a single ZFC proof that works for every machine input. We record the uniform numerical and structural statements needed to transfer a Rosser sentence to a forbidden family.

A function is \emph{primitive recursive} if it is obtained from the zero, successor and projection functions by composition and primitive recursion. We use finite codes, canonical representatives obtained by permutation search, and finite computation records~\cite{Rogers1967,Soare1987}. A hypergraph code records a nonnegative order and a finite edge list, with duplicate edges removed; a weighted hypergraph is evaluated only on a nonempty vertex set.

A malformed finite-family code denotes the empty family. When processing a valid family list, discard malformed hypergraph codes and edgeless hypergraphs. Word codes represent all finite binary words, including the empty word; malformed word codes are sent to $\beta_\infty$. A certificate conversion returns a fixed failure symbol on a malformed certificate or a failed bounded check. All finite parsing and bounded searches are primitive recursive.

Fix $r\ge3$, and write $\tau_r=p_r/b_r$. Let $\mathsf{Num}_r(\beta)$ be the code of the finite family $\cF_{r,\beta}$ produced by the numerical construction.

Let $\Improve_r(f,h)$ mean that $h$ codes a nonempty $\cF_f$-hom-free $r$-graph $G$ with $b_rr!|G|>p_rv(G)^r$. This predicate and the assertion that $K_{*,r}$ is $\cF_f$-hom-free are primitive recursive. If $K_{*,r}$ is $\cF_f$-hom-free, \eqref{eq:finitecertificate} identifies the density equality $\pi(\cF_f)=\tau_r$ with $\forall h\,\neg\Improve_r(f,h)$.

We next formalize the structural conclusions using density-vector limits. For $r\ge r_{\rm str}$, the structural construction takes a finite family of $3$-graphs. Before applying it, delete the isolated vertices from each member and take the homomorphic closure of the resulting family. The construction is defined without a density promise.

For a finite profile $p$ with role proportions $\alpha_i$, let $\mathfrak B_p(V)$ be the finite set of complete labeled realizations of $p$ on $V$ in which role class $i$ has size $\lfloor\alpha_i|V|\rfloor$ or $\lceil\alpha_i|V|\rceil$, with the class sizes summing to $|V|$. For a family code $f$, define the following five assertions.
\begin{enumerate}
\item $S_1(f)$ asserts that $\cE(\cF_f)$ is a singleton.
\item $S_2(f)$ asserts that $\cE(\cF_f)$ is connected.
\item $S_3(f)$ asserts that the separator on $\cE(\cF_f)$ takes both signs.
\item $S_4(f)$ asserts that $\cE(\cF_f)$ is the union of two nonempty regions $\Sep\le-a$ and $\Sep\ge a$ for some rational $a>0$.
\item $S_5(f)$ asserts the existence of a finite profile $p$ whose associated homomorphism-density limit belongs to $\cE(\cF_f)$ and for which
\begin{equation}\label{eq:finitarystability}
\begin{aligned}
 &\forall\epsilon\in\Q_{>0}\ \exists\eta\in\Q_{>0}\ \exists N\ \forall G,\\
 &G\text{ is }\cF_f\text{-free},\ v(G)\ge N,\quad
            |G|\ge\ex(v(G),\cF_f)-\eta v(G)^r\\
 &\hspace{12mm}\Longrightarrow\
   \exists B\in\mathfrak B_p(V(G))\colon\ |G\triangle B|\le\epsilon v(G)^r.
\end{aligned}
\end{equation}
\end{enumerate}
To compare \eqref{eq:finitarystability} with the sequential definition of stability, fix a profile $p$. If the quantified estimate fails, there is an $\epsilon_0\in\Q_{>0}$ such that, for every $j\ge1$, some $\cF_f$-free $G_j$ of order $n_j\ge j$ satisfies
\[
 \ex(n_j,\cF_f)-|G_j|\le \frac{n_j^r}{j},\qquad
 \min_{B\in\mathfrak B_p(V(G_j))}|G_j\triangle B|>\epsilon_0n_j^r.
\]
This is an asymptotically extremal sequence that contradicts the sequential stability condition for $p$. The quantified estimate implies sequential stability directly, so the two definitions are equivalent.

\begin{proof}[Proof of \Cref{thm:independence}]
We first establish the numerical and structural equivalences uniformly in their finite inputs inside ZFC. We then apply them to a Rosser sentence for the given theory.

\begin{claim}
\label{prop:numericalbridge}\label{prop:structuralbridge}
For fixed $r$ the output maps and the numerical certificate conversions are primitive recursive, with ZFC proofs of totality and correctness. The numerical maps are total and computable uniformly in $r,\beta$. ZFC proves that, for every $\beta$, $K_{*,r}$ is $\cF_{r,\beta}$-hom-free and
\begin{equation}\label{eq:internalnumerical}
 \exists t\,[\mathsf U(\beta)\text{ halts within }t\text{ steps}]\quad\text{if and only if}\quad
 \exists h\,\Improve_r(\mathsf{Num}_r(\beta),h).
\end{equation}
For $r\ge r_{\rm str}$, let $f$ code a finite family of $3$-graphs and let $d$ code the output of the structural construction applied to $\cF_f$. Uniformly for codes $f$ such that $K_{*,3}$ is $\cF_f$-hom-free, ZFC also proves
\begin{equation}\label{eq:internalstructural}
 \begin{aligned}
 \pi(\cF_f)=\tau_3&\quad\text{if and only if}\quad S_i(d)\quad(i=1,2,5),\\
 \pi(\cF_f)>\tau_3&\quad\text{if and only if}\quad S_i(d)\quad(i=3,4).
 \end{aligned}
\end{equation}
The structural test hypergraphs and the neutral profile depend only on $r$.
\end{claim}
\begin{proof}
The numerical bounds are explicit in \Cref{lem:unroot,lem:final-singleton-effective-basis,lem:v3-effective-suspension,prop:v3-certificates} and in the preceding output-size estimates. Finite overlap unions, canonical extensions, categorical products, and homomorphic images are computed by bounded operations on finite vertex and edge lists. The parameters involving square roots use integer comparisons to choose the least admissible integer. In the structural construction, the vertex bound for local tests and the order bound for forbidden hypergraphs are given explicitly as functions of the code for the finite inner family. These bounds use integer arithmetic, powers, finite sums and finite maxima. Hypergraph enumeration and the tests of root and type assignments within these bounds are primitive recursive. Their composition therefore gives the output maps, with ZFC-provable totality and the prescribed defaults on malformed inputs.

Given a halting computation record, the forward conversion in \Cref{prop:v3-certificates} produces a certificate with uniform vertex weights and order $D_rA_{\rm hist}(\beta)$. Its verification establishes hom-freeness and the integer inequality defining $\Improve_r$. Conversely, for a finite $r$-graph $G$ passing these two checks, the same proposition gives $T(\beta)\le v(G)$, so simulation for $v(G)$ steps produces a halting record. The verifications use the fixed transition table and the stated size bounds, with $\beta$ and the finite record or hypergraph code as arbitrary variables. Carrying out these verifications in ZFC gives the single universally quantified equivalence \eqref{eq:internalnumerical}.

For the structural construction, the assumed hom-freeness of $K_{*,3}$ gives $\tau_3\le z\le1$. The record reconstruction, erasure estimate, optimization, frame identification and ordinary transfer were proved throughout this range for an arbitrary finite inner family. All are ZFC arguments, including the cited stability theorem~\cite{LMR2023}. The zero-gap case yields $S_1,S_2,S_5$ and the positive-gap case yields $S_3,S_4$ while excluding stability with respect to a single finite profile. These exhaustive alternatives prove \eqref{eq:internalstructural}.
\end{proof}

For the Rosser application, let $a$ be a computable axiom enumerator for $T$. Each nonlogical axiom in a proof is accompanied by a finite computation record showing that $a$ outputs it. The primitive recursive predicate $\mathsf{Prf}_a(u,\ulcorner\varphi\urcorner)$ checks these records and the inference steps in the proof of $\varphi$ coded by $u$, where $\ulcorner\varphi\urcorner$ denotes the sentence code. The effective Rosser construction~\cite{Rosser1936,BoolosBurgessJeffrey2007} gives a sentence $\rho_a$ equivalent to
\[
 \forall u\,[\mathsf{Prf}_a(u,\ulcorner\rho_a\urcorner)
       \Longrightarrow\exists v\le u\,
              \mathsf{Prf}_a(v,\ulcorner\neg\rho_a\urcorner)].
\]
Under consistency, neither $\rho_a$ nor its negation is provable. Moreover $\rho_a$ is true, since no proof of $\rho_a$ exists.

The predicate inside the universal quantifier of the Rosser sentence is primitive recursive, since its search over $v\le u$ is bounded. We convert the sentence into a machine input by the following uniform search construction. For a valid code $d$ of a primitive recursive predicate $R_d(e,n)$, let $\beta(d,e)$ encode the program that searches successively for an $n$ with $\neg R_d(e,n)$. Producing this code only inserts $d,e$ into a fixed search program.

The predicate is evaluated only when the program runs. Finite evaluation trees record the intermediate values; their validity can be checked by a primitive recursive predicate. Induction on the definition encoded by $d$ proves that each evaluation of $R_d(e,n)$ terminates. ZFC proves the equivalence of the following three assertions.
\[
 \forall n\,R_d(e,n),\qquad
 \mathsf U(\beta(d,e))\text{ does not halt},\qquad
 \pi(\cF_{r,\beta(d,e)})=\tau_r.
\]
For an invalid predicate code $d$, set $\beta(d,e)=\beta_\infty$.

Apply this search construction to the Rosser sentence to obtain $\beta_a$, so ZFC proves that $\rho_a$ holds if and only if $\mathsf U(\beta_a)$ does not halt. Use $\cF_{r,\beta_a}$ and, when $r\ge r_{\rm str}$, $\cG_{r,\beta_a}$. At $r=3$, the numerical output is the single $3$-graph $P_{\beta_a}$; denote it by $P_T$. The index transformations are total.

By \Cref{prop:numericalbridge} and the choice of $\beta_a$, density equality and $S_1,S_2,S_5$ are equivalent to $\rho_a$ in $T$, while $S_3,S_4$ are equivalent to its negation. A proof or refutation in $T$ of any of these assertions would yield a proof or refutation of $\rho_a$ in $T$. The baseline assertions are ZFC-provable, and the external truth values follow from the truth of $\rho_a$.
\end{proof}

\section{Effective spectra and finite value descriptions}\label{sec:spectrum}
We now allow a forbidden family to be specified by a membership algorithm or an enumeration, so it may be infinite. For $r\ge3$, such a family can have noncomputable density, although finite-family densities are computable. Finite initial subfamilies give successive upper bounds, leading to a characterization of the possible densities and their comparison complexities. We also study oracle variants. Subsection~\ref{sec:semialgebraic} examines two finite schemes for exact values.

\subsection{Effective spectra}
A \emph{computably enumerable} (\emph{c.e.}) family is supplied by an algorithm listing its members; a \emph{recursive} family by an algorithm that always halts and decides membership. Malformed outputs and codes of edgeless hypergraphs are discarded. A family is an \emph{antichain} if no member is isomorphic to a subgraph of another member. Write $\Pi_\infty^{(r)}$ for the density spectrum of arbitrary forbidden families.

A real is \emph{computable} if it has effective rational approximations of prescribed error, and \emph{right-c.e.} if it is the limit of a computable nonincreasing rational sequence~\cite{Weihrauch2000}. These nonincreasing approximations give upper bounds, but need not provide a computable bound on the remaining error. Write $\mathbb R_{\rm right\text{-}ce}$ for the right-c.e. reals.

\begin{theorem}\label{thm:spectrum}
For every $r\ge2$, with $\cF$ ranging over families of $r$-graphs,
\begin{equation}\label{eq:spectrum}
 \{\pi(\cF)\colon \cF\text{ is recursive}\}
 =\{\pi(\cF)\colon \cF\text{ is c.e.}\}
 =\Pi_\infty^{(r)}\cap\mathbb R_{\rm right\text{-}ce}.
\end{equation}
From an index for a computable nonincreasing rational sequence converging to $\alpha\in\Pi_\infty^{(r)}$, one can uniformly compute an index for a recursive antichain of density $\alpha$. For every $r\ge3$, some recursive antichain has noncomputable density.
\end{theorem}
\begin{proof}
For a c.e. family, truncate its enumeration after $s$ computation steps. The densities of these finite subfamilies form a nonincreasing sequence converging to the density of the full family by \eqref{eq:BS}. Approximate the density at stage $s$ from above within $2^{-s}$ and take running minima, giving a computable nonincreasing rational approximation.

Conversely, fix a computable nonincreasing rational sequence $\boldsymbol q=(q_n)$ converging to $\alpha\in\Pi_\infty^{(r)}$. We construct from this sequence alone a recursive forbidden family of density $\alpha$. Define
\begin{equation}\label{eq:recursivefamily}
 \cR_{\boldsymbol q}\coloneqq\{G\colon v(G)\ge r,\ r!|G|>q_{v(G)}v(G)^r\}.
\end{equation}
On input $G$, reject if $v(G)<r$. Otherwise compute $q_{v(G)}$ using the supplied approximation program, count the edges, and test the displayed rational inequality. The approximation program is total, so this membership test halts.

An $n$-vertex free $r$-graph is not itself in the family, so its density is at most $q_n n^r/(r!\binom{n}{r})\to\alpha$. This gives $\pi(\cR_{\boldsymbol q})\le\alpha$.

For the reverse inequality, choose any family of density $\alpha$ and take its homomorphic closure $\cB$. This family is used only to verify the lower bound; the construction of $\cR_{\boldsymbol q}$ reads only $\boldsymbol q$. If $H$ is $\cB$-free, then every subgraph $J\subseteq H$ with $v(J)\ge r$ has all its balanced blowups $\cB$-free, and hence
\[
 r!|J|\le\alpha v(J)^r\le q_{v(J)}v(J)^r.
\]
Thus $H$ avoids $\cR_{\boldsymbol q}$, proving $\pi(\cR_{\boldsymbol q})\ge\pi(\cB)=\alpha$.

Take the subgraph-minimal members of \eqref{eq:recursivefamily}. Test the $r$-graph and its finitely many proper subgraphs, including vertex deletions, to decide minimality. Every member of $\cR_{\boldsymbol q}$ contains a subgraph-minimal member of $\cR_{\boldsymbol q}$. Hence an $r$-graph avoids the minimal members exactly when it avoids the whole family. The minimal members therefore form a recursive antichain of density $\alpha$. The index transformation is total.

For a noncomputable example, fix a standard effective enumeration $(\psi_e)_{e\in\mathbb N}$ of the partial computable functions for the rest of this section. Let $K\coloneqq\{e\colon\psi_e(e)\text{ halts}\}$ be the diagonal halting set. The real $\xi\coloneqq\sum_{e\in K}4^{-(e+1)}$ is \emph{left-c.e.}, meaning that it has a computable nondecreasing rational approximation. The usual Specker-sequence argument~\cite{Specker1949} proves noncomputability. If $\xi$ were computable, we could determine membership in $K$ successively. After subtracting the already determined terms with index less than $e$, the remainder is at least $4^{-(e+1)}$ if $e\in K$, and at most $4^{-(e+1)}/3$ otherwise. Approximating this remainder to error less than $4^{-(e+1)}/6$ distinguishes the two cases, contradicting the undecidability of $K$.

To realize an affine transform of $\xi$ as a Tur\'an density, use the interval $[1-\delta_r,1]\subseteq\Pi_\infty^{(r)}$ supplied by Liu--Pikhurko~\cite[Theorem 1.1]{LiuPikhurkoIntervals2026} for $r\ge3$ and some $\delta_r>0$. Fix a rational $0<\epsilon_r<\delta_r$, nonuniformly in $r$. Then the right-c.e. noncomputable number $1-\epsilon_r\xi$ lies in this interval and has a recursive antichain presentation.
\end{proof}

\subsection{Comparisons and oracle variants}
For a c.e. family, a rational upper approximation can certify that its density is below a fixed rational threshold. Being above the threshold requires a positive gap to persist through every approximation stage. Thus the two comparisons have different quantifier forms.

A relation is \emph{$\Sigma^0_2$} if it has the form $\exists u\,\forall v\,R(e,u,v)$ for a decidable predicate $R$, and \emph{$\Pi^0_2$} if it has the form $\forall u\,\exists v\,R(e,u,v)$. More generally, for $j\ge1$, $\Sigma^0_j$ and $\Pi^0_j$ use $j$ alternating blocks of number quantifiers, beginning with $\exists$ and $\forall$, respectively, followed by a decidable predicate.

For fixed $r\ge3$, choose rationals $q,c>0$ with $[q-c,q+c]\subseteq\Pi_\infty^{(r)}\cap(0,1)$. Let $\cH_q$ be the fixed recursive family produced by the constant sequence $q_n=q$.
\begin{theorem}\label{thm:cecomparison}
For c.e.-family indices, the threshold relations have complexities
\[
 \begin{array}{c|c}
  \pi(\cF)<q&\Sigma^0_1\text{-complete}\\
  \pi(\cF)\ge q&\Pi^0_1\text{-complete}\\
  \pi(\cF)>q,\ \pi(\cF)\ne q&\Sigma^0_2\text{-complete}\\
  \pi(\cF)\le q,\ \pi(\cF)=q&\Pi^0_2\text{-complete}.
 \end{array}
\]
For two c.e. families, $>,<,\ne$ are $\Sigma^0_2$-complete and $\le,\ge,=$ are $\Pi^0_2$-complete. The second family may be fixed to $\cH_q$ for $>,\ne,=,\le$, and the first for $<,\ge$. Variable families in the reductions can be recursive antichains.
\end{theorem}
\begin{proof}
For the upper bounds, let $x_e$ be the density of the c.e. family indexed by $e$. The proof of \Cref{thm:spectrum} gives rational approximations $u_{e,s}$, uniformly computable in $e,s$ and nonincreasing in $s$, with limit $x_e$. Then
\[
\begin{aligned}
 x_e<q&\quad\text{if and only if}\quad\exists s\ u_{e,s}<q,\\
 x_e>q&\quad\text{if and only if}\quad\exists k\ \forall s\ u_{e,s}\ge q+2^{-k},\\
 x_e>x_f&\quad\text{if and only if}\quad\exists k,N\ \forall s\ge N\colon\ u_{e,s}\ge u_{f,s}+2^{-k}.
\end{aligned}
\]
The first two equivalences follow from $u_{e,s}\ge x_e$ and convergence to $x_e$. For the last equivalence, suppose $x_e>x_f$. Choose $k$ with $2^{-k}<(x_e-x_f)/2$, and then choose $N$ such that $u_{f,s}<x_f+(x_e-x_f)/2$ for all $s\ge N$. Since $u_{e,s}\ge x_e$, we obtain $u_{e,s}>u_{f,s}+2^{-k}$ for all $s\ge N$. Conversely, taking limits gives $x_e\ge x_f+2^{-k}>x_f$.

For comparisons between two families, taking the union of the two strict relations places $\ne$ in $\Sigma^0_2$, and taking complements places $\le,\ge,=$ in $\Pi^0_2$. For the fixed rational threshold $q$, the first equivalence also places $x_e\ge q$ in $\Pi^0_1$.

For first-level hardness, use a computable rational sequence equal to $q$ until $\psi_e(e)$ halts and to $q-c$ thereafter. Apply \Cref{thm:spectrum} to this sequence. We regard each resulting recursive antichain as a c.e. family by dovetailing its membership program on finite $r$-graph codes.

For second-level hardness, we encode whether $\psi_e$ halts on every input. Let
\[
 d_e=\sum_{\substack{x\\\psi_e(x)\text{ does not halt}}}2^{-x-1},\qquad
 d_{e,s}=1-\sum_{\substack{x<s\\\psi_e(x)\text{ has halted by }s}}
                 2^{-x-1}.
\]
The sets in the second sum increase to $\{x\colon\psi_e(x)\text{ halts}\}$. Hence $d_{e,s}$ is nonincreasing and $\lim_{s\to\infty}d_{e,s}=1-\sum_{\psi_e(x)\text{ halts}}2^{-x-1}=d_e$, since $\sum_{x\ge0}2^{-x-1}=1$. Thus $d_e=0$ exactly on the $\Pi^0_2$-complete totality set~\cite{Rogers1967}. Apply \Cref{thm:spectrum} to $q+cd_{e,s}$ to obtain a recursive antichain $\cF_e$ of density $q+cd_e$. All these densities are at least $q$, so comparison with $q$ by $=$ or $\le$ expresses totality, and comparison by $>$ or $\ne$ expresses nontotality. This proves the four second-level hardness claims.

Finally, comparing with $\cH_q$ or reversing the comparison gives the two-family statements.
\end{proof}

An \emph{oracle} $A\subseteq\mathbb N$ permits an algorithm to ask whether an integer belongs to $A$. The terms \emph{$A$-recursive} and \emph{$A$-c.e.} allow these queries in a membership algorithm or enumerator, and $\mathbb R^A_{\rm right\text{-}ce}$ denotes the reals with $A$-computable nonincreasing rational approximations. Let $\zero^{(0)}\coloneqq\varnothing$ and let $\zero^{(d+1)}$ be the halting problem for machines with oracle $\zero^{(d)}$. The constructions give the following spectrum and comparison statements with oracles.

\begin{theorem}\label{thm:oracles}
The following statements hold.
\begin{enumerate}
\item For every oracle $A\subseteq\mathbb N$ and every $r\ge2$, with $\cF$ ranging over families of $r$-graphs,
\[
 \{\pi(\cF)\colon \cF\text{ is }A\text{-recursive}\}
 =\{\pi(\cF)\colon \cF\text{ is }A\text{-c.e.}\}
 =\Pi_\infty^{(r)}\cap\mathbb R^A_{\rm right\text{-}ce}.
\]
From an index of an $A$-computable nonincreasing rational sequence converging to $\alpha\in\Pi_\infty^{(r)}$, one can compute an index of an $A$-recursive antichain of density $\alpha$. The transformation does not query $A$.
\item Fix $r\ge3$ and the rational interval used in \Cref{thm:cecomparison}. For $A=\zero^{(d)}$, $d\ge0$, consider indices of $A$-c.e. families. The threshold comparisons listed in the four rows of the table in \Cref{thm:cecomparison} have respective complexities $\Sigma^0_{d+1}$, $\Pi^0_{d+1}$, $\Sigma^0_{d+2}$ and $\Pi^0_{d+2}$, with completeness under ordinary many-one reductions and $A$-recursive antichain outputs. For two $A$-c.e. families, the relations $>,<,\ne$ are $\Sigma^0_{d+2}$-complete, and the relations $\le,\ge,=$ are $\Pi^0_{d+2}$-complete.
\end{enumerate}
\end{theorem}
\begin{proof}
For an $A$-c.e. family, run its enumeration for $s$ steps using $A$ and approximate the density of the resulting finite subfamily from above within $2^{-s}$, using \eqref{eq:approx}. Taking running minima gives an $A$-computable nonincreasing rational sequence. By \eqref{eq:BS}, its limit is the density of the full family.

Conversely, let $\boldsymbol q$ be an $A$-computable nonincreasing rational sequence with limit $\alpha\in\Pi_\infty^{(r)}$. The density bounds for $\cR_{\boldsymbol q}$ in the proof of \Cref{thm:spectrum} give $\pi(\cR_{\boldsymbol q})=\alpha$. On input $G$, the generated membership program tests whether $G$ is a subgraph-minimal member of $\cR_{\boldsymbol q}$. It calls the supplied $A$-computable approximation program for each required term of $\boldsymbol q$; edge counting and finite-subgraph enumeration do not use $A$. This gives an $A$-recursive antichain with the same free class and density. Producing its index only inserts the approximation program's code into this fixed procedure, without querying $A$.

For the upper bounds on comparison complexity, use the quantified formulas in the preceding proof with $A$-computable nonincreasing rational approximations. For first-level hardness, the approximation program starts at $q$ and drops to $q-c$ as soon as a bounded simulation finds that $\psi_e^A(e)$ halts. For second-level hardness, use the approximation $q+cd_{e,s}$ from the preceding proof with $\psi_e$ replaced by $\psi_e^A$. Both reductions produce program indices by inserting $e$ into these fixed approximation programs and applying the index transformation in part (i). Oracle queries occur when the generated programs run, so the reductions themselves are ordinary computable maps.

For $A=\zero^{(d)}$, Post's theorem~\cite{Post1944,Soare1987} identifies $\Sigma_i^{\zero^{(d)}}$ and $\Pi_i^{\zero^{(d)}}$ with the displayed arithmetical levels.
\end{proof}
\subsection{Two finite schemes for exact values}\label{sec:semialgebraic}
For finite forbidden families, \Cref{thm:finitecomparison} gives approximations to any prescribed accuracy. We now ask whether exact densities can be described by finite rules applied repeatedly within vertex classes or by finite polynomial optimization systems. The two schemes below yield only algebraic values.

For $m\ge1$, let $\Delta_m\coloneqq\{\x\in[0,1]^m\colon\sum_{i=1}^m x_i=1\}$ be the probability simplex. A \emph{pattern} $P=(m,E,R)$ consists of a set $E$ of $r$-multisets on $[m]$ and $R\subseteq[m]$. Partition the vertices into $V_1,\ldots,V_m$. For each $D\in E$, include all $r$-sets meeting each $V_i$ in $D(i)$ vertices. Repeat the construction inside $V_i$ for $i\in R$. The top-level density polynomial is $\lambda_E(\x)\coloneqq r!\sum_{D\in E}\prod_{i=1}^m\frac{x_i^{D(i)}}{D(i)!}$, where $\x$ gives the part proportions.

A bound $z$ on the densities within the recursive parts also bounds the total density for every choice of part proportions whenever $z$ belongs to
\[
 S_P\coloneqq\{z\in[0,1]\colon\lambda_E(\x)+z\sum_{i\in R}x_i^r\le z
                      \text{ for every }\x\in\Delta_m\}.
\]
The first scheme consists of the numbers $\min S_P$ for nonempty $S_P$, and finite maxima of such numbers. It includes the values of proper recursive-pattern constructions of Pikhurko~\cite{Pikhurko2014}.

Nonemptiness of $S_P$ already excludes an edge supported entirely on a part indexed by $R$, since evaluating at that part's unit vector would give $1+z\le z$. This is the relevant admissibility requirement preventing a top-level edge from also being counted inside a recursive part.

Fix $P$ with $S_P\ne\varnothing$. To show that finite constructions approach $\min S_P$, let $V_d$ be the supremum of limiting densities over all partition proportions when recursion stops after $d$ levels and the remaining recursive parts are left edgeless. Then $V_0=0$. For fixed top-level proportions $\x$, the top-level edges contribute $\lambda_E(\x)$, and the recursive parts contribute at most $V_d\sum_{i\in R}x_i^r$, with this bound approached by choosing near-optimal constructions in those parts. Thus
\[
 V_{d+1}=\max_{\x\in\Delta_m}
       \left(\lambda_E(\x)+V_d\sum_{i\in R}x_i^r\right).
\]
For fixed $d$, the limiting density is a polynomial in the proportions used in the finitely many partitions. By continuity, choose rational proportions for which this value is arbitrarily close to $V_d$. Taking the total order through multiples of the product of these denominators realizes every part size exactly. The densities of the resulting finite $r$-graphs then converge to the chosen polynomial value.

For each fixed $\x$, the coefficient $\sum_{i\in R}x_i^r$ of $V_d$ is nonnegative. Since $V_1\ge V_0=0$, induction in the recurrence gives $V_{d+1}\ge V_d$. For every $z\in S_P$, induction gives $V_d\le z$, so $V_d\to v\le\min S_P$. For each fixed $\x\in\Delta_m$, pass to the limit in $\lambda_E(\x)+V_d\sum_{i\in R}x_i^r\le V_{d+1}$. This gives $v\in S_P$, and therefore $v=\min S_P$.

For mixing constructions, Liu and Pikhurko~\cite[Lemma~3.4]{LiuPikhurko2025} proved that the maximum asymptotic density equals the largest of the single-pattern values. Thus the finite maxima in the first scheme also cover their mixing constructions.

The second scheme determines several values simultaneously, allowing the optimization for each coordinate to depend on the whole vector. A set is \emph{semialgebraic over $\mathbb Q$} if it is defined by a finite Boolean combination of polynomial equalities and inequalities with rational coefficients.

Fix $s\ge1$. For each coordinate $i\in[s]$, choose a nonempty finite index set $A_i$, nonempty compact parameter sets $K_{i,a}$ semialgebraic over $\mathbb Q$ for $a\in A_i$, and polynomials $p_{i,a}$ and $q_{i,a,j}$, $j\in[s]$, with rational coefficients. Define
\[
 (\mathcal B\boldsymbol z)_i\coloneqq
 \max_{a\in A_i}\max_{x\in K_{i,a}}
       \left(p_{i,a}(x)+\sum_jq_{i,a,j}(x)z_j\right).
\]
Assume $\mathcal B$ maps $[0,1]^s$ into itself and has a unique fixed point there. A designated coordinate of this fixed point is a value of the second scheme.

Real quantifier elimination rewrites quantified polynomial conditions as finite Boolean combinations of polynomial equalities and inequalities.
\begin{proposition}\label{prop:algebraicity}
Every value in either scheme is algebraic. For every sufficiently large $r$, some finite forbidden family has a Tur\'an density that admits neither kind of exact-value description.
\end{proposition}
\begin{proof}
The nonempty set $S_P$ is closed in $[0,1]$, so its minimum exists. Real quantifier elimination~\cite{Tarski1951}, in the form described in \cite[Chapter~14]{BasuPollackRoy2006}, shows that $S_P$ is semialgebraic over $\mathbb Q$. In one variable, a finite Boolean combination of rational polynomial sign conditions can change truth value only at a root of one of the nonzero polynomials. Every finite boundary point is therefore algebraic, including $\min S_P$. Finite maxima preserve algebraicity.

For the second scheme, a vector $\boldsymbol z\in[0,1]^s$ is fixed precisely when, for every coordinate $i\in[s]$,
\[
 \begin{aligned}
 &\forall a\in A_i\ \forall x\in K_{i,a},\quad
 p_{i,a}(x)+\sum_jq_{i,a,j}(x)z_j\le z_i,\\
 &\exists a\in A_i\ \exists x\in K_{i,a},\quad
 p_{i,a}(x)+\sum_jq_{i,a,j}(x)z_j=z_i.
 \end{aligned}
\]
The maxima are attained because the parameter sets are compact and the objective functions are continuous. Since each $A_i$ is finite, its universal and existential quantifiers expand into finite conjunctions and disjunctions. Quantifier elimination therefore shows that the fixed-point set is semialgebraic over $\mathbb Q$. By uniqueness, its projection to any coordinate axis is a singleton, again semialgebraic over $\mathbb Q$ by quantifier elimination. This point is a boundary point, so the preceding one-dimensional argument shows that the coordinate is algebraic.

Since both schemes produce only algebraic values, neither represents the finite-family transcendental densities of Li--Liu--Liu~\cite[Theorem 1.1]{LiLiuLiu2026}.
\end{proof}

\section*{Declaration on the use of AI}
The authors used generative AI tools to discuss proof organization, test local calculations, check proofs, and improve exposition.

\end{document}